\newcommand{\MA}{\mathfrak{A}}
\newcommand{\Tr}{\mbox{Tr}}
\newcommand{\MsL}{\mathscr{L}}
\newcommand{\Lve}{\lVert}
\newcommand{\Rve}{\rVert}

\newcommand\BH{\mathbb{H_\rho}}

\newcommand{\mU}{\mathcal{U}_{\mathfrak{A}}(K_1, K_2,K_3) }

\newcommand{\baray}{\begin{array}{rcl}}
\newcommand{\earay}{\end{array}}
\newcommand{\barray}{\begin{array}{rcl}}
\newcommand{\earray}{\end{array}}

\newcommand{\bE}{\mathbb{E}}

\newcommand{\X}{\mathbb{X}}
\newcommand{\DeltaA}{A}
\newcommand{\D}{\mathbb{D}}

\newcommand\comca[1]{{\color{blue} #1}} 
\newcommand\dela[1]{}

\newcommand{\bcase}{\begin{cases}}
\newcommand{\ecase}{\end{cases}}

\newcommand\Leb{\lambda }

\newcommand\del[1]{}

\del{

\newcommand\del[1]{}
}

\newcommand{\lk}{\left}
\newcommand{\lqq}{\lefteqn}
\newcommand{\rk}{\right}
\newcommand{\la}{{\langle}}
\newcommand{\ra}{{\rangle}}

\newcommand{\LL}{{\rm I \kern -0.2em L}}

\newcommand{\ep} {\varepsilon }

\newcommand{\be} {\begin{enumerate} }
\newcommand{\ee} {\end{enumerate} }
\newcommand{\BB} {X}

\newcommand{\CO}{{{ \mathcal O }}}
\newcommand{\CK}{{{ \mathcal K }}}

\newcommand{\CC}{{{ \mathcal C }}}

\newcommand{\SCT}{{{ \mathscr{T} }}}

\newcommand{\CH}{{{ \mathcal H }}}

\newcommand{\MsV}{{{\mathscr{V}}}}

\newcommand{\CM}{{{ \mathcal M }}}

\newcommand{\BF}{{{ \mathbb{F} }}}
\newcommand{\CF}{{{ \mathcal F }}}

\newcommand{\CN}{{{ \mathcal N }}}

\newcommand{\RR}{{\mathbb{R}}}

\newcommand{\NN}{\mathbb{N}} 

\newcommand{\PP}{{\mathbb{P}}}

\newcommand{\EE}{ \mathbb{E} }

\newcommand{\DEQS}{\begin{eqnarray*}}
\newcommand{\EEQS}{\end{eqnarray*}}
\newcommand{\DEQSZ}{\begin{eqnarray}}
\newcommand{\EEQSZ}{\end{eqnarray}}
\newcommand{\DEQ}{\begin{eqnarray}}
\newcommand{\EEQ}{\end{eqnarray}}

\newcommand{\Dcal} {{\mathcal D}}

\newcommand{\Fcal} {{\mathcal F}}

\newcommand{\Mcal} {{\mathcal M}}

\newcommand{\Ocal} {{\mathcal O}}

\newcommand{\Vcal} {{\mathcal V}}

\newcommand{\Xcal} {{\mathcal X}}

\newcommand{\Afrak} {{\mathfrak A}}

\newcommand{\R}{\mathbb{R}}

\renewcommand{\P}{\mathbb{P}}
\newcommand{\Eb}{\mathbb{E}}

\renewcommand{\v}{v}

\newcounter{gg1q1}
\newenvironment{listinitial}
{\begin{list} {(\roman{gg1q1})$\,$}
{\usecounter{gg1q1}
\setlength{\leftmargin}{0.4cm}
\setlength{\topsep}{0.0cm}
\setlength{\itemsep}{0.0cm}
\setlength{\parsep}{0.3cm}
\setlength{\itemindent}{0.9cm}
\setlength{\parskip}{0.2cm}}}
{\end{list}}

\newcounter{con}
\newcommand{\Con}{C_{\arabic{con}}}
\newcommand{\tCon}{\tilde C_{\arabic{con}}}
\newcommand{\Cdelta}{\delta_{\arabic{con}}}
\newcommand{\SC} {\stepcounter{con}}

\newcounter{coz}
\newcommand{\cz}{c_{\arabic{coz}}}

\newcommand{\Cdz}{\rho_{\arabic{coz}}}
\newcommand{\SCz} {\stepcounter{coz}}

\documentclass[11pt, leqno]{amsart}

\usepackage{xcolor}
	\definecolor{apricot}{rgb}{0.98, 0.81, 0.69}
	\definecolor{babyblue}{rgb}{0.54, 0.81, 0.94}

\usepackage{color}
\usepackage{amsmath}   
\usepackage{amsfonts,amssymb, color}
\usepackage{dsfont}
\usepackage{latexsym}
\input colordvi

\usepackage{mathrsfs}
\usepackage{enumitem}

\usepackage{enumitem}

\usepackage{amsfonts}

\theoremstyle{plain}

\begin{document}

\title[The stochastic Gierer-Meinhardt system]{The stochastic Gierer-Meinhardt system}

\author{Erika Hausenblas}
   \address{Department of Mathematics,
	Montanuniversit\"at Leoben,
	Austria.}
\email{erika.hausenblas@unileoben.ac.at}
\author[Akash A. Panda]{Akash Ashirbad Panda}
\address{Mathematisches Institut,
	Universit\"at T\"ubingen,
	D-72076 T\"ubingen,
	Germany.}
\email{panda@na.uni-tuebingen.de}

\begin{abstract}
The Gierer-Meinhardt system occurs in morphogenesis, where the development of an organism from a single cell is modelled. One of the steps in the development is the formation of spatial patterns of the cell structure, starting from an almost homogeneous cell distribution.
Turing proposed different activator-inhibitor systems with varying diffusion rates in his pioneering work, which could trigger the emergence of such cell structures.
Mathematically, one describes these activator-inhibitor systems as a coupled systems of reaction-diffusion equations with hugely different diffusion coefficients and highly nonlinear interaction. One famous example of these systems is the Gierer-Meinhardt system.
These systems usually are not of monotone type, such that one has to apply other techniques.

The purpose of this article is to study the stochastic reaction-diffusion Gierer–Meinhardt system with homogeneous Neumann boundary conditions on a one or two-dimensional bounded spatial domain. To be more precise, we perturb the original Gierer-Meinhardt system by an infinite-dimensional Wiener process and show under which conditions on the Wiener process and the system, a solution exists. In dimension one, we even show the pathwise uniqueness.
In dimension two, uniqueness is still an open question.

\end{abstract}

\date{\today}

\maketitle

\section*{Acknowledgements} This work was partially  supported by the Austrian Science Foundation (FWF) project number P28010 and P28819.

\smallskip

\pagestyle{myheadings} \markright{\today}\pagenumbering{arabic}

\theoremstyle{plain}

\numberwithin{equation}{section}



\newtheorem{theorem}{Theorem}[section]
\newtheorem{notation}{Notation}[section]
\newtheorem{claim}{Claim}[section]
\newtheorem{lemma}[theorem]{Lemma}
\newtheorem{corollary}[theorem]{Corollary}
\newtheorem{hypo}[theorem]{Hypothesis}
\newtheorem{example}[theorem]{Example}
\newtheorem{assumption}[theorem]{Assumption}
\newtheorem{tlemma}{Technical Lemma}[section]
\newtheorem{definition}[theorem]{Definition}
\newtheorem{remark}[theorem]{Remark}
\newtheorem{hypotheses}{H}
\newtheorem{proposition}[theorem]{Proposition}
\newtheorem{Notation}{Notation}
\renewcommand{\theNotation}{}

\renewcommand{\labelenumi}{\alph{enumi}.)}

\textbf{Keywords and phrases:} {Gierer-Meinhardt system, Pattern formation, Coupled system, Activator-Inhibitor system, Stochastic Partial Differential Equations, stochastic systems, Wiener process, Biomathematics}

\textbf{AMS subject classification (2020):} {Primary 60H15, 35A08, 35D30, 35G25, 35K55, 92C15;
Secondary 60G57, 35Q92, 35K45, 92C40, 35K57, 47H10, 92B05.}


\section{Introduction}
Pattern formation is a phenomenon based on the interaction of different components, possibly under the influence of their surroundings. Alan Turing, a cryptographer and a pioneer in computer science, developed algorithms to describe complex patterns using simple inputs and random fluctuation. In his seminal paper \cite{turing} in 1952, he proposed that the interaction between two biochemical substances with different diffusion rates have the capacity to generate biological patterns. In his mathematical framework, there is one activating protein (activator) that activates both itself and an inhibitory protein (inhibitor), which only inhibits the activator. He detected that a stable homogeneous pattern could become unstable if the inhibitor diffuses more rapidly than the activator. The interplay between the concentrations of these substances forms a pattern whose spatiotemporal evolution is governed by coupled reaction-diffusion systems (activator-inhibitor model). This phenomenon is called diffusion-driven instability or \emph{Turing instability}.
Thus, the most fundamental phenomenon in pattern-forming activator-inhibitor systems is that a slight deviation from spatial homogeneity has vital positive feedback leading to increase further. The presence of nonlinearities in the local dynamics, for example, due to the inhibitor concentration, saturates the Turing instability into a stable and spatially inhomogeneous pattern.
Pattern formation via diffusion-driven instabilities play an essential role in, e.g., biology, chemistry, physics, ecology, and population dynamics.

One well--established activator-inhibitor system was suggested by Gierer and Meinhardt in $1972$ to model the (re)generation phenomena in a hydra (see \cite{GM1}, and \cite{GM2}). This model, called the Gierer–Meinhardt model, describes the interaction of two biochemical substances, a slowly diffusing activator $U$ and a rapidly diffusing inhibitor $V$.
In particular, $U$ activates its production and the production of $V$. The inhibitor $V$ represses the production of $U$ and diffuses more rapidly than $U$.

\smallskip

\subsection{The Deterministic Gierer-Meinhardt Model}
For a domain, $\CO\subset \RR^d$, $d=1,2$,  the model introduced by Gierer and Meinhardt reads as
\DEQSZ\label{equ1s}
&&
\lk\{ \barray
 \dot {u}(t)  &=& r_u \, \Delta u(t) 
 + \kappa _u\,\frac {u^2(t)}{\v(t)} - \mu_u u(t), \quad  t> 0, 
\\ u(0) &=& u_0 ,  \phantom{\Big|}
\earray\rk.
\EEQSZ
and
\DEQSZ\label{eqv1s}
&&\lk\{ \barray
 \dot{\v}(t) &=& r_\v \,\Delta \v(t) + \kappa_\v \,u^2(t) - \mu_\v \v(t),\quad   t> 0,
\\
{\v}(0) &=& \v_0 ,\phantom{\Big|}
\earray\rk.
\EEQSZ
subjected to Neumann boundary conditions.
Here, the unknowns $u$ and $v$ stand for the concentrations $U$ and $V$, with a source distribution $\kappa_u$ and $\kappa_v$ respectively. Here $\kappa_u$ and $\kappa_v$ are positive constants. Furthermore,  the constants $r_u>0$ and $r_v>0$ are the diffusion coefficients, and the constants $\mu_u>0$ and $\mu_\v>0$ are given  decay rates.

\medskip

This system explains that starting from a uniform condition (i.e., a homogeneous distribution with no pattern), they could spontaneously self-organise their concentrations into a repetitive spatial pattern. The parameters in the system can be tuned in such a way, that some interesting phenomena such as Turing instability and peak steady states occur.

The Gierer-Meinhardt system has been studied extensively by many authors, both in the biological and physical communities; more recently also in mathematics, to elucidate its role in pattern formation. There are several works worth mentioning on the deterministic Gierer-Meinhardt system, e.g.,\ Masuda and Takahashi \cite{MT_1987} showed the existence
and boundedness of solutions. Gonpot {\em et al.} \cite{GCS_2008} studied the roles of diffusion and Turing Instability in the formation of spot and stripe patterns investigated by performing a nonlinear bifurcation analysis. Kelkel and Surulescu \cite{Kelkel+Surulescu_2009} proved the existence of a local weak solution for general initial conditions and parameters upon using an iterative approach. Chen {\em et. al.} \cite{CSWZ_2017} have carried out Bifurcation analysis, including theoretical and numerical analysis. Kavallaris and Suzuki \cite{KS_2018} focused on the derivation of blow-up results for the shadow Gierer-Meinhardt system. Recently, Wang {\em et. al.} \cite{GuGoWa_2020} investigated the stability of the equilibrium and the Hopf bifurcation of the Gierer-Meinhardt system of the Depletion type.
In addition, we refer for a more general background on the model to the book of Britton \cite{britton}, Meinhardt \cite{meinhardt}, the books of Murray \cite{murray1,murray2}, the book of Ghergu and Ruadulescu \cite{ghergu}, Perthame \cite{perthame1}, and of Wei and Winter \cite{wei1}.

\del{
Alens Gierer and Heins Meinhardt \cite{GM1}, and subsequently in \cite{GM2, GM3}, twenty years after the popular paper of Turing, proposed a model showing this kind of behaviour. Gierer and Meinhardt used this model, known as Gierer-Meinhardt system (abbreviated as GM-system), to explain the results of experiments on head regeneration and transplantation in the freshwater polyp hydra.
}


\subsection{Why Study Stochasticity?}

The deterministic model, i.e., \ the macroscopic system of equations, is derived from the microscopic behavior studying the limit behavior. From the microscopic perspective, one interprets the movements of the molecules as a result of microscopic irregular movement. Taking the limit and passing from the microscopic to the macroscopic equation, one neglects the fluctuations around the mean value.
Biological systems are frequently subject to noisy environments, inputs, and signalling. These stochastic perturbations are crucial when considering the ability of such models to reproduce results consistently. Murray \cite{murray1} addressed the possible impact of the noise, where he mentioned that to study the effect of stochasticity would be illuminating.

For a more realistic model, it is necessary to consider
features of the natural environment that are non-reproducible and, hence, modelled by random spatiotemporal forcing.
An appropriate mathematical approach to establish more realistic models is the incorporation of stochastic processes.
The randomness leads to a variate of new phenomena and may have a highly non--trivial impact on the behavior of the solution.
In Karig {\em et al.} \cite{karig}, the authors explore whether the stochastic extension leads to a broader range of parameters with Turing patterns by a genetically engineered synthetic bacterial population in which the signalling molecules form a stochastic activator{\textendash}inhibitor system.
Biancalani {\em et al.} \cite{noise2} studied the impact of noise on Turing pattern on several examples and showed that with a random noise, the range of parameter where
Turing patterns may appear, are enlarged.
Kolinichenko and Ryashko \cite{noise4}, respectively, Bashkirtseva {\em et al.} \cite{noise5} addresses multistability and noise-induced transitions between different states.
In Kolinichenko {\em et al.} \cite{noise5}, scenarios of noise-induced
pattern generation and stochastic transformations are
studied using numerical simulations and modality
analysis.
In summary, the stochastic term gives rise to a new type of behavior and, therefore, often leads to a more practical description of natural systems than their deterministic counterpart.

\subsection{The Stochastic Model}
Modelling randomness is done by adding a  random forcing term to the system  \eqref{equ1s} or the system \eqref{eqv1s} leading to a stochastic partial differential equation (SPDE).
 The internal noise  results in a multiplicative term;
since the white noise is an approximation of a continuously fluctuating noise with finite memory being much shorter than the dynamical timescales, the integration with respect to the noise is modelled by the Stratonovich integral.

\medskip

Let $\mathfrak{A}=(\Omega,\CF,(\CF_t)_{t\in[0,T]},\PP)$ be a filtered probability space. Let $W_j$, $j=1,2$, be two independent Wiener processes defined on $L^2(\CO)$ over the probability space $\mathfrak{A}$ with covariances $Q_1$ and $Q_2$ respectively, which will be explained in detail later.
We are now interested in the existence and uniqueness of the solution of the following stochastic Gierer-Meinhardt system
\DEQSZ\label{equ1ss}
 d {u}(t) & =&  \Big[r_u   \Delta u(t) + \kappa_u\, \frac {u^2(t)}{\v(t)} - \mu_u u(t) \Big]\, dt +\sigma_u u(t)\circ  dW_1(t),\quad t>0,
\EEQSZ
and
\DEQSZ\label{eqv1ss}
d{\v}(t) &=& \Big[ r_\v \Delta \v(t) + \kappa_v\, u^2(t) - \mu_\v \v(t) \Big]\,dt +\sigma_{\v} \v(t)\circ d W_2(t),\quad t>0,
\EEQSZ
again subjected to
Neumann boundary conditions (or if $\CO$ is a torus to periodic boundary conditions) and initial conditions $u(0)=u_0$ 
and $\v(0)=\v_0$. 
Besides,   $\sigma_u, \sigma_{\v} >0$ and  $\mu_u > 0$, $\mu_{\v} >0$. In dimension one we can show existence and uniqueness of the solution, in dimension two we can only show uniqueness.

\subsection{Novelty of the Paper}

The main difficulty in the treatment of \eqref{equ1s} and \eqref{eqv1s} is the lack of a variational structure. Jiang \cite{jiang} have shown the global existence of solutions to the deterministic Gierer-Meinhardt system. However, the method cannot be transferred to the stochastic case due to the randomness's peculiarity.
Due to the lack of variational structure, one can not apply standard methods to show the existence and uniqueness of solutions to stochastic partial differential equations. A way out is to use a stochastic Schauder-Tychonoff type Theorem which we present in Section \ref{schauder}.
In dimension one, we were able to show pathwise uniqueness, see Section \ref{sec-uniqueness}. By standard arguments based on the Yamada-Watanabe Theorem gives for 1D, the existence and uniqueness of a strong solution.

There are very few works known to the authors investigating the Gierer-Meinhardt system with full stochasticity.
 Li and Xu \cite{gierer1} consider the shadow Gierer-Meinhardt System with random initial data; Winter {\em et al.} \cite{gierer2} investigate the stochastic shadow Gierer-Meinhardt system.
Kelkel and Surulescu \cite{Kelkel+Surulescu_2010} prove locally in time a pathwise unique mild solution via an iterative method, which is done for positive saturation constants (that is biologically expedient); whereas our paper shows the existence of a martingale solution (global-in-time) without the positive saturation restriction (i.e. for nonnegative constants).

\bigskip

The article is structured as follows.
In Section \ref{mainresult}, we list the hypotheses and the main result is presented, i.e., \ the existence of a martingale solution to the stochastic counterpart of the system \eqref{equ1ss} and \eqref{eqv1ss}.
In Section \ref{proof} the actual proof is stated, in Section \ref{technics}, we prove several technical propositions necessary for Section \ref{proof}. The pathwise uniqueness in dimension one has been proved in Section \ref{sec-uniqueness}. In Section \ref{schauder}, we present a stochastic version of a Schauder-Tychonoff type Theorem which is used to prove the main result.

\smallskip

\noindent
Let us introduce few notations of functional spaces which will be used throughout the paper.

\del{
Since the white noise is an approximation of a continuously fluctuating noise with finite memory being much shorter than the dynamical timescales, the representation of the stochastic integral as a Stratonovich stochastic integral is appropriate.
In order to show the existence of a solution of the original system,
the linear parts have to be incorporated, which can be done by modifying the proof given here. The structure of the proof will remain even adding a gradient term in the equation.
}



\begin{notation}
For a Banach space $E$ and $0\le c<d<\infty$, let $C^{\delta}_b(c,d;E)$ denote a set of all continuous and bounded functions $u:[c,d]\to E$ such that
$$
\| u\|_{C_b^{\delta}(c,d;E)} :=\sup_{c\le t\le d} |u(t)|_E +\sup_{c\le s,t\le d\atop t\not= s} \frac{ |u(t)-u(s)|_{E}}{|t-s|^\delta},
$$
is finite. The space $C_b^{\delta}(c,d;E)$  endowed with the norm $\| \cdot\|_{C_b^{\delta}(c,d;E)}$ is a Banach space.
\end{notation}

\begin{notation}\label{notation2}
For $1< p < \infty,$  let $W^1_p(\CO)$ be the Sobolev space defined by (compare \cite[p.\ 263]{Brezis})
\DEQS
\lqq{W^{1}_p(\CO)
:=\lk\{ u\in L^p(\CO)\mid \exists g_1,\cdots,g_d\in L^p(\CO)\mbox{ such that }\phantom{\bigg|}\rk.}
\\
&&{}\lk. \qquad \quad \int_\CO u(x)\frac{\partial \phi(x)}{\partial x_i} \, dx =-\int_\CO g_i(x)\phi(x)\, dx \quad \forall \phi\in C^\infty_c(\Omega), \forall i=1,\ldots ,d \rk\}
\EEQS
equipped with norm
$$
|u|_{W^1_p}:=|u|_{L^p}+\sum_{j=1}^d\lk|\frac{\partial u}{\partial x_j}\rk|_{L^p},\quad u\in W^1_p(\CO).
$$
Given an integer $m\ge2$ and a real number $1\le p<\infty$, we define by induction the space
\DEQS
W^{m}_p(\CO) :=\lk\{ u\in W^{m-1}_p(\CO)\mid D u\in W^{m-1}_p(\CO) \rk\}
\EEQS
equipped with norm
$$
|u|_{W^m_p}:=|u|_{L^p}+\sum_{|\alpha| \leq m} \lk|D^\alpha u\rk|_{L^p},\quad u\in W^m_p(\CO).
$$
Let $H_2^m(\CO):=W^m_2(\CO)$, and for $\rho\in(0,1)$ let $H_2^\rho (\CO)$ be the real interpolation  space given by $H^\rho _2(\CO):=(L^2(\CO),H^1_2(\CO))_{\rho ,2}$.
In addition, let $H^{-1}_2(\CO)$ be the dual space of $H^1_2(\CO)$ and for $\rho \in(0,1)$ let $H^{-\rho }_2(\CO)$ be the  real interpolation  space given by $H^{-\rho }_2(\CO):=(L^2(\CO),H^{-1}_2(\CO))_{1-\rho ,2}$.
Note, by Theorem 3.7.1 \cite{bergh}, $H^{-\rho }_2(\CO)$ is dual to $H^\rho _2(\CO)$, $\rho \in(0,1)$. Furthermore, we have  $(H^{-\rho }_2(\CO),H^{\rho }_2(\CO))_{\frac 12,2}=L^2(\CO)$ and $(H^\alpha_2(\CO),H^\beta_2(\CO))_{\rho,2}=H^\theta_2(\CO)$ for $\theta=\alpha(1-\rho)+\beta\rho$, $\rho\in(0,1)$ and $|\alpha|,|\beta|\le 1$.
\end{notation}

\smallskip

\section{Hypotheses and the Main Result}\label{mainresult}

Let us denote $\mathfrak{A}=(\Omega,\CF, \mathbb{F},\PP)$ be a complete probability space with the filtration $\mathbb{F} = \{{{\mathcal{F}}}_t:t\in [0,T]\}$ satisfying the usual conditions i.e., $\PP$ is complete on $(\Omega, \CF)$, for each $t \geq 0, \CF_t$ contains all $(\CF, \PP)$-null sets, and the filtration $\BF$ is right-continuous. Let in case $d=1$ the domain $\CO$ be an interval and in case $d=2$, let $\CO\subset \RR^2$ be a bounded domain with $C^\infty$  boundary.
Let $W_1$ and $W_2$ be two independent Wiener processes  in $\CH:=L^2(\CO)$ defined over the probability space $\mathfrak{A}$.
Before introducing the  hypothesis on the noise let us introduce the unbounded operator $A$ given by the Laplace operator $-\Delta$ in $L^2(\CO)$ with Neumann boundary conditions. In particular, let
\DEQS
\lk\{ \begin{array}{rcl}
D(A) &=&\{ u\in H^2_2(\CO): \frac {\partial}{\partial n}u(x)=0,\, x\in\partial \CO\},
\\
Au&=&-\Delta , \quad u\in D(A).
\end{array}\rk.
\EEQS
Here, $n$ denotes the outward normal vector of the boundary $\partial \CO$.{
\begin{hypo}\label{wiener}
We assume that $W_1$ and $W_2$ are two independent Wiener processes defined on  $\CH:=L^2(\CO)$. In particular, we have
$$
W_j(t):=\sum_{k \in\NN} \beta^j_k(t)\,(1+\lambda_k)^{-\gamma_j/2}\, e_k,\quad t\ge 0, \,j=1, 2.
$$
where  $\gamma_j>d$, $j=1,2$, $\{\beta^j_k:k\in\NN\} $, $j=1,2$, are two independent families  of mutually independent real-valued Brownian motions over $\MA$,  $\{e_k:k\in\NN\}$ are the eigenfunctions of $A$ and $\{\lambda_k:k\in\NN\}$ with $(0 = \lambda_1 < \lambda_2 < \cdots)$ are the corresponding eigenvalues.

In the case of one dimension, we have  $\lambda_k = 4 \pi^2 k^2$, $k\in\NN$, and if $\CO$ is a square domain $(i.e. \ \CO = \{ (x,y) \in \mathbb{R}^2 : 0<x<a, 0<y<b\})$, then we have
$$
\lambda_k = \lambda_{l, m}= \lk(\frac{l \pi}{a}\rk)^2 + \lk(\frac{\phantom{l}\!\!m \pi}{b}\rk)^2, \quad l,m\in\NN .
$$

\del{Furthermore, let us define the mapping $\Sigma(z)$ that depends linearly on $z$, i.e.\ 
$$
\Sigma (z)h=  z \,(-\Delta)^{-\gamma_j/2}h,\quad h\in \CH, \ u \in \CH,
$$}

\end{hypo}
}
\noindent
For later on, let us define
the covariances $Q_j$ by $Q_j := (\mbox{Id} + A)^{-\gamma_j}$, $j=1,2$. Next, we define the notion of solution to the system \eqref{equ1ss}--\eqref{eqv1ss}.
The notion of solution will depend on a parameter $\rho$, later in Hypothesis \ref{init} we will see that
$\rho$ belongs to the interval $(1, \frac 65)$ if $d=2$ and $[1,\frac 65)$ if $d=1$.

\del{\begin{remark}
If $\CO=[0,1]^d$ is a rectangular domain, then the condition $\gamma_j >\frac 32 d, \,j=1,2,$ can be relaxed to $\gamma_j >\frac d2, \,j=1, 2.$ $($see \cite[p.\ 7]{BDPR2016}$)$.
\end{remark}}

\begin{definition}
We say a pair of progressively measurable processes $(u, \v)$ a strong solution to the system \eqref{equ1ss}--\eqref{eqv1ss} on $[0,T]$, $T>0$, for initial data $(u_0,\v_0)$ if for $\rho>0$ we have $\PP$-a.s.
{\[ u\in C_b([0,T]; H_2^{1-\rho}(\CO)))\cap L^{2}([0,T];H^{2-\rho}_2(\CO)),\]
and}
\[\v\in C_b([0,T];L^2(\CO)))\cap L^2([0,T];H_2^{1}(\CO)), \]
$u$ and $\v$ are $(\Fcal_t)_{t\ge 0}$-adapted, and satisfies for all $t\in [0,T]$, $\P$-a.s.
\begin{align}\label{equ1i}
u(t)=u_0+\int_0^t \Big[r_u   \Delta u(s) + \kappa_u\, \frac {u^2(s)}{\v(s)} -  \mu_u u(s) \Big] \,ds+\sigma_u\int_0^t u(s) \circ dW_1(s),
\end{align}
and
\begin{align}\label{eqv1i}
\v(t)=\v_0+\int_0^t \Big[r_\v \Delta \v(s) + \kappa_v \,u^2(s) -  \mu_\v  \v(s) \Big]\,ds+ \sigma_{\v}\int_0^t \v(s) \circ dW_2(s).
\end{align}
\end{definition}

Note that $r_u, r_v, \kappa_u, \kappa_v, \mu_u, \mu_v, \sigma_{u}, \sigma_{\v} >0$. As mentioned before, in the proof of the main result, we are using compactness arguments, which causes the loss of the original probability space. This means the solution will only be a weak solution in the probabilistic sense.

\begin{definition}\label{Def:mart-sol}
A {\sl  martingale solution}   to the problem
\eqref{equ1ss}--\eqref{eqv1ss} is a system
\begin{equation*}
\left(\Omega ,{{\mathcal{F}}},\mathbb{P},{\mathbb{F}},
(W_1,W_2), (u,\v)\right)
\end{equation*}
such that
\begin{itemize}[label=$\circ$]
\item  $\mathfrak{A}:=(\Omega ,{{\mathcal{F}}},{\mathbb{F}},\mathbb{P})$ is a complete filtered
probability space with a filtration ${\mathbb{F}}=\{{{\mathcal{F}}}_t:t\in [0,T]\}$ satisfying the usual condition;
\item $W_1$ and $W_2$ be two independent Wiener processes in $L^2(\CO)$, defined over the probability space $\mathfrak{A}$ with covariances $Q_1$ and $Q_2$; 
\item $u:[0,T]\times \Omega \to H^{1-\rho}_2(\CO)$ for $\rho>0$ and $\v:[0,T]\times \Omega \to L^2(\CO)$  are two  ${\mathbb{F}}$-progressively
measurable 
 process such that the couple $(u,\v)$ is a solution to the system \eqref{equ1ss}--\eqref{eqv1ss} over the probability space $\mathfrak{A}$.
 \end{itemize}
\end{definition}

\noindent
The range of $\rho$ will be specified in the next hypothesis. In this section, we list all the hypotheses on the  initial data and on the Wiener process.
%
%
%
%
\begin{hypo}
\label{init}
We assume that the initial data $(u_0,v_0)$ satisfies
\begin{listinitial}
\item $u_0\ge 0$ and $v_0 >0$ on $\CO$;
\item
there exists a number with
$$
\rho\in\bcase (1, \frac 65)&\mbox{ if }d=2,
\\ [1,\frac 65)&\mbox{ if }d=1,
\ecase
$$
such that
$ \EE |u_0|_{H_2^{1-\rho}}^2 < \infty$;
\item suppose that $ \EE|v_0|_{L^2}<\infty$; 
\item $\EE |\xi_0|^{p}_{L^{p}}<\infty$ for $p=31/7$ and  $\EE \big( \int_\CO\ln(\xi_0(x))\, dx \big)<\infty$, where $\xi_0=v_0^{-1}$.
\end{listinitial}
\end{hypo}




\begin{remark}
The constant $31/7$ and $\rho$ are mainly given by Lemma \ref{con-p1}. 
\end{remark}
\noindent
We assume that the noise is of linear type. With these hypotheses, we  are now ready to state the main result of our paper. 
\begin{theorem}\label{mainresulttheo}
Under the Hypotheses \ref{wiener} and \ref{init},
there exists a martingale solution to system \eqref{equ1ss}--\eqref{eqv1ss} in the sense of Definition \ref{Def:mart-sol}, satisfying the following properties
\begin{enumerate}[label=(\roman*)]
\item
The solutions $u(t, x) \geq 0$ and $\v(t, x) > 0,$ $\mathbb{P} \otimes \mbox{Leb}$ a.s.;
\item We have $\PP$--a.s. $u\in C_b(0,T;H^{1-\rho}_2(\CO))$ and $v \in C_b(0,T;L^2(\CO))$;
\item
There exist   constants $C_1,C_2>0$ such that
   \DEQS 
\EE  \sup_{t\in[0, T]} |u(t)|_{H^{1-\rho}_2}^{2} +2 \EE \int_0^T |\nabla u(t)|^2_{H^{1-\rho}_2}\,dt   \le  C_1(1+\EE|u_0|_{H^{1-\rho}_2}^2),
\EEQS
and
   \DEQS 
\EE  \sup_{t\in[0, T]} |v(t)|_{L^{2}} +2 \EE \Big(\int_0^T |\nabla v(t)|^2_{L^2}\,dt \Big)^\frac 12 \le  C_2(1+\EE|v_0|_{L^2}).
\EEQS
  \item
There exists a  constant $C_3>0$ such that we have for $\xi=\v^{-1}$ and  $p=31/7$
   \DEQS 
\EE \sup_{t\in[0, T]} |\xi(t)|_{L^{p}}^{p}
+ p(p+1) \EE \int_0^T \int_{\CO} \xi^{p+2}(t,x) \,(\nabla \v)^2(t,x) \,dx \,dt \le  {C}_3(1+ \bE |\xi_0|_{L^p}^p).
\EEQS
\end{enumerate}
\end{theorem}

The proof of this main result involves many steps and several tedious calculations. For the convenience of the reader, we assign a different section for the proof, which will constitute many lemmas and propositions. For the proof, we will use the stochastic version of the Schauder-Tychonoff type Theorem, which is stated and proved in Section \ref{schauder}. In Section \ref{proof} the actual proof of the main theorem is stated, and in Section \ref{technics}, the propositions used in the proof are postponed. The statement and the proof of the pathwise uniqueness in one dimension is obtained in Section \ref{sec-uniqueness}.

\section{Proof of the Main Result}\label{proof}

\subsection{Technical Preliminaries}

As mentioned before, the proof is an application of the stochastic Schauder-Tychonoff fixed point theorem and consists of several steps.
\begin{itemize}[label=$\circ$]
\item
In the first step, we specify the underlying Banach spaces;
\item
In the second step, we construct the operator $\mathscr{T}$;
\item
In the third step, we formulate the main claims and show that the assumptions of the stochastic Schauder-Tychonoff Theorem are satisfied;
\item
In the fourth step, we conclude the proof of the
main result.
\end{itemize}
To keep the proof itself simple, the proof will use many technical propositions, which are the content of Section \ref{technics}.
Before starting with the actual proof, we will set up the notation we will use.

\medskip

First, we are considering a linear noise, that is (see \cite[Example 2.1.2]{Barbu+DaPrato+Rockner_2008}).
%
The enumeration is chosen in increasing order by counting the multiplicity.
Here,
one can give
the following estimate of the asymptotic behavior of the eigenvalues: there exist two numbers $c,C>0$ such that
we have (e.g. see \cite[Subsection 1.2.1]{BDPR2016} and references therein)
\DEQSZ\label{EVsup}
c k^\frac 2d \le \lambda_k\le C k^\frac 2d, \quad  k\in\NN.
\EEQSZ
In addition, there exists some constant $c>0$ such that  (see \cite[p.\ 7]{Barbu+DaPrato+Rockner_2008})
\DEQSZ\label{EFsup}
\sup_{x\in\CO} |e_k(x)|\le c\,\lambda_k^\frac {d-1}2,\quad k\in\NN.
\EEQSZ

\medskip

A drawback of the Stratonovich stochastic integral is that it is not a martingale, and therefore, one does not have the Burkholder-Davis-Gundy inequality at its disposal. Hence, we changed within the proof from the Stratonovich stochastic integral to the It\^ o stochastic integral.
For a detailed discussion, we refer to the book of Duan and Wang \cite{duan} or to the original work of Stratonovich \cite{stratonovich,stratonovich1}.
The step from the Stratonovich integral to the It\^o integral can be done by adding a correction term, i.e., by adding drift to the It\^o form.
In this way, we will end up with a similar equation.
In particular, we will consider the following system
\DEQSZ \label{equation_uv_strat}
\lk\{ \barray
 d {u} (t)& =& \Big[ r_u\Delta   u(t)+\kappa_u \frac {u^2(t)}{v(t)} -\mu_u u(t)\Big]\, dt + \sigma_u u(t) \,dW_1(t),\quad t>0,
  \phantom{\Big|}
\\
d{v}(t) & = &\Big[r_v \Delta v(t)  +\kappa_v u^2(t) -\mu_v v(t)\Big]\, dt +\sigma_v v(t) \,dW_2(t) ,\quad t> 0,
 \phantom{\Big|}
\earray\rk.
\EEQSZ
where the scalar $\mu_u$ and $\mu_v$ are replaced by the operators
 $\Upsilon_u:=\mu_u \mbox{Id}-\sigma_u (\mbox{Id}+A)^{-\gamma_1}$ and  $\Upsilon_v:=\mu_v\mbox{Id}-\sigma_v (\mbox{Id}+A)^{-\gamma_2}$, respectively.
Observe, since $\gamma_1,\gamma_2\ge 0$,  we know that $\Upsilon_u$ and $\Upsilon_v$ are bounded operator on $L^p(\CO)$, $1\le p<\infty$. Since $\mu_u \ge 0$ and $u$ is non-negative,  $\mu_u u$ is non-negative. Since $\sigma_u>0$, we are only concern about the term $(\mbox{Id}+A)^{-\gamma_1} u$. Let $u$ be non--negative. If $z:=(\mbox{Id}+A)^{-\gamma_1} u$, then  $(\mbox{Id}+A)^{\gamma_1} z = u \geq 0$ in $\mathcal{O}$ and $\frac{\partial z}{\partial \nu} = 0$ on $\partial \mathcal{O}$. By the maximum principle, we obtain $z \geq 0$ in $\mathcal{O}$, i.e.\ $\Upsilon_u \geq 0$ in $\mathcal{O}$.  Similarly, one can show if $v$ is positive, the $\Upsilon_v v> 0$ in $\mathcal{O}$.


\medskip

Now, we can start with the actual proof of Theorem \ref{mainresulttheo}.

\begin{proof}[Proof of Theorem \ref{mainresulttheo}:]
$\phantom{mmmm}$

\paragraph{{\bf Step I}}  {\bf Definition of the functional spaces:}

Our main aim is to construct an integral operator, denoted by $\mathscr{T}$, whose fixed point is the solution of the coupled system \eqref{equation_uv_strat}. In this step, we define the required functional spaces for the operator $\mathscr{T}$ to act upon. Let the probability space $\mathfrak{A}=(\Omega,\CF,\mathbb{F},\PP)$ be given and let $W_1$ and $W_2$ be two independent cylindrical Wiener processes in $\CH=L^2(\CO)$ defined over $\mathfrak{A}$. Let $W=(W_1,W_2)$, $H=\CH\times \CH$ and let $Q$ be a covariance operator
defined by
$$
Q=\lk(\begin{matrix} (\mbox{Id}+A)^{-\gamma_1} & 0 \\0 &  (\mbox{Id}+A)^{ -\gamma_2}\end{matrix}\rk).
$$
\del{Let us put \comca{(to be clarified later)}
\begin{align}
\tilde{K} = ,
\end{align}
where is the constant appearing in the Proposition \ref{}.
\\
Let $K_1$ be so large such that
\begin{align}
 \leq K_1,
\end{align}
where is the constant appearing in the Proposition \ref{}.\\
Let $K_2$ be so large such that
\begin{align}
 \leq K_2,
\end{align}
where are the constants appearing in the Proposition \ref{}.}

\noindent
Let us define for $\rho \in (1, \frac 65)$ the space
\begin{equation*}
\begin{split}
\mathcal{M}_{\mathfrak{A}} := \mathcal{M}_{\mathfrak{A}}(0, T) := \bigg\{ &(\chi, \eta): \Omega \times [0, T] \times \CO \to \big( \mathbb{R}^{+}_{0} \times \mathbb{R}^{+}_{0}\big) :
\\ & \ (\chi, \eta) \mbox{ is } \mathbb{F}\mbox{-progressively measurable and }
\\ & \ \mathbb{E} \sup_{s\in[0, T]}|\chi(s)|^2_{H^{1-\rho}_2} + \mathbb{E} \sup_{s\in[0, T]}|\eta(s)|_{L^2} < \infty \bigg\},
\end{split}
\end{equation*}
equipped with the semi-norm
$$|(\chi, \eta)|_{\mathcal{M}_{\mathfrak{A}}} := \Big( \mathbb{E} \sup_{s\in[0, T]}|\chi(s)|^2_{H^{1-\rho}_2} \Big)^{1/2} + \Big( \mathbb{E} \sup_{s\in[0, T]}|\eta(s)|_{L^2} \Big), \quad \mbox{for} \ (\chi, \eta) \in \mathcal{M}_{\mathfrak{A}}.
$$
Let $K_1,K_2$ and $K_3$ be fixed constant.
Let us fix
the set
\begin{equation*}
\begin{split}
\mathcal{U}_{\mathfrak{A}}(K_1, K_2,K_3) := \bigg\{ &(\chi, \eta) \in \mathcal{M}_{\mathfrak{A}} : \chi(t, x) \geq 0 \mbox{ and } \eta(t, x) >0 \,\mbox{ Leb} \otimes \mathbb{P}-a.s.
\\ & \ \forall (t, x) \in [0, T] \times \mathcal{O}:\, \mathbb{E}\, \mathscr{L}_1(\chi, \eta) \leq K_1
\\ & \ \EE\,\MsL _2(\chi,\eta )
 \leq K_2 , \quad   \mbox{ and }\quad   \sup_{0\le t\le T} \EE\,\MsL _3(\chi(t),\eta(t))\le K_3\bigg\},
\end{split}
\end{equation*}
where the Lyapunov functionals $\mathscr{L}_1(\chi, \eta)$, $\mathscr{L}_2(\chi, \eta)$, and $\mathscr{L}_3(\chi, \eta)$ are  defined by
\DEQS
\mathscr{L}_1(\chi, \eta) &:= &  \sup_{s \in [0, T]} |\chi(s)|_{L^2}^2  + \Big(\int_0^ T |\nabla \chi(s)|^2_{L^2}\, ds\Big)+ \sup_{s \in [0, T]} |\xi(s)|_{L^p}^p  ,
\\
\mathscr{L}_2(\chi, \eta) &:=  &
\Big( \int_0^T \int_\CO \chi^2(s,x) \xi(s,x)\, dx \, ds\Big)^2  +\Big( \int_0^T \int_\CO \xi^2(s,x)\chi^2(s,x)\,dx\, ds\Big) ,
\EEQS
and
\DEQS
\mathscr{L}_3(\chi, \eta)(t) &:= &|\xi(t)|_{L^p}^p + 
|\xi(t)|_{L^1}
+
\Big(\int_\CO \ln(\xi(t,x)\, dx\Big)^2.
\EEQS
Here, the process $\xi$ is given by $\xi:= v^{-1}$, where $v$ depends on $\chi$ and solves
\begin{equation*} 
\lk\{\begin{split}
 d{\v}(t) &= \big[ r_\v \Delta \v(t) + \kappa_v\, \chi^2(t) -  \Upsilon_v \v(t)  \big] \,dt +\sigma_{\v} \v(t) \,dW_2(t)  ,\quad   t> 0,
\\
{\v}(0) &= \v_0 .\phantom{\Big|}
\end{split}\rk.
\end{equation*}

\medskip

\paragraph{{\bf Step II}}  {\bf Construction of the operator $\mathscr{T}$:}
In this step we construct an integral operator whose fixed point is the solution of the system \eqref{equ1ss}-\eqref{eqv1ss}.
Let us fix $K_1, K_2, K_3,$ to be some positive real numbers. First, let us define the following operator:
\DEQS
\mathscr{T}:= \mathscr{T}_{\mathfrak{A}} : \mathcal{U}_{\mathfrak{A}}(K_1, K_2,K_3)& \longrightarrow& \mathcal{M}_{\mathfrak{A}}(0, T);
\\ (\chi,\eta)&\mapsto& \mathscr{T}[(\chi,\eta)]:=(u,v),
\EEQS
where 
the pair $(u, \v)$ solves the following system
\begin{equation}\label{eqv1sto}
\lk\{\begin{split}
 d{\v}(t) &= \big[ r_\v \Delta \v(t) + \kappa_v\, \chi^2(t) -  \Upsilon_v \v(t)  \big] \,dt +\sigma_{\v} \v(t) \,dW_2(t)  ,\quad   t> 0,
\\
{\v}(0) &= \v_0 ,\phantom{\Big|}
\end{split}\rk.
\end{equation}
and
\begin{equation}\label{eqv2sto}
\lk\{\begin{split}
 du(t) &= \big[r_u \Delta u(t) + \kappa_u \frac {\chi^2(t)}{\v(t)} - \Upsilon_u  u(t) \big]\,dt + \sigma_u u(t)\, dW_1(t),\quad   t> 0,
\\
u(0) &= u_0.\phantom{\Big|}
\end{split}\rk.
\end{equation}
%
\begin{remark}
As one can see, $\eta$ will not be used in the definition of the operator. However, we keep $\eta$, since then the setting fits to the setting of the Schauder-Tychonoff-type Theorem. 
\end{remark}
The operator $\mathscr{T}$ is well-defined. In fact, from Theorem \ref{theou1}, stated in  the next section, we infer that for given $(\chi,\eta) \in \mathcal{U}_{\mathfrak{A}}(K_1, K_2,K_3),$ the existence of a non--negative unique solution $\v$ to \eqref{eqv1sto} such that we have
\DEQS
\mathbb{E}  \sup_{0 \leq s \leq T} \big|\v(s)\big|_{L^2} + \mathbb{E}\Big( \int_0^T \big|\nabla \v(s) \big|_{L^2}^2 \,ds\Big)^\frac12  < \infty.
\EEQS
From Theorem \ref{sol_w} we infer that there exists a unique solution $u$ to \eqref{eqv2sto} such that
\DEQS
\mathbb{E} \sup_{0 \leq s \leq T} \big|u(s)\big|^2_{H^{1-\rho}_2} + \mathbb{E} \int_0^T \big|u(s) \big|_{H^{2-\rho}_2}^2 \,ds < \infty.
\EEQS

\medskip

\paragraph{{\bf Step III}} {\bf Verification of the assumptions of the Schauder-Tychonoff type Theorem:}
In this step, we verify the assumptions of the Schauder-Tychonoff type Theorem \ref{schauder} by
proving the  three Claims: \ref{firstclaim}, Claim \ref{secondclaim}, and Claim \ref{thirdclaim}.
\begin{claim}\label{firstclaim}
 There exist $K_1, K_2 >0$, and $K_3>0$  such that $\mathscr{T}$ maps $\mathcal{U}_{\mathfrak{A}}(K_1, K_2,K_3)$ into $\mathcal{U}_{\mathfrak{A}}(K_1, K_2,K_3)$.
\end{claim}


\begin{proof}[Proof of Claim \ref{firstclaim}:]
First, we will start to show the existence of numbers $K_2,K_3>0$ (independent of $K_1$) such that, given $(\chi,\eta)\in\mU$, then
\DEQSZ\label{toshowk1}
 \EE\,\MsL _2(u, \v)\leq K_2,\quad \mbox{and}\quad  \sup_{0\le t\le T}\,  \EE\,\MsL _3(u(t), \v(t))\leq K_3.
\EEQSZ
\SCz
In fact, by Proposition \ref{pr-xi_est} we know that there exist constants $\cz>0$ and $\Cdz>0$ such that
\eqref{toshowk1}  is satisfied if
\DEQSZ
\label{condk2} K_2\ge \cz e^{\Cdz T}\big(1+\EE|\xi_0|_{L^2}^2 \big),
\EEQSZ
and 
\DEQSZ
\label{condk3} K_3\ge \cz e^{\Cdz T}\Big(1+\EE|\xi_0|_{L^p}^p+\EE|\xi_0|_{L^1}+\EE\Big(\int_\CO \ln(\xi_0(x))\, dx\Big) ^2\Big).
\EEQSZ
%
It remains to show that one can find a number $K_1>0$ such that if $(\chi,\eta)\in\mU$ then
\DEQSZ\label{toshowk2}
\mathbb{E} |\mathscr{L}_1(\chi, \eta)| \leq K_1 .
 \EEQSZ
\SCz
In Proposition \ref{prop74} 
 it is shown that 
 there exist constants $\cz>0$ and $\Cdz>0$ such that
we have for any $T>0$
\DEQSZ \label{u^2Ito-8}  \nonumber
\lqq{ \EE \sup_{s \in [0, T]} |u(s)|_{L^2}^2  +\EE \Big(\int_0^ T |\nabla u(s)|^2_{L^2}\, ds\Big)} &&
\\
& \le & \cz e^{\Cdz T}\Big(  \EE |u_0|_{L^2}^2  + 2\kappa _u \EE \int_0^ t \int_\CO u(s,x)\chi^2(s,x)\xi(s,x)\, dx \, ds\Big).
\EEQSZ
By Corollary \ref{cor01} we know that for any $\ep_1,\ep_2>0$ we have
\DEQS 
\lqq{\bE \int_0^ t  \int_{\CO} \chi^2(s, x) \xi(s, x) u(s, x) \,dx\, ds
} &&
\\
& \le & c_1\EE|\xi_0|_{L^1}+c_2\EE|u_0|^2_{L^2} +c_3\bE|v_0|_{L^1}
 +c_4 \EE \int_0^ t |\xi(s)|_{L^2}^2\, ds
\\&&{}+c_5\bigg( \EE\sup_{0\le s\le T}|\eta(s)|_{L^2}^2 +  \EE\int_0^ T |\eta(s)|_{H^1_2}^2\, ds\bigg)^\frac 14
 \\&&{}+\ep\lk( \EE \sup_{0\le s\le T}|u(s)|_{L^2}^2+ \EE\int_0^ T|u(s)|_{L^2}^2\, ds\rk).
\EEQS
Collecting together both estimates,
choosing $\ep\in(0,\frac 12] $ and diving by $(1-\ep)$  gives for some $c_1,\ldots,c_5>0$
\DEQSZ \label{vorher01}  \nonumber
\lqq{ \lk\{  \EE \sup_{s \in [0, T]} |u(s)|_{L^2}^2  +\EE \Big(\int_0^ T |\nabla u(s)|^2_{L^2}\, ds\Big)\rk\}} &&
\\
& \le & c_1\EE|\xi_0|_{L^1}+c_2\EE|u_0|^2_{L^2} +c_3\bE|v_0|_{L^1}
 +c_4 \EE \int_0^ t |\xi(s)|_{L^2}^2\, ds
\\&&{}+c_5\bigg( \EE\sup_{0\le s\le T}|\eta(s)|_{L^2}^2 +  \EE\int_0^ T |\eta(s)|_{H^1_2}^2\, ds\bigg)^\frac 14.\nonumber
\EEQSZ
Note that due to the definition of $\mathscr{L}_1(u,v)$ we have
$$ \EE \sup_{0\le s\le T}|u(s)|_{L^2}^2+ \EE\int_0^ T|u(s)|_{L^2}^2\, ds\le K_1.
$$
Substituting above in \eqref{vorher01} and using the estimate in Proposition \ref{pr-xi_est}-(a) we obtain
for some new constants $c_5,\delta_5>0$
\DEQS 
\lqq{ \EE \sup_{s \in [0, T]} |u(s)|_{L^2}^2  +\EE \Big(\int_0^ T |\nabla u(s)|^2_{L^2}\, ds\Big)+\EE \sup_{0\le s\le T} |\xi(s)|_{L^p}^p
} &&
\\
& \le & c_1\EE|\xi_0|_{L^1}+c_2\EE|u_0|^2_{L^2} +c_3\bE|v_0|_{L^1}
+c_4 K_1^\frac 14 +c_5e^{\delta_5T}\EE|\xi_0|_{L^p}^p.
\nonumber
\EEQS
Since we are aiming to  proof
$$\EE\sup_{0\le s\le T}|\eta(s)|_{L^2}^2 +  \EE\int_0^ T |\eta(s)|_{H^1_2}^2\, ds+\EE \sup_{0\le s\le T} |\xi(s)|_{L^p}^p\le K_1
$$
we chose $K_1>0$ according to
\DEQS 
K_1
& \ge & c_1\EE|\xi_0|_{L^1}+c_2\EE|u_0|^2_{L^2} +c_3\bE|v_0|_{L^1}
+c_4 K_1^\frac 14+c_5e^{\delta_5T}\EE|\xi_0|_{L^p}^p.\nonumber
\EEQS
With this choice, inequality \eqref{toshowk2} is satisfied.
 In this way, we have found three constants $K_1,K_2>0$ and $K_3>0$ such that,
if $(\chi,\eta)\in \mU$, then $(u,v):=\mathscr{T}[(\chi,\eta)]\in\mU$.
\end{proof}

\medskip

\begin{claim}\label{secondclaim}
For any 
$K_1, K_2>0$, and $K_3>0$, 
the map
$$\mathscr{T} :\mU \to \mathcal{M}_{\mathfrak{A}}
$$ is continuous.
\end{claim}

\begin{proof}[Proof of Claim \ref{secondclaim}:]
The continuity follows by a combination of the Proposition \ref{help_cont}, Proposition \ref{help_cont1}, and Proposition \ref{help_cont4}.
First, it follows by Proposition \ref{help_cont} that there exists some  $\gamma\in(0,1]$ and $\delta_1,\delta_2>0$such that we have
\begin{equation*}
\begin{split}
&\mathbb E \sup_{0 \leq s \leq T}  \big| u_1(s) - u_2(s) \big|_{H^{1-\rho}_2}^2
\\ &\leq C(K_1,K_2,T) \times \Big[ \Big\{ \mathbb{E} \big[ \sup_{s \in [0, T]} \big|\chi_1(s) - \chi_2(s) \big|^2_{H^{1-\rho}_2} \big] \Big\}^{\delta_1} + \Big\{ \mathbb{E} \big[ \sup_{s \in [0, T]} \big| \xi_1(s) - \xi_2(s) \big|_{L^1}^\gamma  \big] \Big\}^{\delta_2}\,\Big].
\end{split}
\end{equation*}
Next, by Proposition \ref{help_cont1} we have for $m>q=1$, $p\ge \frac {2mq}{m-q}$ and $r=\frac 1{1-\gamma}$
\DEQS
 \mathbb{E} \big[ \sup_{s \in [0, T]} \big| \xi_1(s) - \xi_2(s) \big|_{L^1}^\gamma \big]\le
C(K_2,K_3) \Big( \EE  \big[ \sup_{0\le s\le T}  \left| {\v_1(s)} -{\v_2(s)}\right|_{L^m}  \big] \Big)^\frac 12 .
 \EEQS
 Finally, Proposition \ref{help_cont4} gives for $\rho>2d(1-\frac 1m)$
 \begin{equation*}
\begin{split}
&\mathbb E \sup_{0 \leq s \leq T}  \big| v_1(s) - v_2(s) \big|_{L^2}
\leq C(K_3,T) \times \Big[ \mathbb{E}  \sup_{s \in [0, T]} \big|\chi_1(s) - \chi_2(s) \big|^2_{H^{1-\rho}_2}  \Big]^{\gamma}.
\end{split}
\end{equation*}
Collecting altogether and taking into account our choice of $p=20$ gives  the assertion.
\end{proof}

\medskip
\begin{claim}\label{thirdclaim}
For any 
$K_1, K_2>0$, and $K_3>0$, 
the map
$$\mathscr{T} :\mU \to \mathcal{M}_{\mathfrak{A}}
$$ is compact. 
\end{claim}

\begin{proof}[Proof of Claim \ref{thirdclaim}:]
First, let us define  the convolution operator $\mathfrak{C}$ by
$$ (\mathfrak{C}f)(t):= \int_0^t e^{-(t-s)(A- \Upsilon_vI)}f(s) \, ds,
$$
and the convolution operators $\mathfrak{S}_j$, $j=1,2$, by
$$ (\mathfrak{S}_jf)(t):= \int_0^t e^{-(t-s)(A- \Upsilon_vI)}f(s) \, dW_j(s),\quad j=1,2.
$$
We can write the solution $(u,v)$ for $t>0$ as
\DEQS
u(t)&=&e^{-t(A- \Upsilon_vI)}u_0+\int_0^t e^{-(t-s)(A- \Upsilon_vI)}\chi^2(s)\xi(s) \, ds+\int_0^t e^{-(t-s)(A- \Upsilon_vI)}u(s)\, dW_1(s)
\\
&=& e^{-t(A- \Upsilon_vI)}u_0+\big( \mathfrak{C}[\chi^2\xi] \big)(t) +\big( \mathfrak{S}_1[u] \big)(t),
\EEQS
and
\DEQS
v(t)&=&e^{-t(A-\Upsilon_vI)}v_0+\int_0^t e^{-(t-s)(A- \Upsilon_vI)}\chi^2(s)\, ds+\int_0^t e^{-(t-s)(A- \Upsilon_vI)}v(s)\, dW_2(s)
\\
&=& e^{-t(A- \Upsilon_vI)}v_0+\big( \mathfrak{C}[\chi^2] \big)(t) +\big( \mathfrak{S}_2[v] \big)(t).
\EEQS
%
%
Since $v$ depends only on $\chi$, we study this process at first. Let us first examine the perturbation by $\chi^2$.
Since $(\chi,\eta)\in\mU$, by Proposition \ref{interp_rho}, for  $p=\frac 52$
  there exists a $C>0$ such that
  $$\EE \|\chi\|_{L^p(0,T;L^4)}^2\le C \EE \lk(\|\chi\|^2 _{L^\infty(0,T;H^{1-\rho}_2)}+\|\chi\|^2_{L^r(0,T;H^{2-\rho}_2)}\rk).
  $$
In particular, we know for $f=\chi^2$, that there exists a constant $ C(K_1,K_2,K_3)>0$ such that
$$\EE \|f\|_{L^\frac 54(0,T;L^2)}\le   C(K_1,K_2,K_3).
$$
It follows by
by Lemma \ref{L:reg} with $\alpha=1$, $q=\frac 52$, $X=L^2(\CO)$, $\gamma=0$, $\delta>0$, and $\beta>0$ with $\beta+\frac\delta2<1-\frac 1q$, that we have
\DEQS 
\mathfrak{C}[\chi^2]\in C^\beta_b(0,T;H^{\delta}_2(\CO)),
\EEQS
and
\DEQSZ
\EE \|\mathfrak{C}[\chi^2]\|_{C^\beta_b(0,T;H^{\delta}_2)} =\EE \|\mathfrak{C}[{f}]\|_{C^\beta_b(0,T;H^{\delta}_2)} \nonumber
\le C\EE \|{f}\|_{L^q(0,T;L^2)}\le C(K_1,K_2,K_3).
\EEQSZ
%
Since we do not have the integrability of $v$ for $p>2$, we cannot rely on the factorisation method and have to use another argument.
Observe, firstly, we know by Theorem \ref{theou1} that for all $(\chi,\eta)\in \mU$,  $\EE\sup_{0\le s\le T}  |v(s)|_{L^2}\le C(K_1,K_2,K_3)$.
By Theorem 6.3 \cite{Ethier+Kurtz} we have to show that
\begin{enumerate}
\item[(a)] for every number $t\in\mathbb{Q}\cup[0,T]$ for any $\ep>0$
there exists a compact set $\Gamma^\ep_t\subset L^2(\CO)$ such that
for any  $(\chi,\eta)\in \mU$ $\PP(v(t)\in \Gamma^\ep_t)\ge 1-\ep$, where $v$ denotes the solution of \eqref{eqv1sto};
\item[(b)] for any $T>0$ and any $\ep>0$ there exists a $\delta>0$ such that we have  $\PP\lk( \omega(\mathfrak{v},\delta )\ge \ep\rk)\le \ep$
for all  $(\chi,\eta)\in \mU$ (and $v$ denotes the solution of \eqref{eqv1sto}), where $\omega(\mathfrak{v},\delta)=\sup_{|t-s|\le \delta}|\mathfrak{v}(t)-\mathfrak{v}(s)|$.
\end{enumerate}
Fix $\delta>0$ and $\Gamma^\ep_t:=\{ w\in H^\delta_2(\CO): |w|\le t^{-\delta}|v_0|/\ep\}$. Due to the compact embedding $H_2^\delta(\CO)\hookrightarrow L^2(\CO)$ and since
straightforward calculations gives that $e^{-t(A- \Upsilon_vI)}v_0\in \Gamma_t$, we know $t\mapsto e^{-t(A- \Upsilon_vI)}v_0$ satisfies (a).
Secondly, straightforward calculations gives $\sup_{t\in(0,T]}\EE |\mathfrak{S}v(t)|_{H^\delta_2}<\infty$. The Chebycheff inequality gives (a).
Finally, we have to estimate the modulo of continuity, i.e. to verify condition (b). However,
splitting the difference we arrive at the following sum
$$
\EE \sup_{t\le s\le t+h}\int_t^s e^{-(s-r)(A- \Upsilon_v)}v(r)dW_2(r)+ \EE  \sup_{t\le s\le t+h}(I-e^{-s(A- \Upsilon_v)})\mathfrak{S}v(t),
$$
straightforward estimates and using estimates (4) of \cite{seidler} gives the desired result. %
Next, we will investigate the regularity of $u$.
First, we will show that there exists a constant such that there exists some $l=\frac 52$ 
such that  for all $(\chi,\eta)\in\mU$  and $(u,v)=\mathscr{T}(\chi,\eta)$
and for $r=2$
$$\EE \Big( \int_0^ T |\chi^2(s)\xi(s)|_{L^1}^r\, ds \Big)^\frac 1l\le C.
$$
Here, let us note that we get for $q=4$ and $q'=\frac 43$ conjugate 
\DEQS
\lqq{ \EE \Big( \int_0^ T |\chi^2(s)\xi(s)|_{L^1}^r\, ds \Big)^\frac 1l\le \EE \Big( \int_0^ T |\chi(s)|^{2r}_{L^{2q'}} | \xi(s)|_{L^q}^r\, ds \Big)^\frac 1l}
&&
\\
&\le &\EE \sup_{0\le s\le T}  |\xi(s)|^\frac  rl_{L^q}\, \Big(\int_0^ T|\chi(s)|_{L^{2q'}}^{2r}\, ds \Big)^\frac 1l
\\
&\le &
\Big\{\EE \sup_{0\le s\le T}  |\xi(s)|^4_{L^4} \Big\}^\frac 1{2l}\Big\{ \EE  (\int_0^ T|\chi(s) |^{2r}_{L^{2q'}}\, ds)^\frac 2{2l-1} \Big\}^\frac {2l-1}{2l}  .
\EEQS
Due to Proposition \ref{pr-xi_est}-(a)  we know that the first term is bounded.
By Proposition \ref{interp_rho} and the fact that $\frac d2-(1-\rho)\le \frac 2{2r}+\frac d{q'}$ we
have
\DEQSZ\label{finiteint}
\|\chi\|_{L^{2r}(0,T;L^{2q'})}\le  C\lk(\|\chi\| _{L^\infty(0,T;H^{1-\rho}_2)}+\|\chi\|_{L^r(0,T;H^{2-\rho}_2)}\rk).
\EEQSZ
Now, the second term is bounded due to the definition of  $\mU$, \eqref{finiteint}. In particular, we have
$$
\EE \|\chi\|_{\mathbb{H}_{1-\rho}}^2\le K_1.
$$
If  $l$ is chosen sufficiently large, i.e. $\frac 2{2l-1}\le \frac 12$.
Therefore, we know that for any $\delta>0$  there exists a constant $C(K_1,K_2,K_3)>0$ such that
$$
\EE \|\chi^2\xi\|^\frac 1l_{L^{r}(0,T;H^{-1-\delta}_2                                                                                                   )}\le C(K_1,K_2,K_3), \quad (\chi,\xi)\in \mU.
$$
Now,
 by Lemma \ref{L:reg} with $\alpha=1$, $q=r$, and $\gamma=0$ gives that for any $\beta+\delta<1-\frac 1r $, $0<\delta<\frac {\rho-1}4$
$$
\mathfrak{C}[\chi^2\xi]\in C^\beta_b(0,T;H^{1-\rho}_2(\CO)).
$$
In particular, there exists a constant $C>0$ such that  for $(\chi,\eta)\in\mU$  and
we have
\DEQSZ\label{ll3}
\EE\|\mathfrak{C}[\chi^2\xi]\|^\frac 12_{C^\beta_b(0,T;H^{\delta-1}_2)}\le C.
\EEQSZ
Secondly, we know since $(u,v)\in\mU$ that
$$\EE \sup_{0\le s\le T} |u(s)|_{H^{1-\rho}_2}^2+\EE \int_0^ T  |\nabla u(s)|_{H^{1-\rho}_2}^2\, ds \le K_1.
$$
%
Then, Corollary \ref{brz:con} gives for    $\beta+\frac 1p=\frac 12$, $X=H^{1-\rho}_2(\CO)$, $\delta=0$, $\nu=\frac 14$
$$
\mathfrak{S}_1[u]\in C^\beta_b(0,T;H^{1-\rho}_2(\CO)).
$$
In particular, there exists a constant $C>0$ such that  for $(\chi,\eta)\in\mU$  and $(u,v)=\mathscr{T}(\chi,\eta)$ we have
\DEQSZ\label{ll4}
\EE\|\mathfrak{S}_1[v]\|^\frac 12_{C^\beta_b(0,T;L^2)}\le C.
\EEQSZ

\del{
Similarly, since $(\chi,\eta)\in\mU$, by interpolation,  for any ${\tilde \rho}\in(0,1)$ and $1\le p\le \frac1{\tilde \rho}$  there exists a $C>0$ such that
  $\EE \|\chi^2\|_{L^p(0,T;H^{\tilde \rho}_2)}\le C$.
 By Lemma \ref{L:reg} with $\alpha=1$, $q=p$, and $\beta+\frac\delta2<1-\frac 1q$, we have
$$
\mathfrak{C}[\chi^2]\in C^\beta_b(0,T;H^\delta_2(\CO)).
$$
Finally,
Corollary \ref{brz:con} gives for  $\beta+\frac \delta2+ \frac 1p=\frac 12$,
$$
\mathfrak{S}_2[v]\in C^\beta_b(0,T;H^\delta_2(\CO)).
$$}
\noindent
Collecting together  \eqref{ll3}, and \eqref{ll4},
 there exist $\beta,\delta>0$ such that
$$
u\in  C^\beta_b(0,T;H_2^{\delta-1}(\CO)) 
\cap L^\infty(0,T; L^2(\CO)). 
$$
In particular, there exists a $C>0$ such that  for all $(\chi,\eta)\in\mU$ 
$$
\EE\Big( \|u\|_{ C^\beta_b(0,T;H_2^{\delta -1})}+\| u\|_{ L^\infty(0,T; L^2)} \Big)^\frac 1l \le C.
$$
The Chebycheff inequality 
and the Prohorov Theorem
gives compactness of the operator $\mathscr{T}$.

\end{proof}

\paragraph{\bf Step IV} {\bf Conclusion:}

By Theorem \ref{schauder}, there exists a probability space
$\tilde{\mathfrak{A}} = (\tilde{\Omega}, \tilde{\mathcal{F}}, \tilde{\mathbb{F}}, \tilde{\mathbb{P}}),$ a Wiener process $\tilde{W}=(\tilde{W}_1, \tilde{W}_2)$ defined on the probability space $\tilde{\mathfrak{A}}$ and an element
$(\tilde{u},\tilde v) \in \mU$
such that $\tilde{\mathbb{P}}$-a.s.
$$\mathscr{T}_{\tilde{\mathfrak{A}}}\big[  ( \tilde{u},\tilde v)\big] (t) = (\tilde{u}(t),\tilde v(t)) \quad \forall t \in [0, T].$$
By the construction of the operator $\mathscr{T}$, the pair $(\tilde{u}, \tilde{\v})$ solves the system \eqref{equ1i}-\eqref{eqv1i}.
In addition,
the solution belongs to $\mU$, hence the estimate in item (iii) 
is satisfied.
\end{proof}

\smallskip

\section{Auxiliary Results and The Technical Propositions}\label{technics}

In this particular section, we will prove various results and energy bounds which are crucial for the proof of the main result in Section \ref{proof}. We assume that the Hypotheses \ref{init} and \ref{wiener} are satisfied. We split this section into four main subsections. In Subsection \ref{PoS-v}, we discuss the properties of the equation \eqref{eqv1sto} and prove several uniform bounds for $\v$. Next, in Subsection \ref{PoS-u}, we analyse the properties of the equation \eqref{eqv2sto} and show various energy bounds for $u$. Next, in Subsection \ref{cont}, we state the necessary results to prove the continuity of the map $\mathscr{T}$. Finally, in Subsection \ref{comp}, we state the necessary results to prove the compactness of the mapping $\mathscr{T}$.

\subsection{Properties of the system \eqref{eqv1sto}}\label{PoS-v}

In the first part of this section, we will focus on the equation \eqref{eqv1sto}. Let us recall the system for the convenience of the reader. We are interested in the solution of the following system for $x\in \CO,\, t> 0$:
\begin{align}\label{eq-vd}
\lk\{ \barray
d{\v}(t) &=& \big[ r_\v \Delta \v(t) + \kappa_v\, \chi^2(t) - \Upsilon_v \v(t) \big]\,dt +\sigma_{\v} \v(t) \,d W_2(t),
\\
\v(0) &=& \v_0. \phantom{\Big|}
\earray\rk.
\end{align}
At first, we will show the existence of a unique solution to the system \eqref{eq-vd} and then we will prove the solution is positive,
suppose $(\chi,\eta)\in \mU$.

\medskip

\begin{theorem}\label{theou1}
Suppose  $K_1, K_2 >0$, $K_3>0$, and $T>0$ are fixed.
For any $\v_{0}\in L^{2}(\Omega; L^2(\CO))$ and any $(\chi,\eta)\in \mU$, the system \eqref{eq-vd} has a unique solution $\v: \Omega \times [0, T] \times \CO \to \mathbb{R}$ such that
$$\EE \sup_{0\le s\le T} |\v(s)|_{L^2} +\EE\Big(\int_0^ T |\nabla {\v}(s)|_{L^2}^2\, ds\Big)^\frac 12 <\infty.
$$
Moreover, this solution is positive provided that $\v_0 >0, \,\mathbb{P} \times Leb$-a.s.
\end{theorem}

\begin{proof}
\del{Since $(\chi,\eta)\in\mU$, by Proposition \ref{interp_rho},  for any ${\tilde \rho}\in(0,1)$ and $p=2$
  there exists a $C>0$ such that
  $$\EE \|\chi\|_{L^p(0,T;L^4)}^2\le \EE C\lk(\|\chi\|^2 _{L^\infty(0,T;H^{1-\rho}_2)}+\|\chi\|^2_{L^r(0,T;H^{2-\rho}_2)}\rk).
  $$
In particular, we know for $f=\chi^2$,
$$\EE \|f\|_{L^2(0,T;L^2)}\le  \EE C\lk(\|\chi\|^2 _{L^\infty(0,T;H^{1-\rho}_2)}+\|\chi\|^2_{L^r(0,T;H^{2-\rho}_2)}\rk)<\infty.
$$
It follows}
Since $(\chi,\eta)\in\mU$,  $\chi$ is adapted. Next, by Proposition \ref{interp_rho},  it follows that
for $\frac d2\le \frac 1{2}+\frac d{l_2}$
there exists a constant $C>0$ such
\begin{equation*}
\EE \|\chi\|_{L^{l_1}(0,T;L^{l_2})}^2\le C\,\EE \lk(\|\chi\|^2 _{L^\infty(0,T;L^2)}+\|\chi\|^2_{L^2(0,T;H^{1}_2)}\rk).
\end{equation*}
Since $(\chi,\eta)\in\mU$,  we know that there exists a constant $C_1>0$ such
\begin{equation}
\EE \|\chi\|_{L^{l_1}(0,T;L^{l_2})}^2\le C_1.
\label{eq:v-bound-L1}
\end{equation}
Due to the embedding $L^\frac {l_2}2(\CO)\hookrightarrow H^{-\gamma }_p(\CO)$ for $p\ge \frac {l_2}2$ and $\gamma\ge 0$ with $\gamma+\frac dp\ge \frac {2d}{l_2}$
we have
\DEQS
\EE \|\chi^2\|_{L^2(0,T;H^{-\gamma}_p)}\le C\, \EE \|\chi^2\|_{L^2(0,T;L^\frac {l_2}2)}\le C \EE \|\chi\|_{L^{4}(0,T;L^{l_2})}^2 \le CC_1.
\EEQS
Let
$$\mathfrak{C}(\eta)(t):= \int_0^ t e^{-(t-s)\Delta}\eta(s)\, ds,
$$
be the deterministic convolution. Then, we have by Lemma \ref{L:reg} for any $ p\ge  \frac {l_2}2$, $\gamma\in\RR$, and $\delta\in (0,2)$ 
$$
\|\mathfrak{C}(\eta)\|_{L^2(0,T;H^{\delta-\gamma}_p)}\le \|\eta\|_{L^2(0,T;H^{-\gamma }_p)}.
$$
Now,
{{with $\eta=\chi^2$}} and 
for any $p\ge 2$ 
we have
$$ \EE \|\mathfrak{C}(\chi^2)\|_{L^2(0,T;H^1_p)} \le C \EE \|\chi\|_{L^{4}(0,T;L^{l_2}/2)}^2<
CK_1.
$$
Let $\alpha=1$, $q=2$, and $\gamma=0$. Now,   for any $1\le p<6$ we have 
$$
\mathfrak{C}[\,\chi^2\,]\in C^\beta_b(0,T;L^p(\CO)).
$$
In particular,  there exists a constant $C(K_1,K_2,K_3)>0$  such that  for $(\chi,\eta)\in\mU$  and $(u,v)=\mathscr{T}(\chi,\eta)$ we have
\DEQSZ\label{ll3}
\EE\|\,\mathfrak{C}[\,\chi^2]\,\|_{C^\beta_b(0,T;L^2)}\le C(K_1,K_2,K_3).
\EEQSZ
It remains to tackle $\Upsilon_v v$, but since $\Upsilon_v :L^2(\CO)\to L^2(\CO)$ and
$\Upsilon_v :H^1_2(\CO)\to H^1_2(\CO)$is bounded and linear, the term can be tackled by standard calculations.
To give an estimate of the stochastic convolution term in $L^2(0,T;H^1_2(\CO))$ we use the fact that
$\sum_{j=1}^\infty \lambda_j^{-2\gamma_2} |e_j|^2_{L^\infty}$ is bounded, i.e.\ $\gamma_2>\frac d2$.
To give an estimate of the stochastic convolution term we use estimate (4) of the main Theorem in \cite{seidler}.
%
A fix point argument leads the existence and uniqueness of a solution.
%
For the positivity of the solution to \eqref{eq-vd} we refer to the work of Barbu, Da Prato and R\"ockner (Section 2.6 in \cite{Barbu+DaPrato+Rockner_2008} see also \cite{BDPR2016}).  For the positivity
we refer to Section \ref{positivity}, where we have shown  in Theorem \ref{positivitytheorem}  the positivity of
the solution of $w(t)=w_0+\int_0^ t e^{A(t-s)}w(s)\, dW_2(s)$ along the proof of Theorem 2.3 in \cite{TessitoreZabczyk1998}.
The positivity of $v$ follows by comparison results, taking into account that $\chi^2$ is non--negative.
 \end{proof}

\subsubsection{{\bf Uniform bounds on $\v$}}\label{SUbv}

In fact, uniform bounds on $\v$ are difficult to achieve, however, we can show several uniform bounds on $\xi$, where
 $\xi = \v^{-1}$.
 Applying the It\^o formula to $\phi_1(v)=v^{-1}$ (see \cite{prato}) we infer that
\DEQSZ 
\nonumber
\lqq{ d \xi(t)= -\xi^2(t) \lk(r_{\v} \Delta \v(t) +\kappa _v\chi^2(t) + \sigma_v \,(\mbox{Id}+A)^{-\gamma_2}\,v(t) \rk)\, dt  +  \mu_{\v} \xi(t) \,dt}
&&\\
&&+ \frac 12 \Tr\lk[D^2 \phi_1(v(t))\lk[ \sigma_v v(t) \sqrt{Q_2}\rk]  \lk[\sigma_v v(t) \sqrt{Q_2}\rk]  ^{\ast} \rk]
- \sigma_{\v} \xi(t)\,dW_2(t). \nonumber
\EEQSZ 
Thus, we have the following equation in $\xi$:
\begin{align}\label{eq-xi}
\lk\{ \barray
d{\xi}(t) &=& -\xi^{2}(t)\big[ r_\v \Delta \v(t) + \kappa_v\, \chi^2(t) + \sigma_v \,(\mbox{Id}+A)^{-\gamma_2}\,v(t) \big] dt
\\ &\qquad+& \big[ \mu_v+\frac 12 \sigma_v^2\,S(\gamma_2)\big]  \xi(t)dt - \sigma_{\v} \xi(t) d W_2(t),
\\
\xi(0) &=& \v^{-1}_0, \phantom{\Big|}
\earray\rk.
\end{align}
where for $\gamma>0$, $S(\gamma)$ denotes the following operator
$$
S(\gamma) f:=\sum_{k=1}^\infty (1+\lambda_k)^{-\gamma}\la e_k,f\ra e_k, \quad j=1, 2,
$$
 which is bounded  due to \eqref{EVsup}.
Here, it should be noted that this formula holds if $v$ is positive.

\begin{proposition}\label{pr-xi_est}
For any $K_1, K_2>0$ and $K_3>0$ there exist
\begin{enumerate}
\item[$(a)$]\stepcounter{con}
 constants $\Con,\Cdelta>0 $ and $\tCon>0$ such that for any $T>0$, for any initial conditions $(u_0, \v_0)$ satisfying Hypothesis \ref{init},
 for any  $(\chi,\eta)\in\mU$ and $(u,\v)=\SCT[(\chi,\eta)]$
we have
{\small{
\begin{equation}\label{xi^2q_est}
\begin{split}
&\bE \sup_{0 \leq s \leq T} |\xi(s)|_{L^{p}}^{p} + 2 p (p+1) r_{\v} \,\bE \int_0^T \int_{\CO} \xi^{p+2}(s, x) |\nabla \v(s, x)|^2 \,dx\,ds
\\ &\quad+ \tCon \bE \int_0^T \int_{\CO} \xi^{p+1}(s, x) \chi^2(s, x) \,dx\,ds \leq \Con e^{\Cdelta T}\, \EE|\xi_0|^{p}_{L^p};
\end{split}
\end{equation}
}}
\item[$(b)$]\SC
constants  $\Con,\Cdelta>0 $ such that for any $T>0$, for any  initial conditions $(u_0,v_0)$ satisfying Hypothesis \ref{init},
 for any  $(\chi,\eta)\in\mU$ and $(u,\v)=\SCT[(\chi,\eta)]$ we have
\begin{equation}\label{eq-6.13}
\begin{split}
& \sup_{0 \leq s \leq T}\mathbb{E}  |\xi(s)|_{L^1} + 4 r_{\v}\, \EE \int_{0}^{T} \int_{\CO} \xi^3(s, x) |\nabla \v(s, x)|^2 \,dx\,ds
\\ &\quad+ \kappa _v \,\EE \int_{0}^{T} \int_{\CO} \chi^2(s, x) \xi^2(s, x) \,dx \,ds \leq \Con e^{\Cdelta T} \,\EE |\xi_0|_{L^1};
\end{split}
\end{equation}
\item[$(c)$]\SC
  a constants  $C_3, \tilde{C}_3>0 $ such that
for any $T>0$, for any   initial conditions $(u_0, \v_0)$ satisfying Hypothesis \ref{init}, for all
 $(\chi,\eta)\in\mU$ and $(u,\v)=\SCT[(\chi,\eta)]$
 we have
\begin{equation}\label{eq-6.131}
\begin{split}
&\bE \Big( \int_{\CO} |\ln \xi(t, x)|\,dx \Big) + r_{\v} \bE \Big(\int_0^t \int_{\CO} \xi^2(s, x) |\nabla \v(s, x)|^2\,dx\,ds \Big)
\\ &\quad{}+ \kappa_v \bE \Big( \int_0^t \int_{\CO} \chi^2(s, x) \,\xi(s, x) \,dx\,ds \Big)
\leq
\bE|v_0|_{L^1}+ C_3 \, \bE \Big( \int_{\CO} |\ln \xi_0(x)| \,dx \Big) +  \tilde{C}_3 T.
\end{split}
\end{equation}

\end{enumerate}
\end{proposition}

\begin{proof}

{\bf Proof of Part (a).}
%
For $k \in \mathbb{N}$ and $p \geq 2,$ let us define the stopping time $\tau_k := \inf\{ t>0 \ | \ |\xi(t)|^p_{L^p} \geq k\}$.
\del{and
\begin{align}
\xi(t):= \begin{cases}
&\xi(t), \mbox{ if } 0\le t\le \tau_k;
 \\ &\mbox{otherwise } \ 0.
 \end{cases}
 \end{align}}
Applying the It\^o formula (see \cite{Kry-Ito}) to the function $\phi_2(z)=|z|^{p}$, for $p\geq 2$, and taking into account that
$$
\Tr\lk[D^2 \phi(\xi(t))\lk[\sigma_v \xi(t) \sqrt{Q_2}\rk]  \lk[\sigma_v \xi(t) \sqrt{Q_2} \rk]  ^{\ast} \rk]
 =p(p-1)\,\sigma_v^2 \,\xi^{p}(t)\,\sum_{k=1}^\infty (1+\lambda_k)^{-\gamma_2} |e_k|^2,
 $$
for ($\phi=\phi_2\circ\phi_1$) we infer for $t\in [0,T\wedge \tau_k]$
{\small{
\begin{equation}\label{eq-dxi^rho}
\begin{split}
d |\xi(t)|^{p}_{L^p} &= \Big[-p \int_{\CO} \Big( r_{\v} \xi^{p+1}(t, x)\, \Delta \v(t, x) + \kappa_v\, \xi^{p+1}(t, x)\, \chi^2(t, x) - \mu_\v  \,\xi^{p}(t, x) \Big)\,dx
\\ &\quad+  \frac{p(p-1)}{2}\sigma_v^2 S(\gamma_2) \int_{\CO} \xi^{p}(t, x)\,dx  - p \int_{\CO} \sigma_v \xi^{p+1}(t, x)\,(\mbox{Id}+A)^{-\gamma_2}\,v(t, x) \,dx \Big]\,dt
\\ &\quad- p \sigma_{\v} \int_{\CO} \xi^{p}(t, x) \,dW_2(t, x).
\end{split}
\end{equation}
}}
Using the positivity of $v$, the maximum principle and $\sigma_v>0$, we conclude that the term $- p \int_{\CO} \sigma_v \xi^{p+1}(t, x)\,(\mbox{Id}+A)^{-\gamma_2}\,v(t, x) \,dx$ can be dropped from the right hand side. Then, applying integration by parts to the first term in right hand side of \eqref{eq-dxi^rho} we get for $t\in[0,T\wedge \tau_k]$
{\small{
\begin{equation}\label{xi^Lrho}
\begin{split}
d |\xi(t)|_{L^{p}}^{p} &\leq \Big[-p (p+1) r_{\v} \int_{\CO} \xi^{p}(t, x) \,\nabla \xi(t, x) \cdot \nabla \v(t, x) \, dx + r_{\v} p \int_{\partial \CO} \vec{n} \cdot \nabla \v(t, x)\, \xi(t, x) \,dx
\\ &\quad- p \int_{\CO} \big[\kappa _v\,\xi^{p+1}(t, x)\, \chi^2(t, x)  \big] \,dx
\\ &\quad+p \,\Big(\mu_\v +\sigma_v^2 \,S(\gamma_2)\,\frac{ (p-1)}{2} \Big) \int_{\CO} \xi^{p}(t, x)\,dx \Big] \, dt- p \sigma_{\v} \int_{\CO} \xi^{p}(t, x)\,dW_2(t, x)
\\ &\leq \Big[ -p (p+1) r_{\v} \int_{\CO} \xi^{p+2}(t, x) (\nabla \v(t, x))^2 dx - p \int_{\CO} \big[\kappa _v\,\xi^{p+1}(t, x) \chi^2(t, x)  \big] \,dx
\\ &\quad+p \,\Big(\mu_\v +\sigma_v^2 \,S(\gamma_2)\,\frac{ (p-1)}{2} \Big)  \int_{\CO} \xi^{p}(t, x)\,dx \Big] \,dt - p \sigma_{\v} \int_{\CO} \xi^{p}(t, x) \,dW_2(t, x).
\end{split}
\end{equation}
}}
\noindent
Here, the positive term $p \int_{\partial \CO} \vec{n} \cdot \nabla \v\, \xi \,dx$ vanishes due to the Neumann boundary conditions.

In the integral form, taking supremum over $s \in [0, t \wedge \tau_k]$ and taking expectation we infer that
\begin{equation}\label{xi^Lq-sup}
\begin{split}
&\bE \sup_{0 \leq s \leq t \wedge \tau_k} |\xi(s)|_{L^{p}}^{p} + p (p+1) r_{\v} \bE \int_0^{t \wedge \tau_k} \int_{\CO} \xi^{p+2}(s, x) |\nabla \v(s, x)|^2 \,dx\,ds
\\ &\quad+ p\,\kappa _v\,\bE \int_0^{t \wedge \tau_k} \int_{\CO} \xi^{p+1}(s, x) \chi^2(s, x) \,dx\,ds
\\ &\leq \bE |\xi_0|_{L^{p}}^{p} + C(\mu_{\v}, S(\gamma_2), \sigma_{\v}, p) \,\bE \int_0^{t \wedge \tau_k} |\xi(s)|_{L^p}^{p} \,ds + \bE \sup_{0 \leq s \leq t \wedge \tau_k} \Big| p \,\sigma_{\v} \int_0^{s \wedge \tau_k} \int_{\CO} \xi^{p}(r, x)\,dW_2(r, x) \Big|.
\end{split}
\end{equation}

\noindent
Using the Burkholder-Davis-Gundy inequality,  the Cauchy-Schwarz inequality, and the Young inequality we infer
{\small{
\begin{align}\label{mbbExi-2}
\begin{split}
&\bE \sup_{0 \leq s \leq t \wedge \tau_k}  \Big| p \sigma_{\v} \int_0^{s \wedge \tau_k} \int_{\CO} \xi^p(r, x)\,dW_2(r, x) \Big| \leq p^2 \sigma^2_{\v} \,C\,\mathbb{E}\bigg(\int_0^{t \wedge \tau_k} \bigg| \int_{\CO} \xi^{p}(s, x)\,dx \bigg|^2 \,ds \bigg)^{1/2}
\\ &\leq p^2 \sigma^2_{\v} \,C\, \mathbb{E}\bigg( \sup_{s \in [0, t \wedge \tau_k]} |\xi(s)|_{L^{p}}^{p}  \int_0^{t \wedge \tau_k} \int_{\CO} \xi^p(s, x) \,dx \,ds \bigg)^{1/2}
\\ &\leq \,C\,\mathbb{E}\bigg( \Big\{\sup_{s \in [0, t \wedge \tau_k]} |\xi(s)|_{L^{p}}^{p}\Big\}^{1/2}  \Big\{ p \sigma_{\v} \int_0^{t \wedge \tau_k} |\xi(s)|_{L^{p}}^{p} ds\Big\}^{1/2} \bigg)
\\ &\leq \frac{\epsilon}{2} \mathbb{E} \sup_{s \in [0, t \wedge \tau_k]} |\xi(s)|_{L^{p}}^{p} +  \frac{p\sigma_{\v}}{2\epsilon} \mathbb{E} \int_0^{t \wedge \tau_k} |\xi(s)|_{L^{p}}^{p} ds.
\end{split}
\end{align}
}}

\noindent
Using this in \eqref{xi^Lq-sup}, choosing $\epsilon=\frac 12$ and finally rearranging we obtain,
\begin{equation}\label{xi^Lq-sup1}
\begin{split}
&\bE \sup_{0 \leq s \leq t \wedge \tau_k} |\xi(s)|_{L^{p}}^{p} + 2 p (p+1) r_{\v} \bE \int_0^{t \wedge \tau_k} \int_{\CO} \xi^{p+2}(s, x) |\nabla \v(s, x)|^2 \,dx\,ds
\\ &\quad+ 2p\,\kappa _v\,\bE \int_0^{t \wedge \tau_k} \int_{\CO} \xi^{p+1}(s, x) \chi^2(s, x) \,dx\,ds
\\ &\leq \bE |\xi_0|_{L^{p}}^{p} + C(\mu_{\v}, S(\gamma_2), \sigma_{\v}, p) \,\bE \int_0^{t \wedge \tau_k} |\xi(s)|_{L^p}^{p} \,ds.
\end{split}
\end{equation}

\noindent
Dropping the second and third terms from the left hand side and applying the Gronwall lemma we infer that, there exists a constant $C>0$ depending on $\mu_{\v}, S(\gamma_2), \sigma_{\v}, p$ such that
\begin{align*}
&\bE \sup_{0 \leq s \leq T \wedge \tau_k} |\xi(s)|_{L^{p}}^{p} \leq \bE |\xi_0|_{L^{p}}^{p} e^{ (T \wedge \tau_k) C(\mu_{\v}, S(\gamma_2), \sigma_{\v}, p) }. 
\end{align*}
By the Chebyscheff inequality, we observe that
as $k \to \infty$,  we have $\PP$--a.s. $T \wedge \tau_k \to T$. Therefore, taking the limit $k\to\infty$ we infer
\begin{equation}\label{estabove}
\begin{split}
&\bE \sup_{0 \leq s \leq T} |\xi(s)|_{L^{p}}^{p} \leq \bE |\xi_0|_{L^{p}}^{p} e^{ T C(\mu_{\v}, S(\gamma_2), \sigma_{\v}, p) }. 
\end{split}
\end{equation}
Substituting \eqref{estabove} in \eqref{xi^Lq-sup1}, we infer that there exist constants $\delta_1, C_1>0$ such that
we get the assertion in part (a).

\bigskip

\noindent
{\bf Proof of Part (b).}
Now we will prove \eqref{eq-6.131}. For clarity, we write the calculation without defining a stopping time $\tau_k := \inf\{ t>0 \ | \ |\xi(t)|_{L^1} \geq k\}$.
Tracing the proof one can see easily that taking into account the stopping time and being precise will not change the result.
Applying the It\^o formula to the function $\phi: L^1(\mathcal{O}) \to \mathbb{R}$ defined by $\phi(z)= \int_{\mathcal{O}} z(x)\,dx$, we infer that
{\small{
\begin{equation*}
\begin{split}
d |\xi(t)|_{L^1} &= -r_{\v} \int_{\CO} \xi^2(t, x) \Delta \v(t, x) \,dx\, dt - \kappa _v\int_{\CO} \xi^2(t, x) \chi^2(t, x) \,dx\, dt
\\ &\quad+ \mu_{\v} \int_{\CO} \xi(t, x) \, dx \,dt + \sigma_{\v} \int_{\CO} \xi(t, x) \, dW_2(t, x) - p \int_{\CO} \sigma_v \xi(t, x)\,(\mbox{Id}+A)^{-\gamma_2}\,v(t, x) \,dx\,dt
\\ &\leq -2 r_{\v} \int_{\CO} \xi^3(t, x) |\nabla \v(t, x)|^2 \,dx\, dt  - \kappa _v\int_{\CO} \xi^2(t, x) \chi^2(t, x) \,dx\, dt
\\ &\quad+ \mu_{\v} \,S(\gamma_2) \int_{\CO} \xi(t, x)\,dx + \sigma_{\v} \int_{\CO} \xi(t, x)\,dW_2(t, x),
\end{split}
\end{equation*}
}}
where the last term in the right hand side of the first equality is dropped using the similar arguments as in part (a)
and used integration by parts to get
{\small{
\begin{equation*}
\begin{split}
-r_{\v} \int_{\CO} \xi^2(t, x) \Delta \v(t, x) \,dx = r_{\v} \int_{\CO} \nabla (\xi^2(t, x)) \cdot \nabla \v(t, x) \,dx = -2r_{\v} \int_{\CO} \xi^3(t, x) \cdot |\nabla \v(t, x)|^2 \,dx.
\end{split}
\end{equation*}
}}

\noindent
Finally, taking integration from $0$ to $t$, then taking expectation, and rearranging we obtain,
\begin{equation*}
\begin{split}
&\mathbb{E} |\xi(s)|_{L^1} + 2 r_{\v} \mathbb{E} \int_{0}^{t} \int_{\CO} \xi^3(s, x) |\nabla \v(s, x)|^2 \,dx\,ds + \kappa _v\mathbb{E} \int_{0}^{t} \int_{\CO} \chi^2(s, x) \xi^2(s, x) \,dx \,ds
\\ &\quad\leq \mathbb{E} |\xi_0|_{L^1} + C(\mu_{\v}, S(\gamma_2))\, \mathbb{E} \int_{0}^{t} |\xi(s)|_{L^1} \,ds.
\end{split}
\end{equation*}

Dropping the second, third and fourth term from the left hand side and by applying the Gronwall lemma we infer that there exists a constant $\delta_2>0$ depending on $\mu_{\v}$ and a constant $C_2>0$ such that we get the assertion in part (b).

\bigskip

\noindent
{\bf Proof of Part (c).}
Again, for clarity, we write the calculation without defining a stopping time $\tau_k := \inf\{ t>0\mid \int\ln(\xi(t,x))\, dx \geq k\}$.
Tracing the proof one can see easily that taking into account the stopping time and being precise will not change the result.
We note that applying the It\^o formula for $\Phi(z)=-\ln(z)$
\begin{equation}
\begin{split}
\label{ca-dlnxi}
 &\ln(\xi(t))-\ln(\xi(0))
 =\int_0^ t \xi^{-1}(s) \big[ d \xi(s) \big]\, ds - \frac{\sigma_v^2}{2} \int_0^t \Tr \big[ D^2 \Phi(v(t)) \big(v(t) \sqrt{Q_2} \big) \big(v(t) \sqrt{Q_2} \big)^{\ast} \big] \,ds.
\end{split}
\end{equation}
From the previous energy estimates we have
\begin{equation*}
\begin{split}
d \xi(t, x) &= -\xi^{2}(t,x)\big[ r_\v \Delta \v(t, x) + \kappa_v\, \chi^2(t, x) +\sigma_v \,(\mbox{Id}+A)^{-\gamma_2}\,v(t,x) \big] dt +  \mu_\v \xi(t, x)dt
\\ &\quad- \frac 12 S(\gamma_2)\sigma_v^2 \xi(t, x)dt - \sigma_{\v} \xi(t, x) d W_2(t, x).
\end{split}
\end{equation*}
We observe that
\begin{equation}\label{lnxi-est}
\begin{split}
d \int_{\CO} \ln \xi(t, x) \,dx &= \Big[ -r_{\v} \int_{\CO} \xi(t, x) \Delta \v(t, x) \,dx - \kappa _v\int_{\CO} \xi(t, x) \chi^2(t, x)\,dx
\\ &\quad - \,\sigma_v  \int_{\CO} \xi(t, x)\,(\mbox{Id}+A)^{-\gamma_2}\,v(t, x) \,dx \Big] \,dt
\\ &\quad+ \int_{\CO} \big( \mu_\v - \frac 12 S(\gamma_2)\sigma_v^2  \big)\, dx \,dt +  \int_{\CO} \sigma_{\v} \,dW_2(t, x).
\end{split}
\end{equation}
By integration by parts we infer that
\begin{align*}
-r_{\v} \int_{\CO} \xi(t, x) \Delta \v(t, x) \,dx = r_{\v} \int_{\CO} \nabla \xi(t, x) \cdot \nabla \v(t, x)\,dx = -r_{\v} \int_{\CO} \xi^2(t, x) |\nabla \v(t, x)|^2\,dx.
\end{align*}
Using this in \eqref{lnxi-est}, rearranging and rewriting in integral form we get
\begin{equation}\label{lnxi-est-1}
\begin{split}
&\int_{\CO} \ln \xi(t, x) \,dx + r_{\v} \int_0^t \int_{\CO} \xi^2(s, x) |\nabla \v(s, x)|^2\,dx\,ds + \kappa _v\int_0^t \int_{\CO} \xi(s, x) \chi^2(s, x)\,dx\,ds
\\ &\leq \int_{\CO} \ln \xi_0(x) \,dx + \int_0^t \int_{\CO} \sigma_{\v} \,dW_2(t, x) + T \big( \mu_\v - \frac 12 S(\gamma_2)\sigma_v^2  \big)\,.
\end{split}
\end{equation}
Observe, firstly that we have for any $n\in\NN$ and $x>0$
 \DEQS
  - \ln(1/x)
\le   1_{(1,\infty)}(x)\ln(x)\le  (x-1)\le x.
 \EEQS
Straight forward calculations using Burkholder-Davis-Gundy inequality and taking expectation  we obtain
\begin{equation}\label{lnxi-est-2}
\begin{split}
&\bE \Big( \int_{\CO} |\ln \xi(s, x)\,dx| \Big) + r_{\v} \bE \Big(\int_0^t \int_{\CO} \xi^2(s, x) |\nabla \v(s, x)|^2\,dx\,ds \Big)
\\ &\quad+ \kappa_v \bE \Big( \int_0^t \int_{\CO} \chi^2(s, x) \,\xi(s, x) \,dx\,ds \Big)
\\
&\leq
\bE|v_0|_{L^1}+ \bE \Big( \int_{\CO} |\ln \xi_0(x)| \,dx \Big) + T C \big( \mu_\v , S(\gamma_2),\sigma_v  \big).
\end{split}
\end{equation}
\del{
Now squaring both sides and taking expectation we obtain
\begin{equation}\label{lnxi-est-2}
\begin{split}
&\bE \Big( \int_{\CO} \ln \xi(s, x)\,dx \Big) + r_{\v} \bE \Big(\int_0^t \int_{\CO} \xi^2(s, x) |\nabla \v(s, x)|^2\,dx\,ds \Big)^2
\\ &\quad+ \kappa_v \bE \Big( \int_0^t \int_{\CO} \chi^2(s, x) \,\xi(s, x) \,dx\,ds \Big)^2 \leq \bE \Big( \int_{\CO} \ln \xi_0(x) \,dx \Big)^2 + T C \big( \mu_\v , S(\gamma_2),\sigma_v  \big).
\end{split}
\end{equation}}
%
\end{proof}

\subsection{Properties of the system \eqref{eqv2sto}}\label{PoS-u}

Given the couple $(\chi, \eta)\in\mU$
 we investigate in this Section the existence and uniqueness of a solution to 
the auxiliary system
\begin{align}\label{gurv}
\begin{cases}
d u(t) &= \big[ r_u \Delta u(t) + \kappa _u\frac{\chi^2(t)}{\v(t)} - \Upsilon_u u(t) \big] dt + \sigma_u u(t) \,dW_1(t),
\\ u(0) &= u_0.
\end{cases}
\end{align}
Note that $v$ is a solution to system \eqref{eqv1sto}.
 Secondly, we investigate the regularity of $u$.

\begin{theorem}\label{sol_w}
Suppose $K_1, K_2>0$ and $K_3>0$ are fixed.
Then for any  $(\chi,\eta)\in\mU$ and $\v$ being a solution to \eqref{eq-vd}
 there exists a unique solution $u$ to system \eqref{gurv}
such that
\[
u \in C([0, T]; H^{1-\rho}_2(\CO)) \cap L^2(0, T;H^{2-\rho}_2(\CO)) \quad \mathbb{P}-a.s.
\]
Moreover, the solution to the system \eqref{gurv} is non-negative provided $u_0(x) \geq 0$ for all $x \in \CO$ and $\v(t, x) > 0,$ for $(t, x) \in [0, T] \times \CO.$
\end{theorem}

\begin{proof}
Since $[0,T]\ni t\mapsto \frac{\chi^2(t)}{\v(t)}$  belongs to $L^1(\Omega\times [0,T];L^1(\CO))$, we know by standard arguments that a solution exists.
In particular, we will apply \cite[Theorem 4.2.4]{weiroeckner} to prove the existence and uniqueness of the system \eqref{gurv}.
Similarly, one can use the setting of the book of Da Prato and Zabczyk \cite[Section 6.3.3]{prato} or use Theorem 4.5 \cite{vanneerven}.
Let us  provide the basic steps of the proof for our case.
Starting with the Gelfand triplet
\[ V' = H^{ -\rho}_2(\CO), \ \mathcal{H} =  H^{1-\rho}_2(\CO) \ \mbox{ and } \ V=H^{2-\rho}_2(\CO),
\]
we define 
$$\mathcal{A}(t): V \to V' \ \mbox{by} \ \mathcal{A} u =\Delta u + \frac{\chi^2(t)}{\v(t)},$$
 where
\[Dom(\Delta)= \big\{ u \in V': \mbox{Neumann  boundary conditions are satisfied} \big\}.\]
By  straight forward calculations one can see that this operator satisfies the Hypotheses of \cite[Theorem 4.2.4]{weiroeckner}. In fact, (H1) or the hemicontinuity, respective (H2) or  the weak monotonicity,  are given by the properties of $\Delta$.
%
%
To prove the coercivity i.e., (H3), 
we consider first
\begin{align}\label{Au_sig}
\langle \mathcal{A}(t) u, u\rangle + |\sigma_u u|^2_{\mathcal{H}} \leq C_1 |u|_{\mathcal{H}}^2 - |u|^2_{V} + \langle G(t), u \rangle,
\end{align}
with $G(t)=\frac{\chi^2(t)}{\v(t)}$.
Duality  for  $\theta=\frac \delta2$, $\delta=\rho-1$,
 the Sobolev embedding $L^1(\CO) \hookrightarrow H_2^{-1-\delta/2}(\CO)$, and the Young inequality for $\epsilon >0$
\begin{align}\label{in_thet}
\langle G(t), u\rangle \leq |G(t)|_{H_2^{-1-\delta/2}} |u|_{H_2^{1-3\delta/2}} \leq |G|_{L^1} |u|^{\theta}_{\mathcal{H}} |u|^{1-\theta}_{V} \leq C_{\epsilon} |G(t)|^2_{L^1} + c_{\epsilon} |u|^{2}_{\mathcal{H}} + \epsilon |u|^{2}_{V}.
\end{align}
Using \eqref{in_thet} in the inequality \eqref{Au_sig} and choosing $\epsilon = \frac{1}{2}$ we obtain,
\begin{align*}
\langle \mathcal{A}(t) u, u\rangle + |\sigma_u u|^2_{\mathcal{H}} \leq C_2 |u|_{\mathcal{H}}^2 - \frac{1}{2} |u|^2_{V} + C |G(t)|^2_{L^1}.
\end{align*}
Since $(\chi,\eta)\in \mU$, it follows by Proposition \ref{pr-xi_est}-(c) that $\chi^2\xi={\chi^2}/{\v} \in L^1((0, T) \times \Omega; L^1(\CO))$.
So, setting $f(t)=|G(t)|_{L^1}^2$ for all $t \in [0, T]$ we obtain the Hypothesis (H3).
The boundedness, i.e.\ (H4) follows since
\begin{align}
\langle \mathcal{A}(t) u, w\rangle \le   |u|_{V} |w|_V+  |G|_{L^1}|w|_V.
\end{align}
The non-negativity of the solution can be proved following Section 2.6 in \cite{Barbu+DaPrato+Rockner_2008} or  \cite{BDPR2016}.

\end{proof}

\subsubsection{{\bf Uniform bounds on $u$}}\label{SUbu}

In this subsection, we show several uniform bounds on $u$ under some assumptions.

\begin{proposition}\label{prop74}\SC
For any $K_1, K_2>0$ and $K_3>0,$
there exists
 constants  $\Con,\Cdelta>0 $ 
 and $\tilde{\Con}>0$  such that
 for any $T>0$, for any initial conditions $(u_0, \v_0)$ satisfying Hypothesis \ref{init}, for all
  $(\chi,\eta)\in\mU$, $\v$ being a solution to \eqref{eq-vd},
and $u$ being a solution to system \eqref{gurv}
\begin{equation}\label{u^2--Ito}
\begin{split}
 &\bE \sup_{0 \leq s \leq t} |u(s)|_{L^2}^{2} + 4 \,r_u\, \bE \int_0^t |\nabla u(s)|_{L^2}^{2} \,ds + \tilde{\Con} \,\bE \int_0^t |u(s)|_{L^2}^{2}\, ds
\\ &\leq \Con\, e^{\Cdelta T}\Big[ \bE |u_0|_{L^2}^{2} + C\, \bE \int_0^t \int_{\CO} u(s, x) \chi^{2}(s, x) \xi(s, x) \,dx\,ds \Big].
\end{split}
\end{equation}

\end{proposition}

\begin{proof}
For clarity, we write the calculation without defining a stopping time $\tau_k := \inf\{ t>0 \ | \ |u(t)|^2_{L^2} > k\}$ as done in  Proposition \ref{pr-xi_est}-(a).
Tracing the proof one can see easily that taking into account the stopping time and being precise will not change the result.
Let us put $\psi_1(u)=|u|_{L^p}^{p}.$ For $z, g, g_1, g_2 \in L^2(\CO)$ we infer that
\[
D \psi_1(z)[g] =  p \int_{\CO} |z(x)|^{p-1} \,g(x)\,dx, \quad D^2 \psi_1(z)[g_1, g_2] = p (p-1) \int_{\CO} |z(x)|^{p-2} g_1(x) g_2(x)\,dx.
\]
Applying the It\^o formula to the function $\psi_1(u)$ for $p=2$ we obtain,
\begin{equation}\label{u^2Ito}
\begin{split}
d|u(t)|_{L^2}^{2}  &= \Big[ 2r_u \int_{\CO} u(t, x) \cdot \Delta u(t, x) \,dx + 2 \kappa  _u\int_{\CO} u(t, x) \frac{\chi^{2}(t, x)}{\v(t, x)} \,dx
\\ &\quad - 2 \mu_u\,  | u(t)|^2_{L^2} + \sigma_u \int_{\mathcal{O}} u(t) [(\mbox{Id}+A)^{-\gamma_1} u(t)]\,dx + \frac{1}{2} S(\gamma_1)\,\sigma^2_{u} \,|u(t)|_{L^2}^{2} \Big] \,dt
\\ &\quad+ 2 \sigma_{u}  \big\langle u(t), u(t)\,dW_1(t) \big\rangle_{L^2},
\end{split}
\end{equation}
where we calculated the trace term as before. Using integration by parts we get,
\begin{align*}
2r_u \int_{\CO} u(t, x) \cdot \Delta u(t, x) \,dx = -2r_u \int_{\CO} (\nabla u(t, x))^2 \,dx.
\end{align*}
We may bound the term
\begin{align*}
\sigma_u \int_{\mathcal{O}} u(t) [(\mbox{Id}+A)^{-\gamma_1} u(t)]\,dx &\leq \sigma_u |u(t)|_{L^2}\, |(\mbox{Id}+A)^{-\gamma_1} u(t)|_{L^2}
\\ &\leq \sigma_u\, \|(\mbox{Id}+A)^{-\gamma_1}\|_{\mathcal{L}(L^2)}\, |u(t)|_{L^2}^2 \leq C(\sigma_u) \, |u(t)|_{L^2}^2\, .
\end{align*}
Using the above estimates, we rewrite \eqref{u^2Ito} in the integral form as follows:
\begin{equation}\label{u^2Ito-3}
\begin{split}
 &|u(t)|_{L^2}^{2} + 2r_u \int_0^t |\nabla u(s)|_{L^2}^{2} \,ds + 2 \mu_u \int_0^t |u(s)|_{L^2}^{2}\,ds
\\ &\leq |u_0|_{L^2}^{2} + C\, \int_0^t \int_{\CO} u(s, x) \chi^{2}(s, x) \xi(s, x) \,dx\,ds
\\ &\quad{}+ {C(S(\gamma_1), \sigma_{u}) \int_0^t |u(s)|_{L^2}^{2} \,ds + 2 \sigma_{u} \Big| \int_0^t \int_{\CO} u^2(s, x) \,dW_1(s, x) \Big|},
\end{split}
\end{equation}
where $\xi=v^{-1}$ and the constant $C>0$ depends on $\kappa _u$.
Now taking supremum over $s \in [0, t]$ and then taking expectation we infer,
\begin{equation}\label{u2Ito_3}
\begin{split}
 &\bE \sup_{0 \leq s \leq t}  |u(s)|_{L^2}^{2} + 2r_u \bE \int_0^t |\nabla u(s)|_{L^2}^{2} \,ds + 2 \mu_u \bE \int_0^t |u(s)|_{L^2}^{2} \,ds
\\
 &\leq \bE |u_0|_{L^2}^{2} + C\, \bE \int_0^t \int_{\CO} u(s, x) \chi^{2}(s, x) \xi(s, x) \,dx\,ds
\\
 &\quad{}+ C(S(\gamma_1), \sigma_{u}) \,\bE \int_0^t |u(s)|_{L^2}^{2} \,ds + \bE \sup_{0 \leq s \leq t} \Big| \int_0^t \int_{\CO} u^2(s, x)\,dW_1(s, x) \Big|.
\end{split}
\end{equation}
Applying the Burkholder-Davis-Gundy inequality and following similar steps as in \eqref{mbbExi-2} we infer,
\begin{equation}\label{u2Ito_-3}
\begin{split}
&\bE \sup_{0 \leq s \leq t} \Big| \int_0^t \int_{\CO} u^2(s, x)\,dW_1(s, x) \Big| \leq \frac{\epsilon}{2} \,\bE \sup_{0 \leq s \leq t}  |u(s)|_{L^2}^{2} + \frac{\sigma_u}{2 \epsilon} \,\bE \int_0^t |u(s)|_{L^2}^{2} \,ds,
\end{split}
\end{equation}
where we have applied the Young inequality for $\epsilon>0.$
Choosing $\epsilon=\frac 12 $ in \eqref{u2Ito_-3} and putting in \eqref{u2Ito_3} and rearranging we obtain
\begin{equation}\label{u^2--Ito2}
\begin{split}
 &\bE \sup_{0 \leq s \leq t} |u(s)|_{L^2}^{2} + 4r_u \bE \int_0^t |\nabla u(s)|_{L^2}^{2} \,ds + 2 \mu_u\, \bE \int_0^t |u(s)|_{L^2}^{2} \,ds
\\ &\leq 2 \bE |u_0|_{L^2}^{2} + C\, \bE \int_0^t \int_{\CO} u(s, x) \chi^{2}(s, x) \xi(s, x) \,dx\,ds + C(S(\gamma_1),\sigma_u)\, \bE \int_0^t |u(s)|_{L^2}^{2} \,ds.
\end{split}
\end{equation}
 Dropping the second and third term from the left hand side and by the application of Gronwall's lemma, there exist constants $\delta_4$ and $C_4$ such that we get the assertion.

\end{proof}

\begin{proposition}\label{lem-xi,u}\SC
For any $K_1, K_2>0$ and  $K_3>0,$
there exist constants  $\Con,\Cdelta>0 $
such that  for any $\epsilon >0$, for any $T>0$, for any initial conditions $(u_0, \v_0)$ satisfying Hypothesis \ref{init},  for any   $(\chi,\eta)\in\mU$, $\v$ being a solution to \eqref{eq-vd},
and $u$ being a solution to system \eqref{gurv} we have for any $t\in[0,T]$
\DEQS 
\lqq{  \mathbb{E}  \big( \ln(\xi(t)) \cdot u(t) \big) +
\mu_u \bE \int_0^ t  \big(\ln(\xi(s)) \cdot u(s) \big)\,ds}
&&\\
&&{}+r_{\v} \,\bE \int_0^ t  \int_{\CO} \xi^2(s, x) u(s, x) \cdot |\nabla \v(s, x)|^2 \,dx\, ds + \kappa_u\,\bE \int_0^ t  \int_{\CO} \chi^2(s, x) \xi(s, x) u(s, x) \,dx\, ds
\\
&\le &
 \mathbb{E} \big(\ln( \xi_0) \cdot u_0 \big)   + \epsilon \,\EE \int_0^ t |\nabla u(s)|^2_{L^2}\, ds +
\frac {(r_u+r_{\v})^2}{4\epsilon} \EE \int_0^ t |\xi(s)\nabla v(s)|_{L^2}^2\, ds
\nonumber
\\ &&{}
 +\kappa _u\, \EE \int_0^ t\ln(\xi(s, x)) \chi^2(s, x) \xi(s, x) \,dx \,ds +  \mu_{\v}\,\EE \int_0^ t \int_{\CO} \xi(s,x)\,u(s, x) \,dx\, ds
 \\ &&{}
  +C(S(\gamma_1), \sigma_u) \,\EE \int_0^ t \big|\ln(\xi(s)) \cdot u(s) \big|\, ds + \sigma_u \EE \int_0^ t \int_\CO \ln(\xi(s,x)) \,(\mbox{Id}+A)^{-\gamma_1} u(s,x) \,dx\, ds.
\EEQS

\end{proposition}

\begin{proof}
We will apply the It\^o formula to $\phi(t)= \ln(\xi(t)) \cdot u(t)$.
For clarity, we write the calculation without defining a stopping time $\tau_k := \inf\{ t>0 \mid \big( \ln(\xi(t)) \cdot u(t)\big) \geq k\}$ as done in  Proposition \ref{pr-xi_est}-(a).
Tracing the proof one can see easily that taking into account the stopping time and being precise will not change the result.

Before applying the  It\^o formula, let us note that using integration by parts and the H\"older inequality leads to 
\begin{equation*}
\begin{split}
\int_{\CO} \ln (\xi(s,x))  \Delta u(s,x) \,dx &= -  \int_{\CO} (\nabla \ln( \xi(s,x))) \cdot (\nabla u(s,x)) \,dx 
\\ &=  \int_{\CO} \xi(s,x) (\nabla \v(s,x)) \cdot (\nabla u(s,x)) \,dx
\\ &\leq
 \Big( \int_{\CO} \xi^2 (s,x)|\nabla \v(s,x)|^2 \,dx \Big)^{1/2} \Big( \int_{\CO} |\nabla u(s,x)|^2 \,dx \Big)^{1/2}.
\end{split}
\end{equation*}
By the Young inequality, we get for any $\epsilon>0$
\DEQSZ\label{est_oben}
r_{u} \int_{\CO} \ln (\xi(s,x))  \Delta u(s,x) \,dx &\le &
 \frac {r_{u}^2}{4\epsilon} \int_{\CO} \xi^2 (s,x) |\nabla \v(s,x)|^2 \,dx +\epsilon  \int_{\CO} |\nabla u(s,x)|^2 \,dx.
 \EEQSZ
Applying  integration by parts leads to
\DEQS
\lqq{ -\int_{\CO}\xi(s,x))  \Delta \v(s,x)\, u(s,x) \,dx =   \int_{\CO} \nabla(\xi(s,x) u(s,x) ) \cdot (\nabla \v(s,x)) \,dx }&&
\\
& =&
\int_\CO   \xi(s,x)  (\nabla u(s,x) ) \cdot (\nabla \v(s,x)) \,dx + \int_{\CO} \nabla \xi(s,x)\cdot \nabla \v(s,x) u(s,x) \, dx
\\
&=&\int_\CO   \xi(s,x)  (\nabla u(s,x) ) \cdot (\nabla \v(s,x)) \,dx
- \int_{\CO}  \xi^2(s,x) |\nabla \v(s,x)|^2 u(s,x) \, dx .
\EEQS
%
%
Secondly, that we have by the It\^o-formula 
\begin{equation*}
\begin{split}
& d\ln(\xi(s,x)) \\ &= -\xi(s,x)\big[ r_\v\Delta \v(s, x) + \kappa _u \chi^2(s, x) +\sigma_v \,(\mbox{Id}+A)^{-\gamma_2}\,v(s, x) \big]\, ds
\\ &\quad+ \mu_{\v}  \xi(s,x)\, ds + \sigma_v \sum_{k=1}^\infty(1+\lambda_k)^{-\gamma_2} \la v(s)\xi^2v(s) e_k,e_k\ra d\beta^2_k(s) - \frac {\sigma_v^2}2  \sum_{k=1}^\infty (1+\lambda_k)^{-2\gamma_2}\,\, ds.
\end{split}
\end{equation*}
As pointed out on page \pageref{equation_uv_strat}, $(\mbox{Id}+A)^{-\gamma_2}$ is a positive operator.
Using the positivity of $v$,  we drop the integral term $-\sigma_v\,\xi(s,x)  \,(\mbox{Id}+A)^{-\gamma_2}\,v(s, x) $
at  the right hand side.

\medskip

Now we are ready to calculate $\ln(\xi(t)) \cdot u(t) $. Applying the  It\^o-formula
 and taking expectation we obtain
\begin{equation}\label{ln.xi,u}
\begin{split}
&\EE \big( \ln(\xi(t)) \cdot u(t) \big) -\EE \big( \ln(\xi_0) \cdot u_0)\big)
\\ &\leq -\EE \int_0^ t\int_\CO  {\xi(s,x)} \Big[r_{\v} \Delta {\v}(s,x) +\kappa  _u\chi^2(s,x)\Big]u(s,x)\, dx\, ds
\\ &\quad+  \mu_{\v}\,\EE \int_0^ t\int_\CO \xi(s,x) u(s,x)\, dx\, ds -\frac 12 S(\gamma_1) \sigma^2_{\v} \,\bE \int_0^ t\int_\CO  u(s,x)\, dx\, ds
\\ &\quad+\EE \int_0^ t\int_\CO \ln(\xi(s,x))\Big[ r_u\Delta u(s,x)+\kappa  _u\chi^2(s,x)\xi(s,x)
\\ &\qquad\qquad\qquad\qquad\qquad\quad- \mu_u u(s,x) + \sigma_u \,(\mbox{Id}+A)^{-\gamma_1} u(s,x) \Big]\, dx\, ds.
\end{split}
\end{equation}
Using the previous estimates we get,
\begin{equation}\label{log.xi,u}
\begin{split}
&\EE \big( \ln(\xi(t)) \cdot u(t) \big) -\EE \big( \ln(\xi_0) \cdot u_0 \big)
\\ &\le -r_{\v} \EE \int_0^ t  \int_{\CO} \xi^2(s, x) u(s, x) |\nabla \v(s, x)|^2 \,dx\, ds - \kappa  _u\EE \int_0^ t  \int_{\CO} \chi^2(s, x) \xi(s, x) u(s, x) \,dx\, ds
\\  &\quad+ \mu_{\v}\EE \int_0^ t \int_{\CO} \xi(s,x)\,u(s, x) \,dx\, ds +\frac {(r_u+r_{\v})^2}{4 \epsilon}\, \bE \int_0^t \int_{\CO} \xi^2(s, x) |\nabla \v(s, x)|^2 \,dx\,ds
 \\ &\quad + \epsilon\, \bE \int_0^t |\nabla u(s)|_{L^2}^2 \,ds +\kappa  _u\,\EE \int_0^ t \int_\CO \ln(\xi(s,x))\chi^2(s,x)\xi(s,x)\, dx \,  ds
\\ &\quad-\mu_u \EE \int_0^ t \int_\CO \ln(\xi(s,x))u(s,x)\, dx \, ds + \sigma_u \EE \int_0^ t \int_\CO \ln(\xi(s,x)) \,(\mbox{Id}+A)^{-\gamma_1} u(s,x) \,dx\, ds.
\end{split}
\end{equation}
Rearranging the terms we finally get 
\DEQS
\lqq{\EE \big( \ln(\xi(t)) \cdot u(t)\big) -\EE \big( \ln(\xi_0) \cdot u_0 \big) + \mu_u \EE \int_0^ t \big( \ln(\xi(s)) \cdot u(s)\big)\,ds}
&&\\
&&{}+ r_{\v} \EE \int_0^ t \int_{\CO} \xi^2(s, x) u(s, x) |\nabla \v(s, x)|^2 \,dx\, ds + \kappa _u \EE \int_0^ t  \int_{\CO} \chi^2(s, x) \xi(s, x) u(s, x) \,dx\, ds
\\
&\le &
 \mu_{\v} \EE \int_0^ t \int_{\CO} \xi(s,x)\,u(s, x) \,dx\, ds  +\frac 12 S(\gamma_1)\sigma^2_u \EE \int_0^ t \big( \ln(\xi(s)) \cdot u(s)\big)\, ds
 \\
&&{} +
\kappa _v \EE \int_0^ t \int_\CO \ln(\xi(s,x)) \chi^2(s,x)\xi(s,x)\, dx \, ds +
 \frac {(r_u+r_{\v})^2}{4\epsilon}\,  \bE \int_0^t \int_{\CO} \xi^2(s, x) |\nabla \v(s, x)|^2 \,dx\,ds
\\
&&{} + \epsilon\, \bE \int_0^t |\nabla u(s)|_{L^2}^2 \,ds + \sigma_u \EE \int_0^ t \int_\CO \ln(\xi(s,x)) \,(\mbox{Id}+A)^{-\gamma_1} u(s,x) \,dx\, ds \, .
\EEQS

%
Thus, we get the assertion.

\end{proof}

\begin{corollary}\label{cor01}
For any $K_1, K_2>0$ and  $K_3>0$,
\begin{enumerate}
\item[$(a)$]\SC
and any $n\in\NN$ there exist a constant  $C_4>0 $
such that  for any $\epsilon>0$, for any $T>0$, for any initial conditions $(u_0, \v_0)$ satisfying Hypothesis \ref{init},
 for any  $(\chi,\eta)\in\mU$, $\v$ being a solution to \eqref{eq-vd},
and $u$ being a solution to system \eqref{gurv}
we have
\begin{equation*}
\begin{split}
&- \EE\int_\CO \ln(\xi(t,x))\, u(t,x)\, dx
\\ &\le C_4
\frac n {4\ep}\lk\{\EE\sup_{0\le s\le T}|\eta(s)|_{L^2}^2 +  \EE\int_0^ T |\eta(s)|_{H^1_2}^2\, ds\rk\}^\frac 1n+ \ep  \EE |u(t)|_{L^2}^2.
\end{split}
\end{equation*}

\item[$(b)$] \SC
there exist constants  $c_1,\ldots,c_{5}>0 $
such that  for any $\epsilon>0$,  for any $T>0$, for any initial conditions $(u_0,v_0)$ satisfying Hypothesis \ref{init},
 for any  $(\chi,\eta)\in\mU$, $\v$ being a solution to \eqref{eq-vd},
and $u$ being a solution to system \eqref{gurv}
we have
\del{\begin{equation*}
\begin{split}
& \EE \int_0^ t  \int_{\CO} \chi^2(s, x) \xi(s, x) u(s, x) \,dx\, ds
\\
& \le  
\tilde c_4 e^{\Cdz T} \Bigg\{ \lk\{ \EE |v_0|_{L^1}\rk\}^\frac 12\EE|u_0|_{L^2}^2
\\&{}+ \frac 4{4\ep_0}\lk\{ \EE\sup_{0\le s\le T}|\eta(s)|_{L^2}^2 +  \EE\int_0^ T |\eta(s)|_{H^1_2}^2\, ds\rk\}^\frac 14
\\
&{} +\Big(\frac { \sigma^2_u}{4\ep_2} + \frac 1{4\epsilon_3} \Big) \int_0^ t |\xi(s)|_{L^2}^2\, ds
+c_5 e^{\delta_5 T} \Big( 1+ \EE |\xi_0|_{L^p}^p\Big)\Bigg\}
 \\&{}+\ep\lk( \EE \sup_{0\le s\le T}|u(s)|_{L^2}^2+ \EE\int_0^ T|u(s)|_{L^2}^2\, ds\rk)
.
\end{split}
\end{equation*}}
\DEQS 
\lqq{\bE \int_0^ t  \int_{\CO} \chi^2(s, x) \xi(s, x) u(s, x) \,dx\, ds
} &&
\\
& \le & c_1\EE|\xi_0|_{L^1}+c_2\EE|u_0|^2_{L^2} +c_3\bE|v_0|_{L^1}
 +c_4 \EE \int_0^ t |\xi(s)|_{L^2}^2\, ds
\\&&{}+c_5\bigg( \EE\sup_{0\le s\le T}|\eta(s)|_{L^2}^2 +  \EE\int_0^ T |\eta(s)|_{H^1_2}^2\, ds\bigg)^\frac 14
 \\&&{}+\ep\lk( \EE \sup_{0\le s\le T}|u(s)|_{L^2}^2+ \EE\int_0^ T|u(s)|_{L^2}^2\, ds\rk).
\EEQS

\end{enumerate}

\end{corollary}

\begin{proof}
Let us start to show item (a).
Before starting, observe, firstly that we have for any $n\in\NN$ and $x>0$
 \DEQS
  - \ln(1/x)\le n\ln(x^{-\frac 1n})
\le  n( 1_{(1,\infty)}(x)\ln(x^{-\frac 1n})\le n x^\frac 1n
 \EEQS
 and, secondly,
 \DEQSZ\label{iqeins}
\lqq{- \EE\int_\CO \ln(\xi(t,x)u(t,x))\, dx \le n\EE\int_\CO v^{\frac 1n}(t,x) u(t,x)}
\nonumber
\\
&\le& n\lk\{ \EE \int_\CO  v^\frac 2n(t,x)\, dx\rk\}^\frac 12 \lk\{ \EE |u(t)|_{L^2}^2\rk\}^\frac 12
\nonumber
\\
&\le& \frac n {4\ep}\EE \int_\CO  v^\frac 2n(t,x)\, dx+\ep  \EE |u(t)|_{L^2}^2
\le
\frac n {4\ep}\lk\{ \EE \int_\CO  v(t,x)\, dx\rk\}^\frac 1n+ \ep  \EE |u(t)|_{L^2}^2
.
 \EEQSZ
 By Theorem \ref{theou1}, we know that there exists constants $C_1,C_2>0$ such that
 \DEQSZ\label{iqzwei}
 \EE |v(t)|_{L^1} &\le& C_1\EE |v(t)|_{L^2}\le C_2\lk(\EE\sup_{0\le s\le T}|\eta(s)|_{L^2}^2 +  \EE\int_0^ T |\eta(s)|_{H^1_2}^2\, ds\rk).
 \EEQSZ
The assertion follows by the Jensen inequality.
\SCz

\medskip

Now we will show (b). By Proposition \ref{lem-xi,u} we know that for any $\ep_1,\ep_2>0$ we have
\DEQSZ\label{prop7.5}
\lqq{ \qquad  \mathbb{E}  \big( \ln(\xi(t)) \cdot u(t) \big) }
&&\\\nonumber
&&
{}+
\mu_u \bE \int_0^ t  \big(\ln(\xi(s)) \cdot u(s) \big)\,ds+
\bE \int_0^ t  \int_{\CO} \chi^2(s, x) \xi(s, x) u(s, x) \,dx\, ds 
\\
\nonumber
&\le &
 \mathbb{E} \int_\CO\ln( \xi_0(x))  u_0(x)\, dx    + \epsilon_1 \,\EE \int_0^ t |\nabla u(s)|^2_{L^2}\, ds +
\frac {(r_u+r_{\v})^2}{4\epsilon_1} \EE \int_0^ t |\xi(s)\nabla v(s)|_{L^2}^2\, ds
\nonumber
\\ &&{}
\nonumber
 +\kappa _u\, |\EE \int_0^ t\ln(\xi(s, x)) \chi^2(s, x) \xi(s, x) \,dx \,ds|
 \\
 && {}+ C(S(\gamma_1), \sigma_u) \,\EE \int_0^ t \int_\CO\big( \lk| \ln(\xi(s,x))\rk| \cdot u(s,x) \, dx\big)\, ds
\nonumber
 \\ &&{}
 +  \mu_{\v}\,\EE \int_0^ t \int_{\CO} \xi(s,x)\,u(s, x) \,dx\, ds +
\frac { \sigma^2_u}{4\ep_2} \EE \int_0^ t |\xi(s)|_{L^2}^2 \, ds + \ep_2  \EE \int_0^ t|u(s)|_{L^2}^ 2\, ds .
\nonumber
\EEQSZ
First, note due to estimate \eqref{iqeins} and \eqref{iqzwei} the negative part of
$$ \mathbb{E} \big(\ln( \xi(t)) \cdot u(t) \big)
$$
and the negative part of
$$ \EE \int_0^ t \int_\CO\big( \ln(\xi(s,x)) \cdot u(s,x) \, dx\,ds
$$
can be estimated by
\DEQSZ\label{ie01}
\\
\nonumber
\lk(\frac n{4\ep_0} +\frac n{4\tilde \ep_0}\rk)\lk\{\EE\sup_{0\le s\le T}|\eta(s)|_{L^2}^2 +  \EE\int_0^ T |\eta(s)|_{H^1_2}^2\, ds\rk\}^\frac 1n+\ep_0\int_0^ T\EE|u(s)|_{L^2}^2\, ds+\tilde \ep_0\EE|u(t)|_{L^2}^2,
\EEQSZ
where $\ep_0$ and $\tilde \ep_0$ are arbitrary.
By Proposition \ref{pr-xi_est}-(c) the term $ \EE \int_0^ t |\xi(s)\nabla v(s)|_{L^2}^2\, ds$
can be bounded by
$$
\bE|v_0|_{L^1}+ C_1\, \bE \Big( \int_{\CO} |\ln \xi_0(x)| \,dx \Big) +  C_2 T.
$$
By Proposition \ref{pr-xi_est}-(c) the term
$\EE \int_0^ t\xi^2(s, x) \chi^2(s, x)  \,dx \,ds$
can be bounded by $Ce^{\delta}T\EE|\xi_0|$.
Since there exists a constant $C>0$ such that $-C\le \ln(x)x\le x^2$, the term $\EE \int_0^ t\ln(\xi(s, x)) \chi^2(s, x) \xi(s, x) \,dx \,ds$
can be estimated by  $\EE \int_0^ t\xi^2(s, x) \chi^2(s, x)  \,dx \,ds$.
In this way there exists constants $C_3,\delta_3>0$ such that
\DEQS
\lqq{ \frac {(r_u+r_{\v})^2}{4\epsilon_1} \EE \int_0^ t |\xi(s)\nabla v(s)|_{L^2}^2\, ds}
&&
\\&&{}
 +\kappa _u\, \EE \int_0^ t\ln(\xi(s, x)) \chi^2(s, x) \xi(s, x) \,dx \,ds
 \le C_3e^{\delta_3T}\EE|\xi_0|_{L^1}.
\EEQS
Let us consider the term $\EE \int_0^ t \int_{\CO} \xi(s,x)\,u(s, x) \,dx\, ds$
and the positive part of $\int_0^ t \int_{\CO} \ln(\xi(s,x))\,u(s, x) \,dx\, ds$.
Firstly, since $\ln(x)\le x$ for $x\ge 1$, the positive part given by
$\EE \int_0^ t \int_{\CO} 1_{[1,\infty)}(\xi(s,x))\ln(\xi(s,x))\,u(s, x) \,dx\, ds$
can be estimated by $\EE \int_0^ t \int_{\CO} \xi(s,x)\,u(s, x) \,dx\, ds$.
Applying the Young inequality we know that we have for any $\ep_3>0$
$$
\EE \int_0^ t \int_{\CO} \xi(s,x)\,u(s, x) \,dx\, ds\le {\epsilon_3} \int_0^ t |u(s)|_{L^2}^2\, ds + \frac 1{4\epsilon_3}  \int_0^ t |\xi(s)|_{L^2}^2\, ds.
$$
Note, we have for $x\ge 1$,
$\ln(x)\le 2\ln(x^\frac 12)\le (x-1)^\frac 12\le x^\frac 12+1$.
The Cauchy Schwarz and Young inequality  gives for the initial condition
$$
\EE\int_\CO1_{[1,\infty)}(\xi_0(x))\,\ln(\xi_0(x))\,u_0(x)\, dx\le C_1\EE|\xi_0|_{L^1}+C_2\EE|u_0|^2_{L^2}.
$$%
Collecting altogether we get by estimate \eqref{prop7.5} and estimates \eqref{ie01}, we can conclude that there exist constants
$c_0,\ldots,c_{10},\delta_6>0$ such that we have for all $\ep_0,\tilde \ep_0,\ep_1,\ep_2,\ep_3>0$
\DEQS
\lqq{\bE \int_0^ t  \int_{\CO} \chi^2(s, x) \xi(s, x) u(s, x) \,dx\, ds\le  c_1\EE|\xi_0|_{L^1}+c_2\EE|u_0|^2_{L^2} }
&&\nonumber
\\
\nonumber
&  &{}+c_3\bE|v_0|_{L^1}+ c_4\, \bE \Big( \int_{\CO} |\ln \xi_0(x)| \,dx \Big) +  c_5 T
+ c_6 e^{\delta_6 T} \EE|\xi|_{L^1}
\\
&&{}+c_7 {(\ep_1+\ep_2+\tilde \ep_0)} \sup_{0\le s\le T}|u(s)|_{L^2}^2+
c_8(\ep_0+\epsilon_3) \EE\int_0^ T|u(s)|_{L^2}^2\, ds
\\&&{}+c_9\lk( \frac n{4\ep_0}+ \frac n{4\tilde \ep_0}\rk)\lk\{ \EE\sup_{0\le s\le T}|\eta(s)|_{L^2}^2 +  c_4\EE\int_0^ T |\eta(s)|_{H^1_2}^2\, ds\rk\}^\frac 1n
\\
&&{} +c_{10}\Big(\frac { \sigma^2_u}{4\ep_2} + \frac 1{4\epsilon_3} \Big) \int_0^ t |\xi(s)|_{L^2}^2\, ds.
\EEQS
%
%
%
Summarizing, for all $\ep>0$ gives constants
\DEQS 
\lqq{\lk(\bE \int_0^ t  \int_{\CO} \chi^2(s, x) \xi(s, x) u(s, x) \,dx\, ds
\rk)
} &&
\\
& \le & c_1\EE|\xi_0|_{L^1}+c_2\EE|u_0|^2_{L^2} +c_3\bE|v_0|_{L^1}
\\&&{}+C_3(n)\bigg( \EE\sup_{0\le s\le T}|\eta(s)|_{L^2}^2 +  \EE\int_0^ T |\eta(s)|_{H^1_2}^2\, ds\bigg)^\frac 1n
\\
&&{} +C_4 \EE \int_0^ t |\xi(s)|_{L^2}^2\, ds
+C_5\Big( 1+ \EE |\xi_0|_{L^p}^p\Big)\Bigg\}
 \\&&{}+\ep\lk( \EE \sup_{0\le s\le T}|u(s)|_{L^2}^2+ \EE\int_0^ T|u(s)|_{L^2}^2\, ds\rk),
\EEQS
which leads by taking $n=4$ the assertion.
\end{proof}

\subsection{Proof of continuity of the operator $\mathscr{T}$}\label{cont}

In this section we state all results, that we are using to show that the operator
\DEQS
\mathscr{T}: \mU &\longrightarrow & \CM_{\mathfrak{A}}(0,T);
\\
(\chi,\eta) & \mapsto & (u,v),
\EEQS
 is continuous.
The following proposition shows the continuity of the operator  $\mathscr{T}$  with respect to $u$ in terms of $\chi$ and $\eta$.

\begin{proposition}\label{help_cont}
Let $K_i>0$ for $i=1,2, 3.$ For all  $(\chi_1, \eta_1),(\chi_2, \eta_2)\in\mU $ and $(u_1,\v_1),(u_2,\v_2)\in\mU $ such that  $(u_1,\v_1)=\mathscr{T}[ (\chi_1, \eta_1)]$ and $(u_2,\v_2)=\mathscr{T}[(\chi_2, \eta_2)]$,
 there exists a constant $C=C(K_1,K_2,K_3)>0$, some numbers $\delta_1, \delta_2 >0$ and $\gamma\in(0,1)$ such that
{\small{
\begin{equation*}
\begin{split}
&\mathbb E \sup_{0 \leq s \leq T}  \big| u_1(s) - u_2(s) \big|_{H^{1-\rho}_2}^2
\\ &\leq C(K_1,K_2,T) \times \Big[ \Big\{ \mathbb{E} \big[ \sup_{s \in [0, T]} \big|\chi_1(s) - \chi_2(s) \big|^2_{H^{1-\rho}_2} \big] \Big\}^{\delta_1} + \Big\{ \mathbb{E} \big[ \sup_{s \in [0, T]} \big| \xi_1(s) - \xi_2(s) \big|_{L^1}^\gamma  \big] \Big\}^{\delta_2}\,\Big],
\end{split}
\end{equation*}
}}
where $\xi_i = {v^{-1}_i}$ for $i=1, 2.$
\end{proposition}

Before stating the proof,
we present the following proposition where  the continuity of $\xi$ on $v$ is stated.
Since the proof is short, we present the proof after stating the proposition.
\begin{proposition}\label{help_cont1}
Let $\xi_i= \v_i^{-1}$ for $i=1, 2$.
Then, for any $1\le q<2$ there exists a constant $C>0$ and  numbers  $p=\frac {2mq}{m-q}$, $m>q$, and $r=\frac 2{1-\gamma}$  such that
{\small{
\begin{align*}
\Big( \mathbb E \big[ \sup_{0\le s\le T}  \left| {\xi_1(s)} -{\xi_2(s)}\right|_{L^q}^\gamma \big] \Big)^\frac 1q  \leq  \Big( \EE  \big[ \sup_{0\le s\le T}  \left| {\v_1(s)} -{\v_2(s)}\right|_{L^m}  \big] \Big)^\frac 12  \Big( \EE |\xi_1(s)|_{L^p}^r\Big)^\frac 1r\, \Big(\EE| \xi_2(s)|_{L^p}^r \Big)^\frac 1r.
\end{align*}
}}
\end{proposition}
\begin{proof}[Proof of Proposition \ref{help_cont1}:]
This follows by
\DEQS
\xi_1(s)-\xi_2(s)=\frac {\v_2(s) - \v_1(s)}{\v_1(s) \v_2(s)} = \big(\v_2(s) - \v_1(s) \big)\,\xi_1(s)\,\xi_2(s),
\EEQS
and the H\"older inequality.
\end{proof}

\noindent
Finally, we present the following proposition, stating the continuity of $\v$ on $\chi$.
Again, since the proof is short, we present the proof after stating the proposition.
\begin{proposition}\label{help_cont4}
Let $K_i>0$ for $i=1,2, 3$. Let us assume that  $\rho\in[1,2)$ and $2\ge d(1-\frac 1q)$. For all  $(\chi_1, \eta_1),(\chi_2, \eta_2)\in\mU $ and $(u_1,\v_1),(u_2,\v_2)\in\mU $ such that  $(u_1,\v_1)=\mathscr{T}[ (\chi_1, \eta_1)]$ and $(u_2,\v_2)=\mathscr{T}[(\chi_2, \eta_2)]$,
 there exists a constant $C=C(K_1,K_2,K_3)>0$, some numbers $\delta_1, \delta_2 >0$ and $\gamma\in(0,1)$ such that
\begin{equation*}
\begin{split}
&\mathbb E \sup_{0 \leq s \leq T}  \big| v_1(s) - v_2(s) \big|_{L^q}
\leq C(K_3,T) \times \Big[ \mathbb{E}  \sup_{s \in [0, T]} \big|\chi_1(s) - \chi_2(s) \big|^2_{H^{1-\rho}_2}  \Big]^{\gamma }.
\end{split}
\end{equation*}
\end{proposition}
\begin{proof}[Proof of Proposition \ref{help_cont4}:]
We start with focusing on the non-linear part.
The Minkowski inequality yields
$$
\Big|\int_0^ t e^{-\sigma_v(t-s)\Delta}(\chi_1^2(s)-\chi_2^2(s))\, ds\Big|_{L^q} \le
\int_0^ t \Big|e^{-\sigma_v(t-s)\Delta}(\chi_1^2(s)-\chi_2^2(s))\Big|_{L^q}\, ds.
$$
The smoothing property for the Laplace operator $\Delta$ and the H\"older inequality give for $\frac \delta2<1-\frac 1r$
\DEQS
\lqq{
\lk|\int_0^ t e^{-\sigma_v(t-s)\Delta}(\chi_1^2(s)-\chi_2^2(s))\, ds\rk|_{L^q}
}
&&
\\
&\le &\int_0^ t (t-s)^{-\frac \delta2 }\lk|\chi_1^2(s)-\chi_2^2(s)\rk|_{H^{-\delta}_q}\, ds
\le C(t)\,
\lk\|\chi_1^2-\chi_2^2\rk\|_{L^r(0,t;H^{-\delta}_q)}
.
\EEQS
The embedding $L^1(\CO)\hookrightarrow H^{-\delta}_q(\CO)$ gives for $\delta\ge d\lk(1-\frac 1q\rk)$
\DEQS
\lk|\int_0^ t e^{-\sigma_v(t-s)\Delta}(\chi_1^2(s)-\chi_2^2(s))\, ds\rk|_{L^q}
\le C(t)\,
\lk\|\chi_1^2-\chi_2^2\rk\|_{L^r(0,t;L^1)}
.
\EEQS
The H\"older inequality gives
\DEQS
\lk|\int_0^ t e^{-\sigma_v(t-s)\Delta}(\chi_1^2(s)-\chi_2^2(s))\, ds\rk|_{L^q}
\le C(t)\,\lk(
\int_0^ t \lk|\chi_1(s)-\chi_2(s)\rk|^r_{L^{2}}\lk|\chi_1(s)+\chi_2(s)\rk|^r_{L^{2}}\, ds\rk)^\frac 1r
.
\EEQS
We know by real interpolation (see paragraph \ref{notation2})  that there exists a constant $C>0$ such that
$|w|_{L^2}\le C\, |w|_{H^{1-\rho}_2}^\theta |w|_{H^{2-\rho}_2}^{1-\theta}$ for $\theta=1/\rho$. 
Applying first interpolation, and, subsequently, the H\"older inequality gives
\DEQS
\lqq{ \lk|\int_0^ t e^{-\sigma_v(t-s)\Delta}(\chi_1^2(s)-\chi_2^2(s))\, ds\rk|_{L^q} }
&&
\\
&\le &C(t)\, \sup_{0\le s\le t} \lk|\chi_1(s)-\chi_2(s)\rk|_{H^{1-\rho}_2}^\theta \sup_{0\le s\le t} \lk|\chi_1(s)+\chi_2(s)\rk|_{H^{1-\rho}_2}^\theta
\,
\\
&&{}\times \lk(\int_0^ t \lk|\chi_1(s)-\chi_2(s)\rk|^{(1-\theta)r}_{H^{1}_2} \lk|\chi_1(s)+\chi_2(s)\rk|^{(1-\theta)r}_{H^{1}_2} \, ds\rk)^\frac 1r
.
\EEQS
Taking expectation and applying the Cauchy-Schwarz inequality we get for $r\le\frac 1{1-\theta}$
\DEQS
\lqq{ \EE \lk|\int_0^ t e^{-\sigma_v(t-s)\Delta}(\chi_1^2(s)-\chi_2^2(s))\, ds\rk|_{L^q} }
&&
\\
&\le &C(t)\,\lk\{ \EE \sup_{0\le s\le t} \lk|\chi_1(s)+\chi_2(s)\rk|_{H^{1-\rho}_2}^2\rk\} ^{\frac \theta2}
\lk\{ \EE  \sup_{0\le s\le t} \lk|\chi_1(s)-\chi_2(s)\rk|_{H^{1-\rho}_2}^2\rk\}^\frac \theta2
\\
&&{}\times \lk\{ \EE  \int_0^ t \lk|\chi_1(s)+\chi_2(s)\rk|^{2}_{H^{1}_2} \rk\}^\frac 1r
.
\EEQS
The rest follows by the definition of $\mU$ and estimating the terms on the RHS.

\end{proof}
Before embarking onto the proof of Proposition \ref{help_cont}, we prove the following important inequality.
This Lemma contains the main step in the proof of Proposition \ref{help_cont}.
\begin{lemma}\label{con-p1}
For all  $(\chi_1, \eta_1),(\chi_2, \eta_2)\in\mU $ and $(u_1, \v_1),(u_2, \v_2)\in\mU $ such that  $(u_1, \v_1)=\mathscr{T}[ (\chi_1, \eta_1)]$ and $(u_2, \v_2)=\mathscr{T}[(\chi_1, \eta_1)]$,
 there exists a constant $C=C(K_1,K_2)>0$, some numbers $\delta_1, \delta_2 >0$ and $\gamma\in(0,1)$ such that
\begin{equation}\label{L1-chi2xi}
\begin{split}
&\mathbb{E} \lk(\int_0^T \big| \chi_1^2(s) \xi_1(s) - \chi_2^2(s) \xi_2(s)\big|_{L^1}\,ds \rk) 
\\ &\leq C(K_1, K_2, T) \times \Big[ \big\{ \mathbb{E} \sup_{0\le s\le T} |\chi_1(s) - \chi_2(s)|^2_{{H^{1-\rho}_2}} \big\}^{\delta_1} + \big\{ \mathbb{E} \sup_{0\le s\le T} |\xi_1(s) - \xi_2(s)|^{\gamma }_{L^1} \big\}^{\delta_2} \Big],
\end{split}
\end{equation}
where $\xi_i := v_i^{-1}$ for $i=1, 2.$
\end{lemma}

\begin{proof}[Proof of Lemma \ref{con-p1}:]
Firstly, let us  note that
\begin{equation}\label{split-chi-xi}
\begin{split}
|\chi_1^2(s) \xi_1(s) - \chi_2^2(s) \xi_2(s)| \leq |\chi_1(s) - \chi_2(s)| |\chi_1(s) + \chi_2(s)| |\xi_1(s)| + |\chi_2^2(s)| |\xi_1(s) - \xi_2(s)|.
\end{split}
\end{equation}
Secondly, let us note that for any $n \in \mathbb{N}$ we have
\begin{align*}
|a-b| &\leq \bigg| \sum_{k=1}^n a^{\frac{k-1}{n}} (a^{\frac{1}{n}} - b^{\frac{1}{n}}) b^{\frac{n-k}{n}}\bigg| \leq C (n-1) (a-b)^{\frac{1}{n}} \big( a^{\frac{n-1}{n}} + b^{\frac{n-1}{n}} \big).
\end{align*}
Using the above inequality, we infer that 
\begin{align*}
\int_0^T \big| \chi_1^2(s) \xi_1(s) - \chi_2^2(s) \xi_2(s)\big|_{L^1} \,ds
\leq \int_0^T \big| |\chi_1^2(s) \xi_1(s) - \chi_2^2(s) \xi_2(s) |^{\frac 1n} \cdot \Psi(s) \big|_{L^1} \,ds,
\\
\leq    \int_0^T \big| \Big( (\chi_1^2(s)-\chi_2^2(s) ) \xi_1(s) +\chi_2^2(s) (\xi_1(s)-\xi_2(s) )\Big) ^{\frac 1n} \cdot \Psi(s) \big|_{L^1} \, ds,
\end{align*}
where
\begin{align}\label{Psi_n}
\Psi(s) := (n-1) \Big[ \big(\chi_1^2\, \xi_1 \big)^{\frac{n-1}{n}}(s) + \big(\chi_2^2\, \xi_2 \big)^{\frac{n-1}{n}}(s) \Big].
\end{align}
Now, let us fix some numbers $p_1,q_1,n\ge 1$ with $n\in\NN$, $n\ge 2$, $\frac 1{p_1}+\frac 1{q_1}\le \frac 12$,  
 and
\DEQSZ\label{cond_qq}
 {\frac {q_1(n-1)}{q_1n-q_1-1}}\le 2.
\EEQSZ
%
%
%
We split the right hand side of the above inequality into two summands, i.e., $S_1$ and $S_2$.
In particular, let
\DEQS
S_1&:= &  \int_0^T \big| \Big( (\chi_1^2(s)-\chi_2^2(s) ) \xi_1(s)\Big) ^{\frac 1n} \cdot \Psi(s) \big|_{L^1} \, ds,
\\
S_2&:= &  \int_0^T \big| \Big( \chi_2^2(s) (\xi_1(s)-\xi_2(s) )\Big) ^{\frac 1n} \cdot \Psi(s) \big|_{L^1} \, ds.
\EEQS
By \eqref{split-chi-xi}, using the H\"older inequality we obtain
\begin{equation}\label{S1+S2}
\begin{split}
S_1 &:=
\int_0^T   \big| |\chi_1(s) - \chi_2(s)|^{\frac{1}{n}} \big|_{L^{2n}}  \big| |\chi_1(s) + \chi_2(s)|^{\frac{1}{n}} \big|_{L^{{p_1}n}}  \big| |\xi_1(s)|^{\frac{1}{n}} \big|_{L^{{q_1} n}} |\Psi(s)|_{L^{\frac{n}{n-1}}}
%
%
\\ &\leq\int_0^T    |\chi_1(s) - \chi_2(s)|_{L^2}^{\frac{1}{n}} |\chi_1(s) + \chi_2(s)|_{L^{p_1}}^{\frac{1}{n}} |\xi_1(s)|_{L^{{q_1}}}^{\frac{1}{n}}
 |\Psi(s)|_{L^{\frac{n}{n-1}}}\, ds.
\end{split}
\end{equation}
%
We have
\DEQS
\mathbb{E}[S_1]&\leq& \mathbb{E} \Big(\int_0^T  |\chi_1(s) - \chi_2(s)|_{L^2}^{\frac{1}{n}}\, |\chi_1(s) + \chi_2(s)|_{L^{p_1}}^{\frac{1}{n}}\, |\xi_1(s)|_{L^{{q_1}}}^{\frac{1}{n}} |\Psi(s)|_{L^{\frac{n}{n-1}}} \,ds\Big).
\EEQS
Note that
we have by interpolation between $H^{1-\rho}_2(\CO)$ and $H^1_2(\CO)$
 for $\theta=1/\rho $
 $$
|\chi_1(s)-\chi_2(s)|_{L^2}\le |\chi_1(s) - \chi_2(s)|_{H^{1-\rho}_2}^\theta |\chi_1(s) -\chi_2(s)|^{1-\theta} _{H^{1}_2}.
$$
\newcommand{\m}{2}
In this way, we obtain by again applying the H\"older inequality
\DEQS
\mathbb{E}[S_1] &\leq &\mathbb{E} \bigg\{ \Big[ \sup_{s \in [0, T]} |\chi_1(s) - \chi_2(s)|_{H^{1-\rho}_2}^2\Big]^{\frac \theta{2n} }
 \,
 \Big[ \int_0^T |\chi_1(s) - \chi_2(s)|^2
 _{H^{1}_2} ds \Big]^\frac {1-\theta}{2n}
\\ &&\hspace{-1.5cm}{} \times   \Big[\int_0^T |\chi_1(s) + \chi_2(s)|^{2}_{L^{p_1}} ds\Big]^{\frac{1}{2n}}\Big[\sup_{s \in [0, T]} |\xi_1(s)|_{L^{q_1}}^{q_1} \Big]^{\frac{1 }{q_1n}} \,\Big[\int_0^T |\Psi(s)|^{\frac {2n}{2n-2+\theta}}
_{L^{\frac n{n-1}}} \,ds \Big]^{\frac {2n-2+\theta}{2n}} \bigg\}.
\EEQS
Then, applying again the H\"older inequality
\begin{equation}\label{eq-2.8-cont}
\begin{split}
\mathbb{E}S_1 &\leq \Big[\mathbb{E}  \sup_{s \in [0, T]} |\chi_1(s) - \chi_2(s)|^2_{H^{1-\rho}_2}\Big]^{\frac \theta{2n} } \,\Big[ \EE \Big( \int_0^T \big( |\chi_1(s) |_{H^{2-\rho}_2}+| \chi_2(s)|
 _{H^{2-\rho}_2}\big)^\m \, ds  \Big) \Big]^\frac {1-\theta}{2n} 
 \\ &\quad\times \Big[  \EE \Big(\int_0^T |\chi_1(s) + \chi_2(s)|^{2}_{L^{p_1}} \,ds\Big)\Big]^\frac 1{2n}
 \\ &\quad\times
\Big[ \EE \sup_{s \in [0, T]} |\xi_1(s)|_{L^{q_1}}^{{q_1 }} \Big]^{\frac{1}{q_1n}} \,
\bigg[\EE \Big(\int_0^T |\Psi(s)|^{\frac {2n}{2n-2+\theta}}_{L^{\frac{n}{n-1}}} \,ds\Big)^{\frac {2n-2+\theta}{2n}\frac {q_1n}{q_1n-q_1-1} } \bigg]^{\frac  {q_1n-q_1-1}{q_1n}} .
\end{split}
\end{equation}
From the definition of $\Psi$ in \eqref{Psi_n}, and, since
\begin{align*}
|\Psi(s)|_{L^{\frac{n}{n-1}}} &\leq C\, \Big( |(\chi_1^2\, \xi_1)^\frac {n-1}n (s)| _{L^{\frac{n}{n-1}}} + |(\chi_2^2\, \xi_2)^\frac {n-1}n (s)| _{L^{\frac{n}{n-1}}} \Big)
\\ &\leq C \Big[  \Big(\int_\CO (\chi_1^2\,\xi_1)(t, x)\, dx \Big)^\frac {n-1}n + \Big(\int_\CO (\chi_2^2\,\xi_2)(t, x)\, dx \Big)^\frac {n-1}n \Big],
\end{align*}
we get
\DEQS
\lqq{  \int_0^T |\Psi(s)|^{\frac {2n}{2n-2+\theta}} _{L^{\frac{n}{n-1}}} \,ds} &&
\\
 &\le& C \Big[\int_0^T \Big( \int_\CO (\chi_1^2\,\xi_1)(s, x)\, dx \Big)^{\frac {n-1}n \frac {2n}{2n-2+\theta}} \,ds + \int_0^T \Big( \int_\CO (\chi_2^2\,\xi_2)(s, x)\, dx \Big)^{\frac {n-1}n\frac {2n}{2n-2+\theta}} \,ds \Big]
\\
 &\le& C(\CO,T) \Big[ \Big(\int_0^T  \int_\CO (\chi_1^2\xi_1)(s, x)\, dx  \,ds\Big)^\frac {2(n-1)}{2n-2+\theta} + \Big(\int_0^T  \int_\CO (\chi_2^2\,\xi_2)(s, x)\, dx  \,ds\Big)^\frac {2(n-1)}{2n-2+\theta} \Big],
\EEQS
and
\begin{align*}
&\bigg[ \EE  \big| \Psi(s) \big|_{L^{\frac {2n}{2n-2+\theta}}(0,T;{L^{\frac{n}{n-1}}})}^{{\frac {2n-2+\theta}{2n}\frac {q_1n}{q_1n-q_1-1} } } \bigg]
\\
 &\le C \Big[ \EE  \Big(\int_0^T  \int_\CO (\chi_1^2\xi_1)(s, x)\, dx  \,ds\Big)^  {\frac {q_1(n-1)}{q_1n-q_1-1}} + \EE  \Big(\int_0^T  \int_\CO (\chi_2^2 \xi_2)(s, x)\, dx  \,ds\Big)^  {\frac {q_1(n-1)}{q_1n-q_1-1}} \Big].
\end{align*}
Due to \eqref{cond_qq}
we know
$$
 {\frac {q_1(n-1)}{q_1n-q_1-1}} \le 2.
$$
Using the previous energy inequalities results in \eqref{eq-2.8-cont} we infer that
\begin{equation}\label{S_1-est}
\mathbb{E}[S_1] \leq \Big[ \mathbb{E}  \sup_{s \in [0, T]} |\chi_1(s) - \chi_2(s)|^2_{H^{1-\rho }_2}\Big]^{\frac{1 }{4n} } \times K_1^{\frac{1- \theta}{2n}+\frac 1{2n}+\frac 1{nq_1}} \times K_2.
\end{equation}
In the next step, we deal with $S_2$, where $S_2$ is given by
\DEQS
\mathbb{E}[S_2]:= \EE \int_0^ T  \big| \chi^{\frac{2}{n}}_2(s) \,|\xi_1(s) - \xi_2(s)|^{\frac{1}{n}}\,\Psi(s) \big|_{L^1}\, ds .
\EEQS
Let us choose $n\in\NN$ with $\gamma\in(0,1)$, $q_2=1$, $p_2> 2$, $m_0\in(1,\infty)$, and $r<2$,  such that 
\DEQSZ\label{p1_cond}
\frac 1{n}+\frac \theta {\gamma n}+\frac {1-\theta}{m_0n}+\frac {n}{2(n-1)}\le 1, \quad
 \frac 1{p_2}=\frac \theta {q_2}+\frac {1-\theta}{m_0}, \quad \mbox{and} \quad \frac 1{r}+\frac 1{p_2}\le 1.
\EEQSZ
A short calculations gives that for $\gamma=\frac 9{10}$, $n=5$, $p_2=5/2$, and $r=5/3$, there exists a number $\theta\in(0,1)$ and a number $m_0\ge 1$ ($\theta=9/40, m_0=31/7$) satisfying  the conditions  above.
An application of the H\"older inequality gives 
\begin{align*}
\mathbb{E}[S_2]
&\le  \EE \int_0^ T  | \chi^{\frac{2}{n}}_2(s)|_{L^{rn}}  \big| |\xi_1(s) - \xi_2(s)|^{\frac{1}{n}} \big|_{L^{p_2n}}  |\Psi(s)|_{L^{\frac{n}{n-1}}}\, ds
\\
 &\le  \EE \int_0^ T  | \chi_2(s)|_{L^{2r}}^\frac 2n  \big|\xi_1(s) - \xi_2(s)\big|_{L^{p_2}}^\frac 1n  |\Psi(s)|_{L^{\frac{n}{n-1}}}\, ds.
\end{align*}
Interpolation gives for $\frac 1{p_2}=\frac \theta {q_2}+\frac {1-\theta}{m_0}$ 
$$
  \big|\xi_1(s) - \xi_2(s)\big|_{L^{p_2}}\le   \big|\xi_1(s) - \xi_2(s)\big|_{L^{q_2}}^\theta     \big|\xi_1(s) - \xi_2(s)\big|_{L^{m_0}} ^{1-\theta}.
$$
Hence, applying the H\"older inequality
\DEQS
\lqq{ \mathbb{E}[S_2]\le\EE \Bigg\{  \lk( \int_0^ T  | \chi_2(s)|_{L^{2r}}^2 \, ds \rk)^\frac 1{n} \, \Big[ \sup_{0\le s\le T}\big|\xi_1(s) - \xi_2(s)\big|_{L^{q_2}}^\gamma \Big] ^\frac \theta {\gamma n} }
\\
&&{}\quad\times  \Big(  \sup_{0\le s\le T} |\xi_1(s) - \xi_2(s)|_{L^{m_0}}^{m_0} \Big)^\frac {(1-\theta)}{m_0n}
\, \Big(  \int_0^ T |\Psi(s)|^{\frac n{n-1}}_{L^\frac n{n-1}}\, ds\Big)^\frac {n-1}n
\Bigg\}.
\EEQS
Note, that we have
$$
 \Big(  \int_0^ T |\Psi(s)|^{\frac n{n-1}}_{L^{\frac{n}{n-1}}}\, ds\Big)\le \Big( \int_0^ T \int_\CO\big( \chi^2_1(s,x)\xi_1(s,x)+\chi^2_2(s,x)\xi_2(s,x)\big) dx \, ds \Big)^{\frac n{n-1}}.
 $$
Again applying the H\"older inequality gives
\DEQS
\lqq{ \mathbb{E}[S_2] \le \Big( \EE   \int_0^ T  | \chi_2(s)|^2_{L^{2r}} \, ds \Big)^\frac 1{n} \, \Big( \EE \sup_{0\le s\le T}\big|\xi_1(s) - \xi_2(s)\big|_{L^{q_2}}^\gamma \Big) ^\frac \theta{\gamma  n} }
\\
&&{}\times  \Big( \EE\sup_{0\le s\le T} \big(  |\xi_1(s)|_{L^{m_0}}^{m_0}+| \xi_2(s)|_{L^{m_0}}^{m_0} \big) \Big)^\frac {(1-\theta)}{nm_0}
\\
&&\times  \bigg( \EE \Big(\int_0^ T \int_\CO\big( \chi^2_1(s,x)\xi_1(s,x)+\chi^2_2(s,x)\xi_2(s,x)\big)\, dx \, ds\Big)^2  \bigg) ^{\frac n{2(n-1)}}.
\EEQS
In particular, for our choice of $q$ we get 
\begin{equation}\label{S_2-est}
\mathbb{E}[S_2] \leq K_1^{\frac{1}{n}} \times K_2^{\frac{(1-\theta)}{m_0 n}} \times \Big( \EE \sup_{0\le s\le T}\big|\xi_1(s) - \xi_2(s)\big|_{L^{1}}^\gamma\Big) ^\frac \theta{2 n}.
\end{equation}
Using the estimates \eqref{S_1-est} and \eqref{S_2-est} in \eqref{S1+S2}, there exists $\delta_1, \delta_2>0$ such that the inequality 
 \eqref{L1-chi2xi} holds.

\end{proof}

\noindent
Now we are ready to prove Proposition \ref{help_cont}.

\begin{proof}[Proof of Proposition \ref{help_cont}]
By analyticity of the semigroup we infer that
\begin{equation}\label{an-se}
\begin{split}
\bE \sup_{0 \leq s \leq t} |u_1(s) - u_2(s)|^2_{H^{1-\rho}_2} \,dt &\leq \bE \sup_{0 \leq s \leq t} \int_0^t \big|\chi_1^2(s) \xi_1(s) - \chi_2^2(s) \xi_2(s) \big|_{H_2^{1-\rho}} \,ds
\\ &\hspace{-3cm}{}+ \bE \sup_{0 \leq s \leq t} \left| \sigma_u \int_0^t (u_1(s) - u_2(s)) \,dW_1(s)\right|_{H^{1-\rho}_2}\,,
\end{split}
\end{equation}
where we used the results from previous sections.
Using the Burkholder-Davis-Gundy inequality and the Young inequality we infer,
\begin{equation}\label{n-term}
\begin{split}
&\bE \sup_{0 \leq s \leq t} \left| \sigma_u \int_0^t (u_1(s) - u_2(s)) \,dW_1(s)\right|_{_{H^{1-\rho }_2}}
\\ &\leq \epsilon \,\bE \sup_{0 \leq s \leq t} |u_1(s) - u_2(s)|^2_{_{H^{1-\rho}_2}} + \frac{C(\sigma_u, S(\gamma_1))}{4 \epsilon} \int_0^t \bE \sup_{0 \leq r \leq s} |u_1(r) - u_2(r)|^2_{_{H^{1-\rho}_2}} \,dr.
\end{split}
\end{equation}
{We consider the non-linear term. By Sobolev embedding theorem, we have g $L^1(\CO) \hookrightarrow H_2^{- \rho}(\CO)$.}
Using the previous technical result, i.e., Lemma \ref{con-p1}, 
we infer that there exists a constant $C>0$ and $\delta_1, \delta_2 >0$ such that
{\small{
\begin{equation}\label{n-l_term}
\begin{split}
&\mathbb{E} \int_0^t \big| \chi_1^2(s) \xi_1(s) - \chi_2^2(s) \xi_2(s)\big|_{H^{1-\rho}_2} \,ds
\\ &\le \mathbb{E} \int_0^t \big| \chi_1^2(s) \xi_1(s) - \chi_2^2(s) \xi_2(s)\big|_{L^1} \,ds
\\ &\leq C(K_1, K_2,K_3, T) \, \Big[ \Big\{ \mathbb{E} \big[ \sup_{s \in [0, t]} |\chi_1(s) - \chi_2(s)|^2_{H^{-\rho }_2} \big] \Big\}^{\delta_1} + \Big\{ \mathbb{E} \big[ \sup_{s \in [0, t]} |\xi_1(s) - \xi_2(s)|^\gamma_{L^1} \big] \Big\}^{\delta_2}\Big].
\end{split}
\end{equation}
}}
Substituting the first and second term by the estimates given in  \eqref{n-term} and  \eqref{n-l_term} in the inequality \eqref{an-se}, we obtain
{\small{
\DEQS 
\begin{split}
&\bE \sup_{0 \leq s \leq t} |u_1(s) - u_2(s)|^2_{H^{1-\rho}_2} \,dt
\\ &\leq \epsilon \,\bE \sup_{0 \leq s \leq t} |u_1(s) - u_2(s)|^2_{H^{1-\rho}_2} + \frac{C \sigma_u}{4 \epsilon} \int_0^t \bE \sup_{0 \leq r \leq s} |u_1(r) - u_2(r)|^2_{H^{1-\rho}_2} \,dr
\\ &\quad+ C(K_1, K_2, T)\,\Big[ \Big\{ \mathbb{E} \big[ \sup_{s \in [0, t]} |\chi_1(s) - \chi_2(s)|^2_{H^{1-\rho}_2} \big] \Big\}^{\delta_1} + \Big\{ \mathbb{E} \big[ \sup_{s \in [0, t]} |\xi_1(s) - \xi_2(s)|_{L^1}^\gamma \big] \Big\}^{\delta_2}\Big].
\end{split}
\EEQS 
}}
Suitably choosing $\epsilon$ and using the Gronwall lemma, we get the assertion.
\end{proof}

\smallskip

\section{Pathwise uniqueness in one dimension}\label{sec-uniqueness}

Since in the proof of the existence of solutions compactness arguments are used, the underlying probability space {gets}
lost. In this section we will prove pathewise uniqueness for the Gierer--Meinhardt system, i.e.  system \eqref{equ1ss}--\eqref{eqv1ss},
in dimension one.
Usually, by the Yamada-Watanabe theory we know that pathwise uniqueness for the system \eqref{equ1ss}--\eqref{eqv1ss} implies joint uniqueness in law {and strong existence} for the equation the system \eqref{equ1ss}--\eqref{eqv1ss}, see \cite{cherny,engelbert,jacod,martin1,tappe,chinese,yamada}.
By standard arguments, the existence of a strong solution follows.

\begin{definition}
The equations \eqref{equ1ss}-\eqref{eqv1ss} are said to be {\sl pathwise unique} if, whenever
\\$\big(\Omega,\CF,(\CF_t)_{t\in [0, T]}, \PP, (u_i,v_i), W_1, W_2 \big)$, $i=1, 2$, are solutions to
\eqref{equ1ss}-\eqref{eqv1ss} such that
$$\PP\lk( u _1 (0)=u _2 (0),\, v_1(0)=v_2(0)\rk)=1,
$$
 then
 $$\PP\lk( u _1 (t)=u _2 (t),\, v_1(t)=v_2(t)\rk)=1,$$ for every $0< t\le T$.
\end{definition}

Under certain assumptions, the pathwise  uniqueness of the stochastic Gierer-Meinhardt system can be shown.

\begin{theorem}\label{main_ex_u}
Let $d=1$.
Let be given a filtered probability space $\mathfrak{A}=(\Omega,\CF, \mathbb{F},\PP)$ with the filtration $\mathbb{F} = \{{{\mathcal{F}}}_t:t\in [0,T]\}$ satisfying the usual conditions. Let $W_1$ and $W_2$ be two independent Wiener processes  in $\CH:=L^2(\CO)$, defined over the probability space $\mathfrak{A}$, with covariances $Q_1$ and $Q_2$ and
two couples of solutions $(u_1,v_1)$ and $(u_2,v_2)$ to equation \eqref{equ1ss}-\eqref{eqv1ss} over $\mathfrak{A}$,
on $[0,T]$ such that 
{$(u_1,v_1)$ and $(u_2,v_2)$ are continuous on $L^2(\CO)$.
If there exists some
and the solutions $(u_1,v_1)$ and $(u_2,v_2)$ are belonging $\PP$--a.s.\ to  $C_b([0,T];L^2(\CO))\times C_b([0,T];L^2(\CO))$,}
then $(u_1,v_1)$ and $(u_2,v_2)$ are indistinguishable in $L^2(\CO)$. 
\end{theorem}

\medskip

\noindent
Since $(u_1, v_1)$ and $(u_2, v_2)$ are solutions to the system \eqref{equ1ss}-\eqref{eqv1ss}, $\PP\lk( u _1 (0)=u _2 (0),\, v_1(0)=v_2(0)\rk)=1$, then we can write
\DEQSZ\label{eq:ui}
 d u_i(t) & =&  \Big[r_u   \Delta u_i(t) + \kappa_u\, \frac {u_i^2(t)}{\v_i(t)} - \mu_u u_i(t) \Big]\, dt +\sigma_u u_i(t)\circ  dW_1(t),\quad t>0, i=1, 2,
\EEQSZ
and
\DEQSZ\label{eq:vi}
d{\v_i}(t) &=& \Big[ r_\v \Delta \v_i(t) + \kappa_v\, u_i^2(t) - \mu_\v \v_i(t) \Big]\,dt +\sigma_{\v} \v_i(t)\circ d W_2(t),\quad t>0, i=1, 2.
\EEQSZ
\del{
Using the technical Proposition \ref{interp_rho} we have
\DEQS
\del{\EE\Big[ \int_0^T |u(s)|_{L^5}^\frac {10}3\, ds\Big]^\frac 3 5<\infty, \ \mbox{ and } \ }\EE \lk[\int_0^t\big(|u^2_1(s)|_{L^4}+|u^2_2(s)|_{L^4}\big)^4\, ds\rk]^\frac 12 < \infty.
\EEQS
}

\begin{proof}
In the first step we will introduce a family of stopping times $\{\tau_m:m\in\NN\}$
and show that on the time interval $[0,\tau_m]$ the solutions $u_1$ and $u_2$ are indistinguishable.
In the last step, we will show that $\PP\lk( \tau_m<T\rk)\to 0$  for $m\to\infty$. From this, it follows that $u_1$ and $u_2$ are indistinguishable on the time interval $[0,T]$.

\medskip

\paragraph{\bf Step I}
{\bf Introducing the stopping times:}
 Let us define the stopping times for $m\in \NN$
 \begin{align*}
 \begin{cases}
 \tau^1_m(\xi):=\inf_{t>0} \Big\{ {\sup_{0\le s\le t} |\xi(s)|_{L^8} \ge m} \Big\}, \\
 \tau^2_m(u):=\inf_{t>0} \Big\{ \big(\int_0^ t |u(s)|_{H^1_2}^2\, ds \big)+ \sup_{0\le s\le t} |u(s)|_{L^2}^2 \ge m \Big\}.
\end{cases}
\end{align*}
Let us introduce the stopping times  $\{ \bar \tau^1_m:m\in\NN\}$ and  $\{\bar  \tau^2_m:m\in\NN\}$  given by
 $$\bar \tau^j_m:=\min(\tau^1_m(\xi_j),\tau^2_m(u_j))\wedge T,\quad j=1,2,
 $$
where $\xi_j:=v_j^{-1}$.
The aim is to show
that the couple $(u_1,v_1)$ and $(u_2,v_2)$ are indistinguishable on the time interval $[0,\tau_m]$, {with $\tau_m=\inf(\bar \tau^1_m, \bar \tau^2_m)$.}

\medskip
%
Fix  $m\in\NN$.
To get uniqueness on $[0,\tau_m]$ we first stop the original solution processes at time $\tau_m$ and extend the processes $u_1$ and $u_2$
by other processes to the whole interval $[0,T]$.  For this purpose,
let $y_j$ be a solution to
\DEQSZ\label{eq1}
 y_j(t) &=& e^{-(t-\tau_m)(r_u\Delta-\mu_u)} u_j(\tau_m) 
+ \int_{\tau_m}^t  e^{-(t-s)(r_u\Delta-\mu_u)} y_j (s) \,\circ dW_1(s)
,\quad t\ge \tau_m ,\, j=1,2,
\nonumber \EEQSZ
and let $z_j$ be a solution to
\DEQSZ\label{eq2}
z_j (t) &= & e^{-(t-\tau_m)(r_v\Delta-\mu_v)} v_j (\tau _m)
+
\int_{\tau_m}^t  e^{-(t-s)(r_v\Delta-\mu_v)}z_j (s) \,\circ dW_2(s)
,\quad t\ge \tau_m,\, j=1,2
. \nonumber  \EEQSZ
Since $u_1$ and $u_2,$ and $v_1$ and $v_2$ are continuous in $L^ 2(\CO)$, $u_1(\tau_m),u_2(\tau_m),v_1(\tau_m)$ and $v_2(\tau_m)$ are well defined and
 belong $\PP$--a.s.\ to $L^2 (\CO)$.
Since, in addition, $( e^{-t(r_u\Delta-\mu_v)})_{t\in\RR}$ and  $( e^{-t(r_v\Delta-\mu_v)})_{t\in\RR}$ are analytic semigroups on $L^2 (\CO)$, 
the existence of unique solutions $y_j$ and $z_j $ in $L^2(\CO)$, $j=1,2$,
can be shown by standard methods.

Now, let us define two  processes $\bar u_{1, m} $ and $\bar u_{2, m}$ which are
equal to $u_1 $ and $u_2 $ on the time interval $[0,\tau_m)$ and
follow the linear equations  $y_1$ and $y_2$ respectively afterwards.
In the same way, let us define two  processes $\bar v_{1, m} $ and $\bar v_{2, m}$ which are
equal to $v_1 $ and $v_2 $ on the time interval $[0,\tau_m)$ and
follow the linear equations  $z_1$ and $z_2$ afterwards.

\noindent
In particular, let
$$
\bar u_{j, m}  (t) = \bcase u_j(t) & \mbox{ for } 0\le t< \tau_m,\\
y_j (t) & \mbox{ for } \tau_m\le  t \le T;\ecase
$$
and
$$
\bar v _{j, m}  (t) = \bcase v_j (t) & \mbox{ for } 0\le t< \tau_m,\\
z_j (t) & \mbox{ for } \tau_m\le  t \le T.\ecase
$$
%
Note, that the couples $(\bar u_{1, m},\bar v_{1, m})$ and $(\bar u_{2, m},\bar v_{2, m})$ solve
the following system corresponding  to  \eqref{equ1ss}-\eqref{eqv1ss} for the time interval $t  \in [0, \tau_m)$:
\DEQSZ\label{equ1ss2}
 d\bar {u}(t) & =&  \Big[r_u   \Delta \bar u(t) + \kappa_u\, \frac {\bar u^2(t)}{\bar{\v}(t)} - \mu_u \bar{u}(t) \Big]\, dt +\sigma_u \bar{u}(t)\circ  dW_1(t),
\EEQSZ
and
\DEQSZ\label{eqv1ss2}
d\bar{{\v}}(t) &=& \Big[ r_\v \Delta \bar{\v}(t) + \kappa_v\, \bar{u}^2(t) - \mu_\v \bar{\v}(t) \Big]\,dt +\sigma_{\v} \bar{\v}(t)\circ d W_2(t).
\EEQSZ
In addition, let us note that the process  $\xi_{j,m}:=v^{-1}_{j,m}$ solves \eqref{eq-xi}.

\smallskip
\noindent
\paragraph{\bf Step II}
The aim in this step is to show that the $\EE |\bar u_{1,m}(t)-\bar u_{2,m}(t)|_{L^2}^2 = 0$. The first step is to show that for all $m\in\NN,$ there exists a $C(m)>0$ such that
\DEQS
\EE |\bar u_{1,m}(t)-\bar u_{2,m}(t)|_{L^2}^2 &\le & C(m) \int_0^ t \EE |\bar u_{1,m}(s)-\bar u_{2,m}(s)|_{L^2}^2 \, ds  .
\EEQS
Let us define the subsets of $\Omega$
\begin{align}
\begin{cases}
\Omega_1^m:=\{ \omega\in\Omega:\tau^1_m(\bar \xi_{1,m})\le t,\,\tau^1_m(\bar \xi_{2,m})\le t\},
\\
\Omega_2^m:= \{ \omega\in\Omega:\tau^2_m(\bar u_{1,m})\le t,\,\tau^2_m(\bar u_{2,m})\le t\}.
\end{cases}
\end{align}
In addition, for a stochastic process $\zeta$, let us define the convolution process $\mathfrak{C}$ by
\DEQS
\mathfrak{C}(\zeta)(t):=\int_0^ t e^{-(t-s)\, r_u A } \,\zeta(s)\, ds,\quad t\in[0,T].
\EEQS
For simplicity we omit the index $m$ in the following paragraph and write $u_j$ and $v_j$ instead of $\bar u_{j,m}(t)$ and $\bar v_{j,m}(t)$ for $j=1,2$.
Since the nonlinear term is the most critical term, we start by analysing the nonlinear term. Let us fix $\gamma\in(\frac 12,\frac {14}{15})$.
Analysing
the convolution term with respect to the nonlinear term by splitting it into two parts. By applying the smoothing property of the semigroup and the embedding $L^1(\CO)\hookrightarrow H^{-\gamma}_2(\CO)$ we obtain for $p=2$
\DEQS
\lqq{ |\mathfrak{C}(u_1^2\xi_1\mathds{1}_{\Omega_1}\mathds{1}_{\Omega_2})(t)- \mathfrak{C}(u_2^2\xi_2\mathds{1}_{\Omega_1}\mathds{1}_{\Omega_2})(t)|^p_{L^2}}
&&
\\
&\le& \lk(  \int_0^t(t-s)^{-\frac\gamma 2} |u^2_1(s)\xi_1(s)-u^2_2(s)\xi_2(s)|_{H^{-\gamma}_2}\mathds{1}_{\Omega_1}\mathds{1}_{\Omega_2}  \,ds \rk)^p
\\
&\le &
 \lk( \int_0^t(t-s)^{-\gamma/2} |u_1(s)-u_2(s)|_{L^2} |(u_1(s)+ u_2(s)) \xi_1(s)|_{L^2}\,\mathds{1}_{\Omega_1}\mathds{1}_{\Omega_2}  ds\rk)^p
\\
&&{}+
\lk(  \int_0^t(t-s)^{-\gamma/2}  (|u^2_1(s)|_{L^4} + |u^2_2(s)|_{L^4})  |\xi_1(s)-\xi_2(s)|_{L^{\frac{4}{3}}}\mathds{1}_{\Omega_1}\mathds{1}_{\Omega_2} \, ds\rk)^p
  \\
  &=:& S^p_1(t)+S^p_2(t).
  \EEQS
\noindent
First we consider the term $S^p_1$. Due to the specific choice of $\gamma,$ we apply the H\"older inequality with $k' = \frac {15}{4}$ and $k=\frac {15}{11}$ to obtain
\begin{equation*}
S_1^p(t) \le C(t) \Big( \int_0^t |u_1(s)-u_2(s)|_{L^2}^{k} \big(|u_1(s)|_{L^5}+|u_2(s)|_{L^5}\big)^{k} \,|\xi_1(s)|^{k}_{L^\frac {10}3}\mathds{1}_{\Omega_1}\mathds{1}_{\Omega_2} \, ds \Big)^{p/k}.
\end{equation*}
Again applying the H\"older inequality we get for
 $q=\frac {22}{15}$
\begin{equation*}
\begin{split}
S_1^p(t) &\le C(t) \Big( \int_0^t |u_1(s)-u_2(s)|^{kq}_{L^2}\, ds \Big)^{\frac{p}{kq}}
\\ &\quad \times \Big( \int_0^ t \big(|u_1(s)|_{L^5}+|u_2(s)|_{L^5} \big)^{kq} \mathds{1}_{\Omega_2}\, ds \Big)^{\frac{p}{kq}} \,\sup_{0\le s\le t}\mathds{1}_{\Omega_1} |\xi_1(s)|^p_{L^\frac {10}3}.
\end{split}
\end{equation*}
By our choice $p = 2 = \frac{15}{11} \times \frac{22}{15} = k \times q$. Using the embedding $H^1_2(\CO) \hookrightarrow L^5(\CO)$, then using the definition of subsets $\Omega_1, \Omega_2$ and the definition of stopping times we infer that
\begin{equation}\label{eq-2.10}
\begin{split}
S_1^p(t) &\le C \Big( \int_0^t |u_1(s)-u_2(s)|^{p}_{L^2}\, ds \Big) \Big( \int_0^ t \big(|u_1(s)|^p_{L^5}+|u_2(s)|^p_{L^5} \big) \mathds{1}_{\Omega_2}\, ds \Big)\,\sup_{0\le s\le t}\mathds{1}_{\Omega_1} |\xi_1(s)|^p_{L^8}
\\ &\leq C(t, m) \int_0^t |u_1(s)-u_2(s)|^{p}_{L^2}\, ds.
\end{split}
\end{equation}
Next, we investigate the second summand. From Proposition \ref{interp_rho} we observe that
$$\|\cdot\|_{L^4(0, T; L^8)} \leq \| \cdot \|_{L^2(0,T; H^{1}_2)}+ \|\cdot\|_{L^\infty(0,T; {L^2})}.
$$ Using this result, then applying again the H\"older inequality and using the definition of subsets $\Omega_1, \Omega_2$ and the definition of stopping times we infer that
\begin{equation}\label{eq-2.11}
\begin{split}
S_2^p(t)&\le \lk(  \int_0^t(t-s)^{-\gamma/2} \big(|u^2_1(s)|_{L^4}+|u^2_2(s)|_{L^4} \big) |\xi_1(s)-\xi_2(s)|_{L^\frac 43}\mathds{1}_{\Omega_1}\mathds{1}_{\Omega_2} \, ds\rk)^p
\\ &\le C(t)\lk(  \int_0^t\mathds{1}_{\Omega_2}\big(|u^2_1(s)|_{L^4}+|u^2_2(s)|_{L^4}\big)^4\, ds \rk)^p \sup_{0\le s\le t} \mathds{1}_{\Omega_1}\, |\xi_1(s)-\xi_2(s)|_{L^\frac 43}^p
\\ &\le C(t) \lk(  \int_0^t\mathds{1}_{\Omega_2}\big(|u_1(s)|^4_{L^8}+|u_2(s)|^4_{L^8}\big)\, ds \rk)^p \sup_{0\le s\le t} \mathds{1}_{\Omega_1}\, |\xi_1(s)-\xi_2(s)|_{L^\frac 43}^p
\\ &\le C(t, m) \sup_{0\le s\le t} \mathds{1}_{\Omega_1}\,|\xi_1(s)-\xi_2(s)|_{L^\frac 43}^p.
\end{split}
\end{equation}
Using the above two estimates \eqref{eq-2.10} and \eqref{eq-2.11} we get
$$
S_1^p(t) + S_2^p(t) \leq C(t, m) \sup_{0\le s\le t} \mathds{1}_{\Omega_1}\,|\xi_1(s)-\xi_2(s)|_{L^\frac 43}^p + C(t, m) \int_0^t |u_1(s)-u_2(s)|^p_{L^2}\, ds .
$$
Now, to tackle the stochastic convolution part, for a stochastic process $\bar{\zeta}$, we define the operator
\begin{align*}
\mathfrak{S} (\bar{\zeta})(t)=\int_0^ t e^{-(t-s)\, r_u A}\, \bar{\zeta}(s)\,d W_1(s).
\end{align*}
Here it is important that the noise coefficient depends linearly on $u,$ and we make use of the Corollary \ref{brz:con} in the Appendix to get the desired estimate.

\noindent
Recollecting all the terms we arrive at
\begin{align}\label{xi+u1}
&|\mathfrak{C}(u_1^2\xi_1\mathds{1}_{\Omega_1}\mathds{1}_{\Omega_2})(t)- \mathfrak{C}(u_2^2\xi_2\mathds{1}_{\Omega_1}\mathds{1}_{\Omega_2})(t)|^p_{L^2} \nonumber
\\ &\le C(t, m) \sup_{0\le s\le t} \mathds{1}_{\Omega_1}\,|\xi_1(s)-\xi_2(s)|_{L^\frac 43}^p + C(t, m) \int_0^t |u_1(s)-u_2(s)|^p_{L^2}\, ds .
\end{align}
It remains to estimate the entity $ \EE\sup_{0\le s\le t}|\xi_1(s)-\xi_2(s)|_{L^\frac 43}^p$.
Here, we obtain by the H\"older inequality
\DEQS
|\xi_1(s)-\xi_2(s)|^p_{L^\frac 43}\le |\xi_1(s)\xi_2(s)|^p_{L^4} \,|v_1(s)-v_2(s)|^p_{L^2}.
\EEQS
{Since $|\xi_1(s)\xi_2(s)|^p_{L^4} \leq |\xi_1(s)|^p_{L^8}\,|\xi_2(s)|^p_{L^8}$, we use the definition of the stopping times and the definition of $\Omega_1$
and $\Omega_2$ to handle the term $|\xi_1(s)\xi_2(s)|^p_{L^4}$. In this way we estimate it by a constant depending on $m$. }

We only have to estimate $\sup_{0\le s\le t} |v_1(s)-v_2(s)|^p_{L^2}$. For this, we first calculate the nonlinear term which is the difficult term to handle. Here, using the smoothing property of the semigroup, the embedding $L^1(\CO)\hookrightarrow H^{-\gamma}_2(\CO)$, the H\"older inequality and the definition of stopping times we infer that
\begin{equation}\label{vvv}
\begin{split}
&|\mathfrak{C}(u_1^2 \mathds{1}_{\Omega_1}\mathds{1}_{\Omega_2})(t)- \mathfrak{C}(u_2^2 \mathds{1}_{\Omega_1}\mathds{1}_{\Omega_2})(t)|^2_{L^2}
\\ &\leq \Big( \int_0^t (t-s)^{-\frac{\gamma}{2}} \,|u_1^2(s) - u_2^2(s)|_{H^{-\gamma}_{2}} \,ds \Big)^2
\\ &\le  \lk(\int_0^ t \mathds{1}_{\Omega_1}\mathds{1}_{\Omega_2} (t-s)^{-\frac \gamma 2 }|u_1(s)-u_2(s)|_{L^2}\big( |u_1(s)|_{L^2}+ |u_2(s)|_{L^2}\big)\, ds\rk)^2
\\ &\le  C(t)\, \sup_{0\le s\le t} \mathds{1}_{\Omega_2} \big( |u_1(s)|_{L^2}+ |u_2(s)|_{L^2}\big)^2 \, \int_0^ t |u_1(s)-u_2(s)|^2_{L^2}\, ds
\\ &\le C(t, m)\, \int_0^ t |u_1(s)-u_2(s)|^2_{L^2}\, ds.
\end{split}
\end{equation}
Using the above estimates in \eqref{xi+u1} and taking expectation we get
\DEQS
\lqq{ \EE |v_1(t)-v_2(t)| ^2_{L^2}
} &&
\\
&\le&
 \lk|\mathfrak{C}(u_1^2\mathds{1}_{\Omega_1}\mathds{1}_{\Omega_2})(t)- \mathfrak{C}(u_2^2\mathds{1}_{\Omega_1}\mathds{1}_{\Omega_2})(t)\rk|^2_{L^2}
+\lk|\EE \int_0^t\la  (v_1(s)-v_2(s)),(v_1(s)-v_2(s))dW_2(s)\ra \rk|
\\
&\le &
  C(t,m) \int_0^ t \EE |u_1(s)-u_2(s)|^2_{L^2}\, ds + C \int_0^ t \EE |v_1(s)-v_2(s)|^2_{L^2}\, ds
  .
\EEQS
Finally taking into account the linear part, we can show that for any $\ep>0$ and $m\in\NN$
there exists a constant $C(m,T,\ep)>0$ such that
\DEQS
 \EE |u_1(t)-u_2(t)|^p_{L^2}\le C(T, m,\ep)  \int_0^t\EE|u_1(s)-u_2(s)|^2_{L^2}\, ds.
\EEQS
For $p=2$, applying the Grownwall Lemma, we get
\DEQS
\EE |u_1(t)-u_2(t)|^p_{L^2}\le0.
\EEQS
Similarly, by estimate \eqref{vvv} and taking into account the linear part and stochastic convolution, we obtain
\DEQS
\EE |v_1(t)-v_2(t)|^p_{L^2}\le0.
\EEQS

\medskip

\paragraph{\bf Step III}
We show that $\PP\lk( \tau_m<T\rk) \to 0$ as $m\to\infty$.
Observe, that
we have 
\DEQS
\{ \tau_m\le T \} \subset  \{ |\xi_1|_{L^{\infty}([0,T];L^{8})}\; \mbox{or} \;
 |\xi_2|_{ L^{\infty}([0,T];L^{8})} \ge m \}
 \\
 \cup \{ |u_1|_{L^{\infty}([0,T];L^2)\cap {L^{2}([0,T];H^1_2)}}\ge m \; \mbox{or} \;
 |u_2|_{ L^{\infty}([0,T];L^2)\cap {L^{2}([0,T];H^1_2)}}\ge m \}.
\EEQS
Therefore,
\DEQS
\PP\lk( \tau_m <  T \rk) \le \PP\lk( |u_1|_{L^{\infty}([0,T];L^2)\cap {L^{2}([0,T];H^1_2)}} \ge m \rk) + \PP\lk(  |u_1|_{L^{\infty}([0,T];L^2)\cap {L^{2}([0,T];H^1_2)}}\ge m \rk)
\\
+\PP\lk( |\xi_1|_{L^{\infty}([0,T];L^{8})}\ \ge m \rk) + \PP\lk(  |\xi_2|_{L^{\infty}([0,T];L^{8})}\ge m \rk).
\EEQS
Since $u_1$ and $u_2$ are continuous  in $L^2(\CO)$ and
$\xi_1$ and $\xi_2$ are continuous  in $ L^{8}(\CO)$, and $\PP$-a.s.\
$$ \sup_{0\le s\le T}|u_1(s)|_{L^{\infty}([0,T];L^2)\cap {L^{2}([0,T];H^1_2)}}<\infty,
$$
$$ \sup_{0\le s\le T}|u_2(s)|_{L^{\infty}([0,T];L^2)\cap {L^{2}([0,T];H^1_2)}}<\infty,
$$
$ \sup_{0\le s\le T}|\xi_1(s)|_{L^{\infty}([0,T];L^{8})}<\infty$,
and $ \sup_{0\le s\le T}|\xi_2(s)|_{L^{\infty}([0,T];L^{8})}<\infty$,
 it follows that
as $m\to\infty$, $\PP\lk( \tau_m \le  T \rk)\to 0$.
Hence, both processes $(u_1, v_1)$ and $(u_2, v_2)$ are undistinguishable on $[0,T]$.
\end{proof}

\medskip

\appendix

\section{Some compactness results}\label{comp}

As mentioned in the introduction, we will shortly introduce the results of \cite{brzezniak}, and \cite{brzezniaGatarek}, respectively, which we use within the proof of the main result. For simplicity, we do not state the results in all generality.
First, let us state the setting within our framework.

\begin{assumption}\label{assum-1}
 Let
$E$ be a Hilbert  space and let $A$ be a linear map on $E$ satisfying the following conditions.
\begin{enumerate}
 \item \label{h2-a}
$-A$ is a positive operator\footnote{See Section I.14.1 in Triebel's monograph \cite{Triebel_1995}.} on $E$ with compact resolvent. In particular,
there
exists $M>0$ such that
\[ \Lve
(A+\lambda)^{-1}\Rve \le \frac{M}{1+\lambda}, \mbox{ for any }\lambda\ge 0;\]
 \item \label{h2-b} $-A$ is an
infinitesimal generator of an analytic semigroup of contraction
type in $E$.
\item \label{h3} $A$ has the BIP (bounded imaginary power) property, i.e.\ there exist some
 constants $K>0$ and $\vartheta\in [0,\frac\pi2)$
such that
\begin{equation}
\Vert A^{is} \Vert \le K e^{\vartheta |s|}, \; s \in \mathbb{R}.
\label{2.1}
\end{equation}

\end{enumerate}
\end{assumption}

Let us fix two real numbers $q\in (1,\infty)$ and $\mu\geq 0$.
Similarly to Brze{\'z}niak and G{\c{a}}tarek \cite{brzezniaGatarek}, for a Banach space $\BB$, we
define  a linear operator ${{\mathcal{A}}}$ by the formula
\begin{eqnarray}\label{operator_def}
D({{\mathcal{A}}})&=&\left\{ u\in L ^q(0,T ;\BB): \; A u \in L ^q(0,T
;E)\right\},
\nonumber\\
 {{\mathcal{A}}}u&:=&\left\{(0,T) \ni t \mapsto A (u(t)) \in \BB \right\},\quad u\in D({\mathcal{A}}).
\label{2.8}
\end{eqnarray}
It  is  known, see  \cite{venni}, that  if $A+\mu I$ satisfies
Assumption \ref{assum-1}, then so does
${{\mathcal{A}}} + \mu I$.
With $q$ and $\BB$ as above we define two operators ${\mathcal{B}}$ and $\Lambda$, see \cite{brzezniaGatarek}, by
 \begin{eqnarray} {\mathcal{B}} u&=&u^\prime,\quad u\in
D({\mathcal{B}}):= H_{0} ^ {1,q}(0,T;\BB).
\\
\Lambda&:=& {\mathcal{B}} +{{\mathcal{A}}}, \quad\quad D(\Lambda):=D({\mathcal{B}})\cap
D({{\mathcal{A}}}). \label{2.10}
\end{eqnarray}

Next, by
\cite{venni} and \cite{Giga+Sohr_1991}, since $\Lambda={\mathcal{B}}-\mu I +{{\mathcal{A}}} +\mu I$,  $\Lambda $ is a positive
operator. In particular, $\Lambda$ has a bounded inverse. The
domain $D(\Lambda)$ of $\Lambda$ endowed with the `graph' norm
\begin{equation}
\Vert u\Vert = \left\{ \int_0^T |u^\prime(s)|^q \, ds + \int_0^T
|Au(s)|^q \, ds \right\}^{\frac1q} \label{2.11a}
\end{equation}
is a Banach space.
Let us  present two results on the fractional powers
of the  operator $\Lambda$, see  \cite{brzezniak} for the proof.

\begin{proposition}\label{Prop:2.0} Assume that Assumption \ref{assum-1} is  satisfied.
Then,  for any $ \alpha \in(0, 1]$,  the operator
$\Lambda^{-\alpha}$ is a bounded linear operator in $L^q
(0,T;E)$, and for $\eta \in L^q(0,T;E)$,
\begin{eqnarray}
\left(\Lambda^{-\alpha}\eta\right) (t)&=&\frac1{\Gamma (\alpha)}  \int_0^t (t-s)^{\alpha -1}
e^{-(t-s)A}\eta (s) \, ds, \;\; t \in [0,T].
\label{2.13a}
\end{eqnarray}
\end{proposition}

\begin{lemma}\label{L:reg}
(Compare with Lemma  2.4 of \cite{brzezniaGatarek}) Let Assumption \ref{assum-1} be satisfied.
Suppose that the positive numbers $\alpha, \beta,\delta$ and $q>1$
satisfy
\begin{equation}
\alpha  -\frac1q   +  \gamma>\beta +\delta.
\label{cond:1}
\end{equation}
    If $T  \in (0,\infty)$, then  the operator
\begin{equation}
\Lambda ^ {-\alpha} : L^q (0,T;D(A^\gamma)) \to
{{C}}^\beta_b ([0,T];D(A^\delta)) \; \label{2.13b}
\end{equation}
is bounded.
\del{
If $T=\infty$  and   the semigroup $\{e^{-tA}\}_{t\ge 0} $ is
exponentially bounded on $E$, i.e., for some  $ a>0$, $C>0$
\begin{eqnarray}
 |e^{-tA}|_{L(E,E)} &\le& Ce^{-at}, \; t \ge 0,
\label{exp-bound}
\end{eqnarray}
then  for any $f \in L^q_0(0,T;D(A^\gamma))$ the function
$u =\Lambda ^ {-\alpha} f$ belongs to
${{\mathcal{C}}}^\beta([0,T];D(A^\delta))$.
 Moreover, the operator $\Lambda ^ {-\alpha}$ is a bounded map in the
above spaces.
\begin{equation}
\Lambda ^ {-\alpha} : L^q(0,T;D(A^\gamma)) \to
{{\mathcal{C}}}^\beta([0,T];D(A^\delta)), \label{2.13c}
\end{equation}
is bounded.}
\end{lemma}

Finally, we present slight modifications of Proposition 2.2
and Theorem 2.6 from  \cite{brzezniak}.
\begin{theorem}
\label{Th:compact}
Let Assumption \ref{assum-1} be satisfied.
 Assume that
$ \alpha \in (0,1]$ and $\delta,\gamma\geq  0$ are such that
\begin{equation}
\alpha-\frac1q +\gamma-\delta>-\frac1r.
\label{cond:2}
\end{equation}
 Then the   operator
 \begin{equation}
 A ^{\delta} \Lambda^{-\alpha}A ^ {-\gamma} : L^q (0,T;\BB)
\to  L^r(0,T;\BB)
\label{cond:3}
\end{equation}
is bounded. Moreover, if the operator $A^{-1}:\BB\to \BB$ is   compact, then the
operator in \eqref{cond:3}   is compact.
\end{theorem}
\begin{remark}\label{Prop:2.1-remark-F}
In view of Theorem \ref{Th:compact}
$$\Lambda ^ {-1}:L ^ p(0,T;E)\to L ^ p(0,T;E)
$$
is a well defined bounded linear operator for $p\in [1,\infty)$.
Observe, that $\Lambda^{-1} \eta$ is the solution to the Cauchy problem
$$
\lk\{ \begin{array}{rcl} \dot u(t)&=& A u(t)+\eta(t),\\u(0)&=&0.\end{array}\rk.
$$
\end{remark}
\begin{corollary}\label{C:comp} Assume the first set of assumptions of Theorem
\ref{Th:compact} are satisfied. Assume  that  three  non-negative
numbers $\alpha,  \beta, \delta$ satisfy  the  following condition
\begin{equation}
\alpha - \frac1q>\beta + \delta.
\label{cond:1a}
\end{equation}
 Then the operator
$\Lambda^{-\alpha}: L^q(0,T;\BB ) \to
{C}^\beta_b([0,T];D(A^\delta))$ is bounded. Moreover, if the
operator $A^{-1}:E\to E$ is  compact, then  the
operator $\Lambda^{-\alpha}: L^q(0,T;\BB ) \to {C}^\beta_b([0,T];D(A^\delta))$  is also compact.
\par
{ In particular, if $\alpha > \frac1 q$ and the operator $A^{-1}:\BB\to \BB$ is compact, the map
$\Lambda^{-\alpha}:L^q(0,T;\BB ) \to C([0,T];\BB )$ is compact.
\del{- Compare with the beginning of Section 8.}}
\end{corollary}
Let $\mathfrak{A}=(\Omega, \CF,\BF,\PP)$ be a complete probability space
and filtration $\BF=(\CF_t)_{t\ge 0}$  satisfying the usual conditions.
Let $H$ and $H_1$ be two Hilbert spaces. Let
 $W$ be an $H_1-$valued cylindrical Wiener process, $L_{HS}(H_1,H)$ be the space of linear operators belonging to $L(H_1,H)$ with finite Hilbert-Schmidt norm, %
and
\DEQSZ\label{def.MA}
\lqq{ \CN_{\MA}^p(H):=\Big\{ \xi:[0,T]\times \Omega\to H,\quad \mbox{such that} \quad }
\notag
 &&\\
&& \mbox{$\xi$ is progressively measurable over $\mathfrak{A}$ and}\quad \EE\int_0^T |\xi(s)|^p_{L_{HS}(H_1,H)}\, ds<\infty\Big\}.
 \EEQSZ
To handle the stochastic convolution, let us define for a process $\xi\in \CN_{\MA}^2(H)$
the operator $\mathfrak{S}_A$ by
$$\mathfrak{S}_A\xi(t)=\int_0^ t e^{-(t-s)A} \xi(s)d W(s)
.
$$
\begin {corollary}\label{brz:con}$($Compare with Corollary 3.5 of \cite[p.\ 266]{brzezniak}$)$
If non--negative numbers $\beta,\delta$ and $\nu$ satisfy the following
$$
\beta+\delta+\frac 1p<\frac 12 -\nu,
$$
then a process $x(t)$, $t\in[0,T]$ belonging to  $\CN_{\MA}^q(A^{-\nu}H)$  possesses a modification $\tilde x(t)$, $t\in[0,T]$ such that
$$
\tilde x \in C^{\beta}_b(0,T;D(A^\delta)), \quad a.s.,
$$
and there exists some constant $C_T>0$ such that
$$
\EE \|\tilde x\|^p_{\CC^{\beta}_b(0,T;D(A^\delta))}\le C_T \EE \int_0^T \|A^{-\nu}x(s)\|^p_{L_{HS}(H_1,H)}\, ds.
$$
\end{corollary}

\del{For a function $f\in L^q(0,T;E)$ we define the operator $\Lambda_{T,\alpha}$ by
$$
(\Lambda 
_{T,\alpha} f)(t)=\frac 1 {\Gamma(\alpha)} \int_0^ t(t-s)^{\alpha-1} e^{-(t-s)A}f(s)\, ds ,
t\in (0,T).
$$
Observe, that $\Lambda_{T,a} f$ is the solution to the Cauchy problem
$$
\lk\{ \begin{array}{rcl} \dot u(t)&=& A u(t)+f(t),\\u(0)&=&0.\end{array}\rk.
$$

\begin{lemma}
Assume a Banach space $E$ and a linear operator satisfy the condition \ref{assum-1}. Suppose that the positive numbers $\alpha,\beta,\delta$ satisfy
$$
0<\beta<\alpha-\frac 1q +\gamma-\delta.
$$
Then, if $T\in(0,\infty)$ and $f\in L^q(0,T;D(A^\gamma))$, the function $u=\Lambda_{T,\alpha}f$ satisfies
$$u\in\CC^{(\beta)}_b(0,T;D(A^\delta)).
$$
\end{lemma}}

\medskip

\section{An embedding result}
Let us fix
the   Banach space $\mathbb{H}_\alpha:= L^2(0,T; H^{\alpha+1}_2(\CO))\cap L^\infty(0,T; {H^{\alpha}_2}(\CO))$
equipped with the norm
\begin{equation}\label{eq:BZeq}
\|\xi\|_{\mathbb{H}_\alpha}:= \|\xi\|_{L^2(0,T; H^{\alpha+1}_2)}+ \|\xi\|_{L^\infty(0,T; {H^{\alpha}_2})},\quad\xi\in \mathbb{H}_\alpha.
\end{equation}
Then, we have the following embedding.
\begin{proposition}\label{interp_rho}
Let $l_1,l_2\in(2,\infty)$ with $\frac{d}{2}-\alpha\le \frac 2{l_1}+\frac d{l_2}$.
Then there exists $C>0$ such that
$$
\|\xi\|_{L^{l_1}(0,T;L^{l_2})}\le C\|\xi\|_{\mathbb{H}_\alpha},\quad \xi\in\mathbb{H}_\alpha.
$$
\end{proposition}
\begin{proof}
First, let us note that for $\delta>0$ such that
\begin{equation}\label{eq:deltaembedding}\frac {1}{l_2}\ge \frac{1}{2}-\frac{\delta}{d},
\end{equation}
we have the embedding $H^{\delta}_2(\CO)\hookrightarrow L^{l_2}(\CO)$,
and therefore
$$
\|\xi\|_{L^{l_1}(0,T;L^{l_2})}\le C\|\xi\|_{L^{l_1}(0,T;H^{\delta}_2)},\quad \xi\in\X.
$$
Due to interpolation, compare with \cite[Theorem 5.1.2, p.\ 107 and Theorem 6.4.5, p.\ 152]{bergh}, we have
$$
 \|\xi\|_{L^{l_1}(0,T;H^{\delta}_2)}\le  C\|\xi\|_{L^\infty(0,T;H^{\alpha}_2)}^\theta \|\xi\|_{L^2(0,T;H^{\alpha+1}_2)}^{1-\theta},\quad \xi\in\BH,
$$
for $\theta\in (0,1)$ with,
\begin{equation}\label{eq:delta_theta}\frac  1{l_1}\le \frac{1}{2}(1-\theta),\quad \delta\le\theta\alpha+(1-\theta)(\alpha+1).
\end{equation}
Taking into account that we have
\[\frac{d}{2}-\alpha\le \frac 2{l_1}+\frac d{l_2},\]
gives that $\delta$ and $\theta$ satisfying \eqref{eq:deltaembedding} and \eqref{eq:delta_theta} exist and thus the Young inequality for products yields the assertion.
\end{proof}

\medskip

\section{Positivity of the Solution}\label{positivity}

Let us introduce the unbounded operator $A$ to be  the Laplace operator $-\Delta$ in $L^2(\CO)$ with Neumann boundary conditions. In particular, let
\DEQS
\lk\{ \begin{array}{rcl}
D(A) &=&\{ w\in H^2_2(\CO): \frac {\partial}{\partial n}w(x)=0,\, x\in\partial \CO\},
\\
Aw&=&-\Delta-\mu , \quad w\in D(A).
\end{array}\rk.
\EEQS
Here, $n$ denotes the inward normal vector of the boundary $\partial \CO$.Let us denote $\mathfrak{A}=(\Omega,\CF, \mathbb{F},\PP)$ be a complete probability space with the filtration $\mathbb{F} = \{{{\mathcal{F}}}_t:t\in [0,T]\}$ satisfying the usual conditions i.e., $\PP$ is complete on $(\Omega, \CF)$, for each $t \geq 0, \CF_t$ contains all $(\CF, \PP)$-null sets, and the filtration $\BF$ is right-continuous. Let in case $d=1$ the domain $\CO$ be an interval and in case $d=2$, let $\CO\subset \RR^2$ be a bounded domain with $C^\infty$  boundary.
Let $W$ be a Wiener processes  in $\CH:=L^2(\CO)$, with $\int w(x)\, dx=0$,  defined over the probability space $\mathfrak{A}$.
Let $w$ be a solution to
\DEQSZ\label{eqpos}
\lk\{ \barray
dw^{w_0}(t)&=& Aw^{w_0}(t)+ g(t,w^{w_0}(t))+w^{w_0}(t)\, dW(t),
\\
w^{w_0}(0) &=& w_0.\earray\rk.
\EEQSZ
Given $w_0>0$, we are interested in the strict positivity of the solution $w(w_0)$.
In \cite{TessitoreZabczyk1998}  the positivity of the Ornstein-Uhlenbeck is shown, where the
operator $A$ is the Laplacian on $\RR^d$. Here, we want to show that the result can be transferred where the
operator $A$ denotes the Laplacian on a bounded domain with Neuman boundary conditions.
To show the positivity of the solution $w$ the Gaussian tail estimate given in Theorem 3.1 of \cite{TessitoreZabczyk1998}
is essential. In \cite{veraar,seidler1,seidler2} Gaussian tail estimates are given in case one has a $C_0$ semigroup of contractions.

In particular, for any $C_0$--semigroup of contractions on a Banach space of $M$-type 2,
and any 
$1\le q\le \infty$, there exists a constant $K>0$ such
we have the maximal inequality (see Theorem 5.3 \cite{veraar})
\DEQSZ\label{ete}
\PP\lk( \sup_{0\le t\le T} \lk|\int_0^ t e^{A(t-s)} \psi(s) dW(s)\rk|_{L^q(\CO)}\ge \ep\rk)\le K \exp\lk( -\frac { \ep^2}{\kappa^2}\rk),
\EEQSZ
for all $\ep>0$ and every process $\psi\in\gamma(\mathcal{H},L^q(\CO))$ satisfying
$$
\mbox{ess\,sup}_{\Omega} \int_0^T\|\psi(s)\|^2_{\gamma(U,L^q(\CO) ) }\, ds\le \kappa.
$$
Here $\gamma(\mathcal{H},L^q(\CO))$ denotes the space of all $\gamma$--radonifying operators from $\CH$ to $L^q(\CO)$.

The second ingredient are the following estimates of the heat semigroup with Neumann boundary.
Let $(e^{tA})_{t\ge 0}$ be the the heat semigroup with Neumann boundary  in $\CO$, and let $\lambda_1>0$ be
the first nonzero eigenvalue of $-\Delta$ in $\CO$ under Neumann boundary conditions. Then there exists $C_1,\ldots,C_4$ depending only on $\Omega$ such that we have the following estimates
\begin{enumerate}[label=(\roman*)]
\item
If $1\le q\le p\le\infty $, then
 we have for all $t>0$ and $w\in L^q(\CO)$ with $\int w(x)\, dx =0$
 \DEQS
| e^{tA  }w|_{L^p}\le C_1\lk( 1+t^{-\frac d2(\frac 1q-\frac 1p)}\rk)e^{-\lambda_1 t}|w|_{L^q}.
\EEQS
\item
If $1\le q\le p\le \infty$, then
 we have for all $t>0$ and $w\in L^q(\CO)$
 \DEQS
|\nabla  e^{tA  }w|_{L^p}\le C_2\lk( 1+t^{-\frac 12 -\frac d2(\frac 1q-\frac 1p)}\rk)|w|_{L^q}.
\EEQS
\item
If $2\le p< \infty$ then
 we have for all $t>0$ and $w\in W^1_q(\CO)$
 \DEQS
|\nabla  e^{tA  }w|_{L^p}\le C_3e^{-\lambda_1 t} |\nabla w|_{L^p}.
\EEQS
\item
Let  $1<q\le p< \infty$. Then
 we have for all $t>0$ and $w\in C_b^\infty(\CO)$
 \DEQSZ\label{exten}
|\nabla  e^{tA  }w|_{L^p}\le C_4e^{-\lambda_1 t} \lk( 1+t^{-\frac 12 -\frac d2(\frac 1q-\frac 1p)}\rk)|w|_{L^p}.
\EEQSZ
It follows that the operator $\nabla  e^{tA  }$ has a uniquely defined extension to an operator from $L^q(\CO)$ into $L^p(\CO)$, with norm according from \eqref{exten}.
\end{enumerate}
Next, since $C_b(\CO)\hookrightarrow W_q^1(\CO)$ continuously for $q>2$ and $d\in\{1,2\}$, we know
that for all $2< p<\infty$ there exists a constant $C>0$ such that we have for all $t>0$ and $w\in L^q(\CO)$
 \DEQSZ\label{one-eq}
|  e^{tA  }w|_{C_b}\le C_2\lk( 1+t^{-\gamma}\rk)|w|_{L^q}, \quad \gamma>\frac d{2q}.
\EEQSZ
Similarly, it follows
 \DEQSZ\label{two}
\lqq{ |  e^{tA  }w|_{\gamma(L^2,L^q)}^2\le \sum _{k=1}^\infty |e^{tA } w e_k|_{L^2}^2 }
&&
\\\nonumber
&\le& \sum _{k=1}^\infty |e^{tA } w |_{L^\infty }^2|e_k|_{L^2 }^2
\le  C\lk( 1+t^{-2 \frac dq}\rk)e^{-2\lambda_1 t}|w|^2_{L^q}.
\EEQSZ

After this preliminaries we can present the following  Theorem.

\begin{theorem}\label{positivitytheorem}
Let $q\ge2$ be such large that ${\frac 12>\frac dq+\frac {d}q}$ and $w_0\in L^q(\CO)$ strict positive.  Then, for all $t>0$
$\PP\otimes \mbox{Leb}$-- a.s. the solution  $w^{w_0}(t,x)>0$.
\end{theorem}

\begin{proof}[Proof of Theorem \ref{positivity}]
The proof goes along the lines of Theorem 2.3  \cite{TessitoreZabczyk1998}.
Note, due to the linearity with respect to the initial condition and the comparison principle we can start with an initial condition $w_0=\bar w+w_je_j$, where $e_j$ is the $j^{th}$ eigenfunction of $A$ and $\lambda_j$ the corresponding eigenvalue.
Let $\kappa>0$ and $m_0\in\NN$ be chosen such that
$$\frac 12  e^{-\mu\frac t{m_0}}\lk(1-e^{-\lambda_j \frac t{m_0}}\rk)<\kappa.
$$
Let $\rho:={2\gamma +2\frac dq}$ with $\gamma>\frac dq$ and $\rho<1$. Due to the assumption on the Theorem, such a number $\rho$ exists. Let us put $\delta=  e^{-t^{\rho-1 }}$.
Let us define the process
$$
N^{w_0}(t)=\int_0^ t e^{A(t-s)}w^{w_0}(s) \, dW(s).
$$
Fix $L>0$ and let  $w_L$ the unique solution to
\DEQS
dw^{w_0}_L(t)&=& Aw^{w_0}_L(t)+w^{w_0}_L(t)\, (1\wedge |w^{w_0}_L(t)|_{L^q}^{-1})\, dW(t),
\\
w^{w_0}_L(0) &=& {w_0}.
\EEQS
By standard arguments we can show for $w_0\in L^2(\CO)$,  $w^{w_0}_L\in C_b(0,T;L^2(\Omega,L^q(\CO)))$.
Next, let us define the process
$$
N_L^{w_0}(t)=\int_0^ t e^{A(t-s)} w^{w_0}_L(s,x) (1\wedge |w^{w_0}_L(s,x)|_{L^q}^{-1})\,\, dW(s).
$$
Observe, since $ |w^{w_0}_L(s) (1\wedge |w^{w_0}_L(s)|_{L^q}^{-1})|_{L^q}\le L$
we get by \eqref{one-eq} and \eqref{two} 
\DEQS
\lqq{ 
\int_0^ t\|e^{A(t-s)} w^{w_0}_L(s) (1\wedge |w^{w_0}_L(s)|_{L^q}^{-1})\|_{\gamma(L^2,C_b)}^2\, ds} &&
\\
&\le&
\int_0^ t(t-s)^{-2\gamma} \|e^{A(t-s)/2} w^{w_0}_L(s) (1\wedge |w^{w_0}_L(s)|_{L^q}^{-1})\|_{\gamma(L^2,L^q)}^2\, ds
\\
&\le&
\int_0^ t(t-s)^{-2\gamma} s^{-\frac {2d}q} \, | w^{w_0}_L(s,x) (1\wedge |w^{w_0}_L(s,x)|_{L^q}^{-1}) |_{L^q}^2\, ds
\\
&\le & L^2 \int_0^ t (t-s)^{-2\gamma } s^{-\frac {2d}q} \,  ds\le L^2 t^{1-2\gamma-2\frac {d}q}
\EEQS
\newcommand{\wwn}{{w_0}}
Applying the exponential tail estimate \eqref{ete} we have
\DEQS
\PP \lk( \sup_{s\in[0,t]}|N_L^{\wwn}(s)|_{C_b}>l_0\rk) \le  K \exp\lk( - t^{2\gamma+2\frac dq-1}\frac {\lambda l_0^2}{L^2 }\rk).
\EEQS
Similarly, we have
\DEQS
\PP \lk( \sup_{s\in[0,t]}|N_L^{\wwn}(s)|_{L^q}>l_0\rk) \le  K \exp\lk( - t^{2\frac dq-1}\frac {\lambda l_0^2}{L^2 }\rk).
\EEQS
This gives
\DEQS
\lqq{ \PP \lk( \sup_{s\in[0,t]}|N^{\wwn}(s)|_{C_b}>l|w_0|_{L^q} \rk)}
&&
\\ & \le& \PP \lk( \sup_{s\in[0,t]}|N_L^{\wwn}(s)|_{C_b}>l|w_0|_{L^q}\rk)+\PP \lk( \sup_{s\in[0,t]}|w_L^{\wwn}(s)|_{L^q}>l|w_0|_{L^q}\rk)
\\ & \le& \PP \lk( \sup_{s\in[0,t]}|N_L^{\wwn}(s)|_{C_b}>l|w_0|_{L^q}\rk)+\PP \lk( \sup_{s\in[0,t]}|N_L^{\wwn}(s)+e^{A  t}w_0|_{L^q}>l|w_0|_{L^q}\rk)
\\ & \le& \PP \lk( \sup_{s\in[0,t]}|N_L^{\wwn}(s)|_{C_b}>l|w_0|_{L^q}\rk)+\PP \lk( \sup_{s\in[0,t]}|N_L^{\wwn}(s)|_{L^q}>(l+\|e^{tA}\|_{L(L^q,L^q)})|w_0|_{L^q}\rk).
\EEQS
By the choice of
$$
L=\lk( \sup_{s\in[0,t]} \|e^{s A}\|_{L(L^q,L^q)}+l\rk)|w_0|_{L^q}
$$
we get
\DEQS
 \PP \lk( \sup_{0\in[0,t]}|N^{\wwn}(s)|_{C_b}>l|w_0|_{L^q} \rk)
 & \le& K \exp\lk( - t^{2\gamma+2\frac dq-1}\rk).
\EEQS

\medskip
\del{Observe, due to the Neumann conditions, we have for all $x\in\CO$
$$|e^{tA} w_0(x)-\bar w|\le e^{-\lambda_1t}|w(x)-\bar w|,
$$
where $\bar w=\frac 1{\Leb(\CO)}\int_\CO w(x)\, dx$ and $\lambda_1$ is the first nonzero eigenvalue of
the operator $A$.}
Let $m\ge m_0$ and
$$
\mathscr{E}_k:= \lk\{ w^{w_0}\lk(\frac {kt}m,x\rk)\ge \kappa w_{k-1}(x)\,\mbox{ for almost all } x\in\CO\rk\}
$$
and fix $\delta>0$ such that $ e^{-\lk(\frac tm\rk)^{2\gamma+2\frac dq-1 }}\ge \delta$.  For 
$ w_0\in L^q(\CO)$ and $w_{k-1}:=w^{w_0}\lk(\frac {t(k-1)}m\rk)$ we have
$$
\PP\lk( \Omega\setminus \mathscr{E}_k\mid \mathscr{E}_0\cap \mathscr{E}_1\cap\cdots\cap \mathscr{E}_{k-1}\rk)
\le \PP\lk(  w^{ w_{k-1}}\lk(\frac {kt}m\rk)<\kappa   w_{k-1}\rk)
.
$$
It follows that  we have
$$\Omega\setminus \mathscr{E}_1\subset \lk\{ \exists (x,s)\in\lk(\CO,\lk[0\frac tm\rk]\rk) \mbox{ such that } N^{w_0}(s,x)\ge \kappa w_0(x)\rk\},
$$
and, consequently,
\DEQS
\lqq{ \PP\lk( \Omega\setminus \mathscr{E}_1 \rk)} &&
\\
& \le& \PP \lk( \lk\{ \exists (x,s)\in(\CO,[0\frac tm] \mbox{ such that } N^{w_0}(s,x)\ge \kappa w_0(x)\rk\}\rk)
\le e^{-\lk(\frac tm\rk)^{2\gamma+2\frac dq-1 }}.
\EEQS
Similarly, we have
\DEQS
\lqq{ \PP\lk( \Omega\setminus  \mathscr{E}_k\mid \mathscr{E}_0\cap \mathscr{E}_1\cap \ldots \cap \mathscr{E}_{k-1}  \rk)}
&&
\\
& \le& \PP \lk( \lk\{ \exists (x,s)\in(\CO,[0\frac tm] \mbox{ such that } N^{w_{k-1}}(s,x)\ge \kappa w_{k-1}(x)\rk\}\rk)
\\
& \le& e^{-\lk(\frac tm\rk)^{2\gamma+2\frac dq-1 }} =  e^{-\lk(\frac tm\rk)^{\rho-1 }}=\delta^{m^{1-\rho}}
.
\EEQS
Since we have for $m\to\infty$
$$
\PP\lk( \mathscr{E}_m\rk)\ge 1-\sum_{k=1}^{m} \PP\lk(\Omega\setminus \mathscr{E}_k \mid \mathscr{E}_0\cap \mathscr{E}_1\cap \ldots \cap \mathscr{E}_{k-1}  \rk)\ge 1-\sum_{k=1}^m \delta^{m^{1-\rho}}=1-m \delta^{m^{1-\rho}}\rightarrow 1
$$
the positivity is proven by the limit
$$\PP\lk(  \mathscr{E}_m\rk)\longrightarrow \PP\lk(\lk\{ w^{w_0}(t,x)>0\, \, \mbox{ for almost all } x\in\CO\rk\}\rk). $$
\end{proof}

\medskip

\section{A Schauder-Tychonoff theorem}\label{schauder}
In order to prove the existence of a martingale solution to the system \eqref{equ1ss}--\eqref{eqv1ss}, we need a stochastic Schauder-Tychonoff fixed point theorem on the space of progressively measurable processes. We refer to the work of the first named author and T\"olle  \cite{jonas} for the detailed proof of Schauder-Tychonoff fixed point theorem. For the convenience of the reader, we state here the exact Theorem.

\medskip

Let us fix some notation.
Let
$U$ be a Banach space. Let $\Ocal\subset \R^d$ be an open domain, $d\ge 1$. Let $\X\subset\{\eta:[0,T]\to E\subset\Dcal^{\prime}(\Ocal)\}$
be a Banach function space\footnote{Here, $\Dcal^\prime(\Ocal)$ denotes the space of Schwartz distributions on $\Ocal$, that is, the topological dual space of smooth functions with compact support $\Dcal(\Ocal)=C_0^\infty(\Ocal)$.}, let $\X^{\prime}\subset\{\eta:[0,T]\to E\subset\Dcal^{\prime}(\Ocal)\}$
be a reflexive Banach function space embedded compactly into $\X$. In both cases, the trajectories take values in a Banach function space $E$ over the spatial domain $\Ocal$.
Let $\Afrak=(\Omega,\Fcal,\mathbb{F},\P)$ be a filtered probability space
with filtration $\mathbb{F}=(\Fcal_{t})_{t\in [0,T]}$ satisfying the usual conditions.
Let $H$ be a separable Hilbert space and $(W_t)_{t\in [0,T]}$ be a Wiener process\footnote{That is, a $Q$-Wiener process, see e.g. \cite{prato} for this notion.} in $H$ with a linear, nonnegative definite, symmetric trace class covariance
operator $Q:H\to H$ such that $W$ has the representation
$$
W(t)=\sum_{i\in\mathbb{I}} Q^\frac 12 \psi_i\beta_i(t),\quad t\in [0,T],
$$
where $\{\psi_i:i\in \mathbb{I}\}$ is a complete orthonormal system in $H$, $\mathbb{I}$ a suitably chosen countable index set, and $\{\beta_i:i\in\mathbb{I}\}$ a family of independent real-valued standard Brownian motions on $[0,T]$ modeled in $\Afrak=(\Omega,\Fcal,\mathbb{F},\P)$. Due to \cite[Proposition 4.7, p. 85]{prato}, this representation does not pose a restriction.

For $m\ge1$, define the collection of processes
\begin{equation}\label{eq:MMdef}
\begin{aligned}  \Mcal_{\Afrak}^{m}(\X)
:= & \Big\{ \xi:\Omega\times[0,T]\to E\;\colon\;\\
&\qquad\text{\ensuremath{\xi} is \ensuremath{\mathbb{F}}-progressively measurable}\;\text{and}\;\Eb|\xi|_{\X}^{m}<\infty\Big\},
\end{aligned}
\end{equation}
equipped with the semi-norm
\[
|\xi|_{\Mcal_{\Afrak}^{m}(\X)}:=(\Eb|\xi|_{\X}^{m})^{1/m},\quad\xi\in\Mcal_{\Afrak}^{m}(\X).
\]

For fixed $\Afrak$, $W$, $m>1$, we define the operator $\Vcal=\Vcal_{\Afrak,W}:\Mcal_{\Afrak}^{m}(\X)\to\Mcal_{\Afrak}^{m}(\X)$
for $\xi\in\Mcal_{\MA}^m(\X)$ via $\Vcal(\xi):=w$, where $w$ is  the solution of the following It\^o stochastic partial differential equation (SPDE)
\DEQSZ\label{spdes}
dw(t) &=&\lk(\DeltaA w(t)+ F(\xi(t),t)\rk)\, dt +\Sigma(w(t))\,dW(t),\quad w(0)=w_0.
\EEQSZ
Here, we implicitly assume that \eqref{spdes} is well posed and a unique strong solution (in the stochastic sense) $w\in\Mcal_{\MA}^m(\X)$ exists for $\xi\in\Mcal_{\MA}^m(\X)$.
Here, we shall also assume that $\DeltaA:D(\DeltaA)\subset\X\to \X$ is a possibly nonlinear and measurable (single-valued) operator and $F:\X\times [0,T]\to \X$ a (strongly) measurable map such that
\[\PP\left(\int_0^T\left(\left|\DeltaA w (s)\right|_E+\left|F(\xi(s),s)\right|_E\right)\,ds<\infty\right)=1,\]
and assume that $\Sigma:E \to \gamma(H,E)$ is measurable. Here, $\gamma(H,E)$ denotes the space of $\gamma$-radonifying operators from $H$ to $E$,
as defined in the beginning of Section \ref{technics}, which coincides with the space of Hilbert-Schmidt operators $L_{\textup{HS}}(H,E)$ if $E$ is a
separable Hilbert space.

With these two definitions in mind we now state  a stochastic variant of the (deterministic) Schauder-Tychonoff fixed point theorem(s) from \cite[§ 6--7]{granas}.
The proof can be found in \cite{jonas}.

\begin{theorem}\label{ther_main}
Let $H$ be a Hilbert space, $Q:H\to H$ such that $Q$ is linear, symmetric, nonnegative definite and of trace class, let $U$ be a Banach space, and let us assume that we have a compact embedding $\X^{\prime}\hookrightarrow\X$
as above. Let $m>1$. Suppose that
for any filtered probability space $\Afrak=(\Omega,\Fcal,\mathbb{F},\P)$
and for any $Q$-Wiener process $W$ with values in $H$ that is modeled
on $\Afrak$ the following holds.

There exists a nonempty, sequentially weak$^\ast$-closed, measurable and bounded subset\footnote{Here, the notation $\Xcal(\Afrak)$ means that $\text{Law}(u)=\text{Law}(\tilde{u})$
on $\X$ for $u\in\Xcal(\Afrak)$ and $\tilde{u}\in\Mcal_{\tilde{{\Afrak}}}^{m}$
implies $\tilde{u}\in\Xcal(\tilde{\Afrak})$.} $\Xcal(\Afrak)$ of $\Mcal_{\MA}^m(\X)$ such that
the operator $\Vcal_{\Afrak,W}$ restricted to $\Xcal(\Afrak)$
satisfies the following
properties:
\begin{enumerate}
\item $\Vcal_{\Afrak,W}(\Xcal(\Afrak))\subset\Xcal(\Afrak)$,
\item the restriction $\Vcal_{\Afrak,W}\big| 
_{\Xcal(\Afrak)}$ is
continuous w.r.t. the strong topology of $\Mcal_{\Afrak}^{m}(\X)$,
\item there exist constants $R>0$, $m_0\ge m$ such that
\[
\Eb|\Vcal_{\Afrak,W}(v)|_{\X^{\prime}}^{m_0}\le R\quad\text{for every}\quad v\in\Xcal(\Afrak),
\]
\item $\Vcal_{\Afrak,W}(\Xcal(\Afrak))\subset\D([0,T];U)$ $\P$-a.s.\footnote{Here, $\mathbb{D}([0,T);U)$ denotes the
Skorokhod space of c\`adl\`ag paths in $U$.}
\end{enumerate}
Then, there exists a filtered probability space $\tilde{\Afrak}=(\tilde{\Omega},\tilde{\Fcal},\tilde{\mathbb{F}},\tilde{\P})$
(that satisfies the usual
conditions)
together with a $Q$-Wiener process $\tilde{W}$ modeled on $\tilde{\Afrak}$
and an element $\tilde{v}\in\Mcal_{\tilde{\Afrak}}^{m}(\X)$ such that
for all $t\in[0,T],$ $\tilde{\P}$-a.s.
\[
\Vcal_{\tilde{\Afrak},\tilde{W}}(\tilde{v})(t)=\tilde{v}(t).
\]
\end{theorem}

\begin{theorem}\label{schauder}
Fix a filtered probability space $\mathfrak{A}=(\Omega,\CF,\mathbb{F},\PP)$ with filtration
$\mathbb{F}=\{\CF(t):t\in [0,T]\}$ and  let $W$ be a  cylindrical Wiener processes  in $H$ over the probability space
$\mathfrak{A}$.
Let $E$ be a Banach space and
let $\MsV _{\mathfrak{A},W}: \CM^r_{\mathfrak{A}}(E,v_0)\to\CM^r_{\mathfrak{A}}(E,v_0)$ be an operator such that
there exists a bounded subset
$\mathcal{X}(\mathfrak{A})$  of $\CM^r_{\mathfrak{A}}(E)$ such that
\begin{enumerate}[label=(\alph*)]
\item there exist a constant $K>0$ and two measurable mappings $\theta_i:C_b([0,T];E)\to [0,\infty]$, $i=1,2$, such that
$$ \mathcal{X}(\mathfrak{A})=:\Big\{ \xi \in \CM^r_{\mathfrak{A}}(E,v_0) \mid \EE\theta_1(\xi)\le K,\quad \PP(\theta_2(\xi)<\infty)=1\Big\};
$$
  \item the operator $\MsV _{\mathfrak{A},W}$ maps $\mathcal{X}(\mathfrak{A})$ into itself; 
  \item the operator $\MsV _{\mathfrak{A},W}:\mathcal{X}(\mathfrak{A})\to \mathcal{X}(\mathfrak{A})$ is continuous with respect the topology of $ \CM^r_{\mathfrak{A}}(E,v_0)$;
\item  
the following two items are satisfied: 
\begin{enumerate}[label=(\roman*)]
\item
there exist a $\CF_0$--measurable function $g:\Omega\times [0,T]\to E$,  a space $E_0$ such that the embedding $E_0\hookrightarrow E$ is compact,  a number $C>0$ and a function $\phi:[0,T]\to \RR_0^+$ bounded away from zero with
$$\frac 1{\phi(t)}\, \EE |\MsV _{\MA,W} (\xi)(t)-g(t)|_{E_0}\le C,\quad \forall\, t\in[0,T],\, \forall\,\xi\in \mathcal{X}(\MA);
$$
\item  there exist a constant $C>0$ and a function $h:[0,T]\to\RR_0^+$, such that $\lim_{t\downarrow 0} h(t)=0$, and for any $0< t_1<t_2\le T$ and $\xi\in\mathcal{X}(\MA)$ we have
\DEQS
\EE\sup_{0\le t_1<s\le t_2\le T}| \MsV _{\MA,W} (\xi )(t_1)-\MsV _{\MA,W}(\xi )(t_2)|_{E}^2\le C\, h(|t_2-t_1|)  .
\EEQS
 \end{enumerate}
 \end{enumerate}
Then there exists a complete filtered probability space
 $$
 \tilde{\mathfrak{A}}:=(\tilde \Omega ,\tilde {\mathcal{F}},\tilde {\mathbb{F}},\tilde {\mathbb{P}})
 $$
  with the filtration $\tilde {\mathbb{F}}=\{{\tilde{\mathcal{F}}}_t:t\in [0,T]\}$ satisfying the usual conditions, a cylindrical Wiener process  $\tilde W$ in $H$ defined over the probability space $\tilde {\mathfrak{A}}$ and an element $\tilde v\in\CM_{\tilde {\mathfrak{A}}}^r(E,v_0)$ such that we have for all $t\in[0,T]$ $\tilde \PP$--a.s.\
 $$\MsV _{\tilde{\mathfrak{A}},\tilde W}(\tilde {\mathfrak{A}})(\tilde v)(t)=\tilde v(t).
 $$
\end{theorem}

\medskip

\dela{

\section{Some estimates}

\subsection{Function spaces}
Let us introduce firstly some functions spaces.
For a Banach space $E$ and $0\le a<b<\infty$, let $C^{\beta}_b(a,b;E)$ denote a set of all continuous and bounded functions $u:[a,b]\to E$ such that
$$
\| u\|_{C_b^{\beta}(a,b;E)} :=\sup_{a\le t\le b} |u(t)|_E +\sup_{a\le s,t\le b\atop t\not= s} \frac {|u(t)-u(s)|}{|t-s|^\beta}
$$
is finite. The space $C_b^{\beta}(a,b;E)$  endowed with the norm $\| \cdot\|_{C_b^{\beta}(a,b;E)}$ is a Banach space.
\del{The Theorem of Arzel--Ascoli says the following.
A set $K\subset C^0_b(a,b;E)$ is  compact in $ C^0_b(a,b;E)$, if (see \cite[p.\ 71, Lemma 1]{simon})
\begin{enumerate}
  \item for all $t\in[a,b]$ the set $\{ f(t):f\in K\}$ is compact in $E$;
  \item $K$ is uniformly equicontinuous, i.e., for all $\ep>0$ there exists a $\delta>0$ such that $|f(t_1)-f(t_2)|_E\le \ep$ for all $f\in K$, $0\le t_1<t_2\le T$ with $t_2-t_1\le \delta$.
\end{enumerate}
Let $E_1\hookrightarrow E$ compactly and $h:[0,1]\to\RR^+_0$ with $\lim_{t\downarrow 0} h(t)=0$.  It follows that the Banach space
$$\CK:=\{ f\in  C^0_b(a,b;E): \sup_{a\le t \le b} |f(t)|_{E_1}<\infty,\quad \sup_{a\le t_1<t_2\le b} \frac {|f(t_2)-f(t_1)|_E}{h(t_2-t_1)}<\infty\}
$$
equipped with norm
$$
|f|_{\CK}=\sup_{a\le t \le b} |f(t)|_{E_1}+ \sup_{a\le t_1<t_2\le b} \frac {|f(t_2)-f(t_1)|_E}{h(t_2-t_1)},\quad f\in\CK,
$$
is compactly embedded in $C^0_b(a,b;E)$.
}
\newline
As second point in this section we first recall some well known facts  concerning Nemytskii operators, which are necessary to prove our main result.
Most of the content is taken from Runst and Sickel  \cite{runst}.
\del{ For any $\gamma>0$, the completion
of $E$ with respect to the norm $\lvert A^{-\gamma}\cdot \rvert$ will be  denoted by $D(A^{-\gamma})$.
{ For any $\gamma>0$, the domain of the fractional power operator $A^\gamma$ (in $E$) will be  denoted by $D(A^{\gamma})$.
} {With few exceptions we will only speak about fractional powers of the operator $A$ with respect to the space $E$ and not $X$ nor $B$ and hence the notation $A^{\gamma}$ and $D(A^{\gamma})$ should be unambiguous. Those few exceptions are when we use notation $D(A_Y^\theta)$, for instance in Theorem \ref{Th:general}.
 Finally, let us note that since by
 assumption $A^{-1}$  exists and is bounded (on $E$),   the fractional powers
 $A^{-\gamma}$, $\gamma \ge 0,$ are bounded (on $E$) too. }}
For the definition of the  Triebel Lizorkin spaces $F_{p,q}^s(\RR^d)$ and the Besov spaces $B_{p,q}^s(\RR^d)$
we refer the reader to \cite{runst} or \cite{triebel}.

Let us shortly recall some known identities. The proof can be found, for example, in  \cite{triebel}.
\begin{itemize}
     \item $L_p(\RR^d ) = F_{p,2}^0(\RR^d ) $ for $1<p<\infty$,
     \item $W_p^m (\RR^d ) = F_{p,2}^m(\RR^d ) $ for $1<p<\infty$, $m\in\NN$,
     \item $W_p^s (\RR^d ) = F_{p,p}^s(\RR^d )=B_{p,s}^s(\RR^d )  $ for $1<p<\infty$, $s>0$, $s\not\in\NN$,
     \item $H_p^s (\RR^d ) = F_{p,2}^s(\RR^d ) $  for $1<p<\infty$, $s\in\RR$.
      \end{itemize}

\vskip 0.1 in

We list here some results which are useful to treat nonlinearities.
Assume $s_1<0<s_2$.

\begin{theorem}\label{RS2}(see \cite[p.\ 229]{runst})
Assume $s=s_1\le s_2$, $s_1+s_2>d\cdot\max(0,\frac 1p-1)$, and $q\ge \max(q_1,q_2)$. Then
\begin{itemize}
  \item if $s_2>s_1$, then $F_{p,q_1}^{s_1}(\RR^d) \cdot B_{\infty,q_2}^{s_2} (\RR^d)\hookrightarrow F_{p,q_1}^{s_1}(\RR^d)$;
  \item if $s_1=s_2$, then $F_{p,q_1}^{s_1} (\RR^d)\cdot B_{\infty,q_2}^{s_1}(\RR^d) \hookrightarrow F_{p,q}^{s_1}(\RR^d)$.
\end{itemize}
\end{theorem}

\begin{theorem}\label{RS3}(see \cite[p.\ 238]{runst})
Let $s>0$,
$$
\frac 1 {r_1} =\frac 1d \lk( \frac d{p_1} -s\rk) >0, \quad \mbox{ and } \quad
\frac 1 {r_2} =\frac 1d \lk( \frac d{p_2} -s\rk) >0,
$$
and
$$
\frac 1 {r_1}+\frac 1 {r_2} =\frac 1r =\frac 1d \lk( \frac dp -s\rk) <1.
$$
Then
$$
F_{p_1,q_1}^s (\RR^d ) \cdot F_{p_2,q_2}^s (\RR^d ) \hookrightarrow F_{p,q}^s (\RR^d ) ,
$$
iff $$ \max( q_1,q_2)\le q \le \infty.
$$
In addition,
$$
B_{p_1,q_1}^s (\RR^d ) \cdot B_{p_2,q_2}^s (\RR^d ) \hookrightarrow F_{p,q}^s (\RR^d ) ,
$$
iff $$ \max( q_1,q_2)\le q \le \infty, \quad \mbox{ and } \quad 0<q_1\le r_1,\quad  0<q_2\le r_2.
$$
\end{theorem}
{\bf Step B-(a):}
To do so we modify in addition the system. However, before we start let us make some observations:
\begin{remark}\label{functional_remark}
Later on we will see that calculating the Lyapunov function, some additional terms will be bounded. To be more precise,
we will get
\DEQS
\EE \CF(u,v)+\int_0^ t\int (\nabla \sqrt{u})^2\, dx \, ds +\int_0^ t |\Delta v(s)|^2_{L^2}\, ds \le  C(n^\ast).
\EEQS
\end{remark}

\section{Powers and Multiplication}

Let $\sigma_p:= d \max(0,\frac 1p-1)$ and $\sigma_{p,q}:= d \max(0,\frac 1p-1,\frac 1q-1)$, where $d$ is the dimension.

\subsection{Powers}
 \begin{theorem}\label{RS1}(see 
\cite[p.\ 355]{runst})
Suppose $1<p<\infty$, $1\le q\le \infty$, $0<s<1+\frac 1p$ with $s\not=1$.
Then the following estimates holds
$$
\| \,|f|\,|F_{p,q}^s\|\le\| f|F_{p,q}^s\|, \quad f\in F_{p,q}^s.
$$
Suppose $1\le p\le \infty$, $0<q\le \infty$, $0<s<1+\frac 1p$.
Then the following estimates holds
$$
\| \,|f|\,|B_{p,q}^s\|\le\| f|B_{p,q}^s\|, \quad f\in F_{p,q}^s.
$$

\end{theorem}
Remark: It can be extended to the case:
$$
\frac 12<p<1,\quad \max(\sigma_p,\frac 1p)<s<2.
$$

\begin{theorem}\label{RS11}(see 
\cite[p.\ 365]{runst})
Let $0<\mu<1$.
Suppose
$$\sigma_{p,q}<s<1.
$$
Then there exists a constant $C>0$ such that the following estimates holds
$$
\| \,|f|^\mu\,|F^{\frac s\mu}_{\frac p\mu,\frac q\mu}\|\le\| f|F_{p,q}^s\|^\mu, \quad f\in F_{p,q}^s.
$$
Suppose
$$\sigma_{p}<s<1.
$$
Then there exists a constant $C>0$ such that the following estimates holds
$$
\| \,|f|^\mu\,|B^{\frac s\mu}_{\frac p\mu,\frac q\mu}\|\le\| f|B_{p,q}^s\|^\mu, \quad f\in F_{p,q}^s.
$$
\end{theorem}

\begin{corollary}\label{RS9}(see 
\cite[p.\ 365]{runst})
Let $0<\mu<1$ and $\mbox{supp}(f)\subset \CO$.
Suppose
$$\sigma_{p,q}<s<1.
$$
Then there exists a constant $C>0$ such that the following estimates holds
$$
\| \,|f|^\mu\,|F^{\frac s\mu}_{ p,\frac q\mu}\|\le C\| f|F_{p,q}^s\|^\mu, \quad f\in F_{p,q}^s.
$$
Suppose
$$\sigma_{p}<s<1.
$$
Then there exists a constant $C>0$ such that the following estimates holds
$$
\| \,|f|^\mu\,|B^{\frac s\mu}_{ p,\frac q\mu}\|\le C\| f|B_{p,q}^s\|^\mu, \quad f\in F_{p,q}^s.
$$
\end{corollary}

Assume $s_1<0<s_2$

\begin{theorem}\label{RS1}(see 
\cite[p.\ 222]{runst})
$F_{p,q}^s(\RR^d)$ is a multiplication algebra, iff
$F_{p,q}^s(\RR^d)\hookrightarrow L^\infty(\RR^d)$, or
either $s>\frac dp$, or, if $s=\frac dp$ and $0<p\le 1$.
\end{theorem}

\begin{theorem}\label{RS2}(see \cite[p.\ 229]{runst})
Assume $s=s_1\le s_2$, $s_1+s_2>d\cdot\max(0,\frac 1p-1)$, and $q\ge \max(q_1,q_2)$. Then
\begin{itemize}
  \item if $s_2>s_1$, then $F_{p,q_1}^{s_1}(\RR^d) \cdot B_{\infty,q_2}^{s_2} (\RR^d)\hookrightarrow F_{p,q_1}^{s_1}(\RR^d)$;
  \item and, if $s_1=s_2$ then $F_{p,q_1}^{s_1} (\RR^d)\cdot B_{\infty,q_2}^{s_1}(\RR^d) \hookrightarrow F_{p,q}^{s_1}(\RR^d)$
\end{itemize}
\end{theorem}

\begin{theorem}\label{RS2a}(see \cite[p.\ 190]{runst})
Let $s_1\le s_2$ and $s_1+s_2>d\cdot\max(0,\frac 1p-1)$.
\begin{itemize}
  \item Assume $s_2>\frac dp$ and $q\ge \max(q_1,q_2)$. Then, if $s_2>s_1$,
  $F_{p,q_1}^{s_1}(\RR^d) \cdot F_{p,q_2}^{s_2} (\RR^d)\hookrightarrow F_{p,q_1}^{s_1}(\RR^d)$;
  if $s_2=s_1$,
  $F_{p,q_1}^{s_1}(\RR^d) \cdot F_{p,q_2}^{s_2} (\RR^d)\hookrightarrow F_{p,q}^{s_1}(\RR^d)$;
  \item Let $s_1=s_2=\frac dp$ and $q\ge \max(q_1,q_2)$. If $0<p\le 1$, then  $F_{p,q_1}^{s_1}(\RR^d) \cdot F_{p,q_2}^{s_2} (\RR^d)\hookrightarrow F_{p,q}^{s_1}(\RR^d)$.
  \item If $s_2<\frac dp$, then
  $F_{p,q_1}^{s_1}(\RR^d) \cdot F_{p,q_2}^{s_2} (\RR^d)\hookrightarrow F_{p,q}^{s_1+s_2-\frac dp}(\RR^d)$;
\end{itemize}

\end{theorem}

\begin{theorem}\label{RS3}(see \cite[p.\ 238229]{runst})
Let $s>0$,
$$
\frac 1 {r_1} =\frac 1d \lk( \frac d{p_1} -s\rk) >0, \quad \mbox{ and } \quad
\frac 1 {r_2} =\frac 1d \lk( \frac d{p_2} -s\rk) >0,
$$
and
$$
\frac 1 {r_1}+\frac 1 {r_2} =\frac 1r =\frac 1d \lk( \frac dp -s\rk) <1.
$$
Then
$$
F_{p_1,q_1}^s (\RR^d ) \cdot F_{p_2,q_2}^s (\RR^d ) \hookrightarrow F_{p,q}^s (\RR^d ) ,
$$
iff $$ \max( q_1,q_2)\le q \le \infty.
$$
In addition, then
$$
B_{p_1,q_1}^s (\RR^d ) \cdot B_{p_2,q_2}^s (\RR^d ) \hookrightarrow F_{p,q}^s (\RR^d ) ,
$$
iff $$ \max( q_1,q_2)\le q \le \infty, \quad \mbox{ and } \quad 0<q_1\le r_1,\quad  0<q_2\le r_2.
$$
\end{theorem}

\begin{theorem}\label{RS10}(see \cite[Theorem , p.\ 363]{runst})
Let $\alpha>1$. Suppose $0<s<\min\lk(\frac dp,\alpha\rk)$ and $\alpha\lk( \frac dp-s\rk)<d$.
Put
$$
t={d\over s+\alpha \lk(\frac dp-s\rk)}.
$$
There exists a constant $C>0$ such that
$$
\lk\| |f| ^\alpha  \mid F_{t,q}^s \rk\| \le C \lk\| f \mid F_{p,q}^s \rk\|^\alpha , \quad \forall  f\in F_{p,q}^s
,
$$
provided $q>d/(d+s)$, and
$$
\lk\| |f|^\alpha \mid B_{t,q}^s \rk\| \le C \lk\| f \mid B_{p,q}^s \rk\|^\alpha, \quad \forall  f\in B_{p,q}^s
,
$$
provided $q<d/(\tfrac dp-s)$.
\end{theorem}

\begin{theorem}\label{theorp365}(Runst and sickl  p. 365)
Let $0<\mu<1$
\begin{itemize}
\item suppose $\sigma _{p,q} <s<1$.
 Then there exists a constant $c$ such that
$$
\| |f|^\mu \mid F_{\frac p\mu,\frac q\mu}^{s\mu}\| \le c\| f\mid F_{p,q}^s\|^\mu.
 $$

        \item Suppose $\sigma_p<s<1$. Then  there exists a constant $c$ such that
$$
\| |f|^\mu \mid B_{\frac p\mu,\frac q\mu}^{s\mu}\| \le c\| f\mid B_{p,q}^s\|^\mu.
 $$
 \end{itemize}

\end{theorem}

\begin{remark}\label{rem2p31rs} \cite[p.31, Remark 2]{runstsickel}
One get for $s_0-\frac d{p_0}\ge s-\frac dp$ and $p_0\le p$ $F_{p_0,q_0}^ {s_0} \hookrightarrow F_{p,q}^ {s}$.
\end{remark}

\begin{theorem}\label{theors01}
 \cite[p.31]{runstsickel}\label{theors01}
Let $0<p_0<p<p_1$ and suppose
$$
s_0=\frac d{p_0} = s-\frac dp = s_1-\frac d{p_1}.
$$
Then
$$
B_{p_0,u}^{s_0} \hookrightarrow F_{p,q}^ s \hookrightarrow B_{p_1,v}^ {s_1},
$$
if and only if $0<u\le p\le v\le \infty$.
\end{theorem}

}

\del{ \bibliographystyle{plain} 

\bibliography{spdesover3,GM_bibtex,gs_literatur3}

\end{document}
}


\end{document}